\documentclass[11pt]{article}
\usepackage[T1]{fontenc}
\usepackage[utf8]{inputenc}
\usepackage[margin=1in]{geometry}
\usepackage{lmodern,microtype}
\usepackage{amsmath,amssymb,amsthm,mathtools}
\usepackage{enumitem}
\usepackage[colorlinks=true,linkcolor=blue,citecolor=blue,urlcolor=blue]{hyperref}
\hypersetup{bookmarksdepth=2,pdftitle={Strict lax algebra bundles and coherent descent: Matrix embeddings, prescribed ranks, and higher twists},pdfauthor={Andrei V. Ershov},pdfsubject={Finite matrix models for UHF bundles and relative K-theory twists}}
\newtheorem{theorem}{Theorem}[section]
\newtheorem{proposition}[theorem]{Proposition}
\newtheorem{lemma}[theorem]{Lemma}
\newtheorem*{unnumberedlemma}{Lemma}
\newtheorem{corollary}[theorem]{Corollary}
\theoremstyle{definition}
\newtheorem{definition}[theorem]{Definition}
\theoremstyle{remark}
\newtheorem{remark}[theorem]{Remark}
\newtheorem{example}[theorem]{Example}
\DeclareMathOperator{\Aut}{Aut}
\DeclareMathOperator{\Emb}{Emb}
\DeclareMathOperator{\Iso}{Iso}
\DeclareMathOperator{\Ad}{Ad}
\DeclareMathOperator{\Sing}{Sing}
\DeclareMathOperator{\diag}{diag}
\DeclareMathOperator{\Ex}{Ex}
\DeclareMathOperator{\Sd}{Sd}
\DeclareMathOperator{\id}{id}
\DeclareMathOperator{\hofib}{hofib}
\DeclareMathOperator{\Map}{Map}
\DeclareMathOperator{\Desc}{Desc}
\DeclareMathOperator*{\hocolim}{hocolim}
\DeclareMathOperator*{\holim}{holim}
\newcommand{\C}{\mathcal C}
\newcommand{\K}{\mathcal K}
\newcommand{\B}{\mathrm B}

\newcommand{\op}{\mathrm{op}}

\newcommand{\HCgpd}{\mathsf{HC}^{\mathrm{gpd}}_{k,l}}
\newcommand{\HCraw}{\mathsf{HC}^{\mathrm{raw}}_{k,l}}

\title{Strict lax algebra bundles and coherent descent\\
Matrix embeddings, prescribed ranks, and higher twists}
\author{Andrei V. Ershov}
\date{2 December 2026}
\begin{document}
\maketitle
\vspace{-1.5em}
\begin{abstract}
For coprime $m\ge1$ and $l\ge2$, the category of matrix algebras
$M_{ml^r}$ and unital embeddings has classifying space
$B\Aut(M_m\otimes M_{l^\infty})$. On finite CW complexes it
classifies strict lax algebra bundles modulo morphisms, refinement
and concordance, and coherent cover diagrams up to the specified
refinements and concordances. We identify the same space with the
homotopy quotient of the rank-$m$ component of $\Omega^\infty KU[1/l]$
by the identity component of its multiplicative units. On a connected finite
CW complex, an orientable $KU[1/l]$-twist with a chosen rank coordinate
admits a LAB of size $m$ exactly when its twisted $K^0$-group contains
a class of fiber rank $m$.

For $k>1$ coprime to $l$, tensoring through all sizes $k^n$ gives
a spectrum whose connected cover classifies an orientable twist with
a trivialization after inverting $k$. This cover is the fiber of the map
between the simply connected covers of the Picard spectra induced by
$KU[1/l]\to KU[1/(kl)]$. Forgetting the
trivialization gives a short exact sequence onto the $k$-power
torsion in the Dadarlat--Pennig Brauer group. The sequence splits
on each connected finite CW complex. A Moore space and its bottom
cell obstruct a natural choice of splitting. Complex linearization of finite sets
commutes with forgetting the prescribed section. On $S^2$ its image in the
matrix LAB group $\mathbb Z/m$ is exactly the subgroup killed by $2$.
For an odd prime $p\nmid l$, a nonzero class of size $p$ on $S^{2p-2}$
has a permutation representative using sets of at most $pl^N$ elements
exactly when $l^N\ge p$.

\end{abstract}
\begingroup
\small
\makeatletter
\renewcommand*\l@section[2]{\begingroup\bfseries\@dottedtocline{1}{0em}{1.8em}{#1}{#2}\endgroup}
\makeatother
\setcounter{tocdepth}{1}
\tableofcontents
\endgroup
\clearpage

\section{Introduction}\label{sec:intro}

Locally trivial bundles of matrix algebras provide a geometric form of the
classical topological Brauer group. For a finite CW complex $X$, their
Morita-equivalence classes are classified by the torsion subgroup of
$H^3(X;\mathbb Z)$. The corresponding Dixmier--Douady invariant is also
the source of the familiar degree-three twists of complex $K$-theory
\cite{AS}. From the viewpoint of ring spectra, twists are organized by
the space $BGL_1(KU)$ of $KU$-line bundles, which contains information
beyond the degree-three invariant \cite{ABGHR}. This broader setting
motivates the search for geometric objects built from finite-dimensional
algebras that retain more homotopical information than Morita equivalence.

The starting point here is an embedding problem. Given an $M_k(\mathbb C)$-bundle
$B\to X$, when can one embed it unitally into a trivial matrix algebra
bundle? The size of the target matters. If the allowed targets are
$M_{kl^n}$ with $(k,l)=1$, obstructions can persist as $n$ increases.
The matrix Grassmannians and their tautological bundles make this a
concrete lifting problem \cite{ErshovEmb}. This contrasts with the
familiar embedding of a vector bundle over a compact base into a large
trivial vector bundle: requiring a unital algebra embedding imposes
multiplicative and arithmetic constraints. Section~\ref{sec:prelim}
recalls the geometric constructions needed to formulate these constraints.

The obstruction already occurs on $S^2$. Let $H$ be the Hopf line bundle,
let $k>1$, and form the $M_k$-bundle
$B=\operatorname{End}(H\oplus\underline{\mathbb C}^{\,k-1})$,
where the underlined summand is trivial. As an endomorphism algebra
bundle, $B$ has zero classical Brauer class. Its clutching loop generates
$\pi_1\mathrm{PU}(k)\cong\mathbb Z/k$, since the determinant of the
clutching loop of $H\oplus\underline{\mathbb C}^{\,k-1}$ has winding
number $\pm1$. For every $l\geq2$ coprime to $k$, the stabilized
bundle $B\otimes M_{l^\infty}$ represents a generator of
$\pi_2 B\Aut(M_k\otimes M_{l^\infty})$
(Example~\ref{ex:morita-versus-lab}). By
Corollary~\ref{cor:stable-embeddability}, $B$ therefore admits no unital
embedding into $S^2\times M_{kl^N}$ for any $N\geq0$. The same clutching
loop represents a symmetry of the trivial LAB on $S^1$ in the completed
theory (Section~\ref{subsec:symmetries}). Section~\ref{sec:perm-linearization}
asks whether this class can be built from permutations. It can when $k=2$.
For every $k>2$ the Hopf-line generator lies outside the permutation image.

A lax algebra bundle, or LAB, allows finite-dimensional algebra data to
be assembled locally even when one global matrix algebra bundle is
insufficient. On a cover $\{U_i\}$, it assigns a matrix algebra bundle
$A_S$ to every nonempty finite intersection $U_S$, together with unital
embeddings $A_S|_{U_T}\to A_T$ when $S\subseteq T$. In the strict version
these embeddings compose exactly. Thus the word \emph{lax} refers to using
embeddings of possibly different matrix sizes in place of transition
isomorphisms; associativity in a strict LAB is not weakened. Coherent
variants require a further choice of model. We distinguish raw cover
diagrams in a Segal space from descent data valued in its groupoid
completion. The latter formally inverts all matrix embeddings and retains
the homotopy type of their classifying space. We compare the concordance
classes of both constructions.

We fix $k\geq1$ and $l\geq2$ with $(k,l)=1$, and write
$D=M_k\otimes M_{l^\infty}$. The finite-dimensional category underlying
the construction is
\[
 \operatorname{Ob}\mathcal C_{k,l}=\{M_{kl^n}:n\geq0\},\qquad
 \mathcal C_{k,l}(n,m)=\operatorname{Emb}(M_{kl^n},M_{kl^m}).
\]
The morphism spaces carry their usual point-norm topology. For $k>1$ this
category and its nerve appeared in \cite[Sections~3.1--3.2]{ErshovCategory},
which also constructed a comparison with a monoid of stabilized
endomorphisms. Here the comparison target is the full point-norm
automorphism group of $D$, and the comparison identifies the class of a
LAB with that of its rectified bundle. We consider strict LABs modulo compatible morphisms,
refinement, and concordance, and denote the resulting functor by
$\mathsf F$. The normalization requiring singleton algebras to have rank
$k$ is included in this definition. We write $\HCgpd$ for the completed
coherent theory, defined by homotopy-coherent \v{C}ech descent in
Section~\ref{sec:coherent}. A distinct symbol $\HCraw$ is reserved for
the raw Segal cover-data theory with specified ordinary refinements.

\paragraph{Main results.}
The first result identifies the homotopy type of the finite-dimensional
category:
\[
 \rho:|N\mathcal C_{k,l}|\xrightarrow{\simeq}B\operatorname{Aut}(D).
\]
The automorphism group has the point-norm topology. Theorem~\ref{thm:main}
also identifies the ordinary realization of
$\operatorname{Ex}(N\mathcal C_{k,l})$ with this homotopy type.
This establishes the nerve comparison without assuming a classification
theorem for LABs.

The second result is the classification theorem. For every finite CW
complex $X$, there are natural bijections
\[
 \mathsf F(X)\ \cong\ \HCraw(X)\ \cong\ \HCgpd(X)\
 \cong\ [X,|N\mathcal C_{k,l}|]\
 \cong\ [X,B\operatorname{Aut}(D)].
\]
The last set classifies locally trivial $D$-bundles. The comparison
agrees with categorical rectification of a strict LAB
(Theorem~\ref{thm:intro-classification}).
In particular, every coherent class on a finite complex has a normalized
strict representative, unique up to the stated strict equivalence relation.
For raw cover diagrams the bijection uses paths, specified index-map
refinements, and concordances over $X\times I$
(Definition~\ref{def:raw-coherent}). At a fixed cover, different
restrictions can still lie in different components of the mapping space.

A third description explains what the matrix size measures. Put
$R=KU[1/l]$. Its space of multiplicative units acts on the
underlying space $\Omega^\infty R$. The component of the unit,
$U_R=GL_1(R)_1$, preserves each rank component $Y_m$.
Theorem~\ref{thm:rank-section} gives a weak equivalence
\[
        |N\C_{m,l}|\simeq Y_m//U_R,\qquad
                     Y_m=(\Omega^\infty R)_m .
\]
Thus a LAB of type $(m,l)$ carries an orientable twist and a
section of rank $m$. For a specified orientable twist with its rank
coordinate on a connected finite CW complex, the possible sizes can
be read from the evaluation map in twisted
$K^0$ (Corollary~\ref{cor:prescribed-rank}).

Tensor product takes sizes $m$ and $a$ to size $ma$. For $k>1$,
combine all sizes $k^n$ and take the group completion of the resulting
symmetric monoidal category. The spectrum of its identity component is
\[
 \operatorname{fib}\!\left(
 \tau_{\geq2}\mathrm{pic}(KU[1/l])
   \longrightarrow \tau_{\geq2}\mathrm{pic}(KU[1/(kl)])\right)
\]
by Theorem~\ref{thm:relative-spectrum}. Here $\mathrm{pic}$ records
invertible module spectra. After extending coefficients, dividing a
section of rank $k^n$ by $k^n$ makes it an invertible section of
rank one. Geometrically this gives a trivialization of the extended
twist. For the $S^2$ example with $k=2,l=3$, the twist is trivial,
while the chosen trivialization represents $1/2$ modulo
$\mathbb Z[1/3]$. A five-dimensional Moore space exhibits a nonzero
higher Brauer class and shows that the forgetful map has no
natural section (Examples~\ref{ex:relative-half} and~\ref{ex:relative-moore}
and Corollary~\ref{cor:no-natural-splitting}).

Replacing a finite set by its complex span gives a functor from the
permutation category $\mathcal P_{k,l}$ to $\C_{k,l}$. On prescribed
sections it extends coefficients along $\mathbb S[1/l]\to KU[1/l]$,
both on the line and on its section (Theorem~\ref{thm:perm-linearization}).
At size four this functor detects a twist on $\Sigma\mathbb{RP}^2$
through the sign line on its equator and kills the spherical twist on $S^3$ clutched
by the square of the stable Hopf element
(Examples~\ref{ex:perm-sign} and~\ref{ex:perm-eta}).
For an odd prime $p\nmid l$, the first even sphere with a nonzero
permutation image at size $p$ is $S^{2p-2}$. The first multiplicity
$l^N$ that can carry a $p$-cycle already suffices, and every smaller
maximum object size gives zero (Theorem~\ref{thm:perm-first-size}). For $p=3$
and $l=2$, Proposition~\ref{prop:perm-twelve-cycle} realizes a generator
on $S^4$ using only sets of sizes three and twelve.

\paragraph{Proof strategy.}
Adjoin the limit algebra $D$ to the matrix category, allowing all
embeddings from the finite stages into $D$ and all automorphisms at $D$.
The two full subcategories, the finite stages and the one-object
automorphism group, each induce a weak equivalence on classifying spaces.
For the first inclusion, approximation of compact families reduces the
only infinite-stage comma category to a telescope of categories with
terminal objects. For the second, the stabilizer of an embedding
$M_{kl^n}\to D$ is the contractible group $\Aut(M_{l^\infty})$.
These are the two applications of topological Quillen Theorem~A.

The strict classification has a geometric proof. Compatible stabilized
sections realize any strict diagram, after shrinking its cover, as
nested subalgebras in a genuine $D$-bundle. A relative simultaneous
gauge lemma places a finite nested family into a common matrix stage.
It gives existence of subalgebra diagrams in every $D$-bundle and
uniqueness up to concordance. Finally, normalization by rank-$k$
singleton cores is checked against all generating relations, including
refinements that identify indices. For completed coherent
classification, the cosimplicial space of local maps has totalization
equal to the mapping space out of the diagonal of the \v{C}ech nerve
after applying the singular functor degreewise. The \v{C}ech augmentation
theorem of Dugger and Isaksen identifies this diagonal with $\Sing X$ up
to weak equivalence for every indexed open cover \cite[Theorems~2.1 and~3.1]{DI}.
For raw diagrams, refine two covers by the vertex stars of a common
triangulation. A homotopy between the completed classes can then be
represented on a subdivision of the prism. A comparison of the carriers
of its endpoint simplices with the specified refinement maps supplies
the required raw concordances. Berwick-Evans, Boavida de Brito and
Pavlov prove a related classification using an enriched category of
covers on smooth manifolds \cite[Section~6]{BEBP}.
The tail automorphism group $\Aut(M_{l^\infty})$ is contractible
\cite{DP}. For the fixed-size description, multiplicative operator
$K$-theory sends a matrix algebra to an $R$-line with a section
given by its unit. Evaluation on endomorphism bundles of vector
bundles over even spheres tests this map on all homotopy groups.
Under tensor stabilization the fixed-size models have the same base
$B U_R$. Their fibers form a telescope whose maps multiply the section
by $k$. This identifies the homotopy fiber and proves the spectrum
comparison using the symmetric tensor products of the scalar sections.

\paragraph{Relation to embedded and global bundles.}
The matrix Grassmannian viewpoint leads to the homotopy-fiber sequence
\[
 \operatorname{Gr}_{k,l^\infty}\longrightarrow B\mathrm{PU}(k)
 \xrightarrow{\,B(\alpha\mapsto\alpha\otimes\id)\,}
 B\operatorname{Aut}(D),
\]
proved in Proposition~\ref{prop:grassmannian-fiber}. Consequently an
$M_k$-bundle embeds into a trivial $M_{kl^N}$-bundle for some $N$ if and
only if its stabilization by $M_{l^\infty}$ is trivial.

The classifying space is simply connected, with $\mathbb Z/k$ in every
positive even homotopy degree and zero in every odd degree. Nevertheless,
the image of fixed-rank global matrix bundles is generally smaller than
the full theory, even on spheres. There are also classes on finite
complexes that do not arise from global $M_{kl^r}$-bundles for any $r$.
Section~\ref{sec:frames} proves these assertions and gives explicit
examples. For fixed $k$, the classifying spaces associated to any two
auxiliary integers $l,l'\geq2$ coprime to $k$ are weakly equivalent
by stabilization.

\paragraph{Scope and organization.}
Dadarlat and Pennig relate torsion in their generalized
Dixmier--Douady theory to stabilized algebra bundles
\cite{DPBrauer}. For the UHF algebra $A=M_{l^\infty}$ their Brauer
group allows fibers $M_n\otimes A$ with varying $n$. At fixed $k$
our classification retains the isomorphism class of the unital
$M_k\otimes A$-bundle. On passing to all powers of $k$, forgetting the
localized trivialization maps onto the subgroup of Brauer classes
annihilated by powers of $k$. Theorem~\ref{cor:brauer-forget}
computes its divisible kernel. This relative description complements
the computations of stable twist groups in ordinary cohomology and
connective $K$-theory by Dadarlat, McClure and Pennig \cite{DMP}.
Rigid equivalence without concordance and the strict theory on
arbitrary infinite CW complexes remain separate questions.

After the preliminaries, Section~\ref{sec:conventions} fixes the
homotopy-theoretic conventions and states the results precisely.
Sections~\ref{sec:approx}--\ref{sec:Ex} prove the nerve equivalence.
Sections~\ref{sec:strict}--\ref{sec:normalization} establish the strict
classification. Section~\ref{sec:coherent} proves completed coherent
descent and identifies the concordance classes of raw cover diagrams.
Sections~\ref{sec:frames} and~\ref{sec:limits} return to embedded bundles,
homotopy groups, examples, and graded tensor products.
Section~\ref{sec:group-completion} constructs the tensor group
completion. Section~\ref{sec:relative-picard} constructs the rank
section, proves the prescribed-size criterion, and identifies the
group-completion spectrum. Section~\ref{sec:brauer-relative} compares
it with orientable algebra twists and proves the splitting obstruction.
Section~\ref{sec:perm-linearization} compares the permutation and matrix
nerves over their twist spaces and determines the least object size on
the first odd-primary sphere.

\smallskip

\begin{samepage}
{\noindent \bf Acknowledgments.\;}
	The author is deeply appreciated to Professor Doctor Thomas Schick for hospitality and 
	very illuminating discussions on a number of related topics. The author would also like to express his deep
	gratitude to Professors V.M. Manuilov, A.S. Mishchenko and E.V. Troitsky for long-time support
	and numerous very helpful discussions. 
	I am also grateful to Professor I.I.~Bogdanov and to my former student Michael Dusman for helpful discussions at an early stage of this work.
\par\end{samepage}

\paragraph{Use of generative AI.}
OpenAI's \href{https://chatgpt.com}{ChatGPT} was used during the development and revision of this
work, including the November and December 2026 revisions, to examine mathematical
arguments, identify possible gaps and counterexamples, develop and
refine proof strategies, assist with literature searches, and draft and
revise portions of the exposition and \LaTeX{} source. Responsibility
for mathematical verification, the references, and the final text
remains with the author.
% Before journal submission, supplement this declaration with the actual
% model/version identifier and further usage dates if the journal requires
% them. No unrecorded model/version or completed independent human review
% has been assumed in the statement above.

\section{Preliminaries: matrix algebra bundles and embeddings}\label{sec:prelim}

This section recalls the finite-dimensional geometry behind the theory.
Unless otherwise specified, base spaces are finite CW complexes.

\subsection{Matrix algebra bundles and their centralizers}

An $M_r$-bundle is a locally trivial bundle of unital complex
$C^*$-algebras with fiber $M_r(\mathbb C)$. Its structure group is
$\operatorname{Aut}_*(M_r)=\mathrm{PU}(r)$, so its isomorphism class is
classified by a map to $B\mathrm{PU}(r)$. The larger group
$\mathrm{PGL}_r(\mathbb C)$ consists of complex-algebra automorphisms;
it has the same homotopy type as $\mathrm{PU}(r)$, but is not the group
of $*$-automorphisms. All bundle maps used here preserve the unit and
the involution.

A unital embedding $M_r\to M_t$ exists precisely when $t=rs$ for an
integer $s$. Every such embedding is unitarily conjugate to
$a\mapsto a\otimes1_s$. Its centralizer is a copy of $M_s$, and
multiplication gives the tensor decomposition $M_t\cong M_r\otimes M_s$.
These facts have a bundle version.

\begin{proposition}\label{prop:finite-centralizer}
Let $\mu:B\to E$ be a unital embedding of an $M_r$-bundle into an
$M_{rs}$-bundle. The centralizers form a locally trivial $M_s$-bundle
$Z(\mu)$, and multiplication is a bundle isomorphism
\[
 B\otimes Z(\mu)\xrightarrow{\cong}E,\qquad
 b\otimes z\longmapsto\mu(b)z.
\]
\end{proposition}
\begin{proof}
Trivialize the two bundles locally. The embedding becomes a map into
the homogeneous space of unital embeddings $M_r\to M_{rs}$.
The quotient map from $\mathrm{PU}(rs)$ has local sections, so a
continuous local unitary conjugacy reduces the pair to the standard
tensor-factor inclusion. Its commutant is $1_r\otimes M_s$.
The resulting local commutants agree on overlaps, and the displayed
multiplication is a fiberwise isomorphism in these local charts.
\end{proof}

An \emph{embedded bundle} includes a specified embedding
$\mu:B\hookrightarrow X\times M_{rs}$. The word \emph{embeddable}
asserts only existence of such an embedding. The target size, or the
allowed sequence of target sizes, is part of either specification.

\subsection{Frames, matrix Grassmannians, and tautological bundles}

Write
\[
 \operatorname{Fr}_{r,s}=\operatorname{Emb}(M_r,M_{rs})
 \cong\mathrm{PU}(rs)/(1_r\otimes\mathrm{PU}(s)).
\]
The group $\mathrm{PU}(r)$ acts freely by precomposition.
The quotient forgets the framing and remembers only the image subalgebra:
\[
 \operatorname{Gr}_{r,s}
 =\operatorname{Fr}_{r,s}/\mathrm{PU}(r)
 \cong\mathrm{PU}(rs)/
       \bigl(\mathrm{PU}(r)\otimes\mathrm{PU}(s)\bigr).
\]
Both spaces are compact smooth homogeneous spaces. The projection
$\operatorname{Fr}_{r,s}\to\operatorname{Gr}_{r,s}$ is a principal
$\mathrm{PU}(r)$-bundle. The associated algebra bundle is the tautological
subbundle
\[
 \mathcal A_{r,s}
 =\{(B,b):B\in\operatorname{Gr}_{r,s},\ b\in B\}
 \ \subset\ \operatorname{Gr}_{r,s}\times M_{rs}.
\]
Its centralizer is a tautological $M_s$-bundle $\mathcal Z_{r,s}$, and
Proposition~\ref{prop:finite-centralizer} gives
$\mathcal A_{r,s}\otimes\mathcal Z_{r,s}
\cong\operatorname{Gr}_{r,s}\times M_{rs}$.

\begin{proposition}[Classification with a chosen embedding]\label{prop:embedded-classification}
The set $[X,\operatorname{Gr}_{r,s}]$ classifies pairs $(B,\mu)$
consisting of an $M_r$-bundle and a unital embedding
$\mu:B\hookrightarrow X\times M_{rs}$ with the target held fixed.
Two pairs $(B_0,\mu_0)$ and $(B_1,\mu_1)$ are equivalent if there is a
bundle isomorphism $h:B_0\to B_1$ over $X$ such that $\mu_0$ and
$\mu_1\circ h$ are homotopic through unital bundle embeddings
$B_0\hookrightarrow X\times M_{rs}$.
The forgetful map to $[X,B\mathrm{PU}(r)]$ is induced by the classifying
map $\gamma_{r,s}$ of $\mathcal A_{r,s}$.
\end{proposition}
\begin{proof}
An embedded bundle gives the continuous map
$x\mapsto\mu(B_x)\in\operatorname{Gr}_{r,s}$ and is isomorphic, as an
embedded bundle, to the pullback of $\mathcal A_{r,s}$.
Conversely that pullback constructs the embedded bundle from any map.
A homotopy of maps gives a family over $X\times I$. Bundle homotopy
invariance identifies its underlying algebra bundle with a pullback
from $X$, and hence identifies such families with homotopies of
embeddings after an isomorphism of the endpoint bundles.
\end{proof}

\paragraph{The role of the fixed target.}
Consider the constant subalgebra $M_2\otimes 1_2\subset M_4$ over $S^1$.
A based loop $g:S^1\to\mathrm{PU}(4)$ acts on the trivial target bundle
and carries this subalgebra to the family $g(x)(M_2\otimes 1_2)$.
Its class in $\pi_1\operatorname{Gr}_{2,2}$ is the image of $[g]$
under the orbit map. Indeed, the fibration
\[
 \mathrm{PU}(2)\times\mathrm{PU}(2)\longrightarrow
 \mathrm{PU}(4)\longrightarrow\operatorname{Gr}_{2,2}
\]
induces $(a,b)\mapsto 2a+2b$ from $(\mathbb Z/2)^2$ to $\mathbb Z/4$
on fundamental groups, since each tensor factor occurs with multiplicity
two. The stabilizer is connected, so $\pi_1\operatorname{Gr}_{2,2}$
is the cokernel $\mathbb Z/2$. A generator loop in $\mathrm{PU}(4)$
therefore carries the constant family to a distinct class in
$[S^1,\operatorname{Gr}_{2,2}]$. Allowing automorphisms of the target
bundle would identify these two classes.

\begin{proposition}[The embedding lifting problem]\label{prop:embedding-lift}
There is a homotopy-fiber sequence
\[
 \operatorname{Fr}_{r,s}\longrightarrow\operatorname{Gr}_{r,s}
 \xrightarrow{\gamma_{r,s}}B\mathrm{PU}(r).
\]
If $f_B:X\to B\mathrm{PU}(r)$ classifies $B$, the bundle $B$ embeds
in $X\times M_{rs}$ exactly when $f_B$ has a homotopy lift through
$\gamma_{r,s}$. For the fixed bundle $B$, homotopy classes of embeddings
are the homotopy classes of sections of its associated embedding bundle.
\end{proposition}
\begin{proof}
Put $P=\mathrm{PU}(r)$. The associated bundle
$EP\times_P\operatorname{Fr}_{r,s}\to BP$ has fiber
$\operatorname{Fr}_{r,s}$. Since the action on frames is free, its total
space is weakly equivalent to the quotient
$\operatorname{Gr}_{r,s}$, over $BP$. Pulling this bundle back by $f_B$
gives the bundle with fiber
$\operatorname{Emb}(B_x,M_{rs})$ at $x$. A section is exactly a
continuous bundle embedding. The homotopy lifting interpretation follows.
\end{proof}

These constructions are the finite-dimensional starting point of
\cite{ErshovEmb}. In the fixed type $(k,l)$, put
\[
 \operatorname{Gr}_{k,l^\infty}
 :=\hocolim_{N\geq0}\operatorname{Gr}_{k,l^N},
\]
with stabilization induced by $b\mapsto b\otimes1_l$. The finite-stage
maps are closed smooth embeddings, so the telescope and ordinary
direct-limit models are weakly equivalent. For a finite CW source,
maps into the telescope can be moved to one finite stage. Consequently
this space classifies embeddings after passage to a larger allowed
matrix stage. Proposition~\ref{prop:grassmannian-fiber} will identify it
as the homotopy fiber of $B\mathrm{PU}(k)\to B\operatorname{Aut}(D)$.

\subsection{The classical Brauer invariant and the distinction of equivalences}

For an $M_r$-bundle $B$, the central extension
$U(1)\to U(r)\to\mathrm{PU}(r)$ gives the Dixmier--Douady class
$\delta(B)\in H^3(X;\mathbb Z)$. It vanishes precisely when
$B\cong\operatorname{End}(V)$ for a rank-$r$ complex vector bundle $V$.
It satisfies $r\delta(B)=0$ and is additive under tensor product.
For a finite CW complex the classical topological Brauer group is
\[
 \operatorname{Br}(X)\cong\operatorname{Tor}H^3(X;\mathbb Z),
\]
where matrix algebra bundles are taken modulo Morita equivalence
\cite{AS}. A torsion class is representable by some matrix algebra
bundle; its order alone does not guarantee a representative of that
same matrix size.

Morita triviality, embeddability in the prescribed coprime stages, and
triviality in the LAB concordance theory are different notions.
The last two will be equivalent for the class coming from a global
$M_k$-bundle, by Corollary~\ref{cor:stable-embeddability}.
Morita triviality alone does not imply them, as
Example~\ref{ex:morita-versus-lab} shows.
In particular, the classification developed here does not identify
the fixed-type LAB functor with the classical Brauer group.

\section{Conventions and the theorems}\label{sec:conventions}

Fix integers $k\geq1$ and $l\geq2$ with $\gcd(k,l)=1$. Set
\[
 A=M_{l^\infty},\qquad D_n=M_{kl^n},\qquad
 D=\varinjlim_n D_n=M_k\otimes A,
\]
where the limit is taken in unital $C^*$-algebras with connecting maps
$a\mapsto a\otimes1_l$. Both $A$ and $D$ are UHF algebras.
Arrows in the matrix categories and the structure maps of LABs are unital
$*$-homomorphisms. We also use nonunital full-corner embeddings when applying
operator $K$-theory. Homomorphism spaces carry the point-norm topology.
Unitary groups used for perturbations carry the norm topology. Products and realizations are taken in compactly generated spaces.

The topological category $\C$ has discrete object set $\mathbb N_0$ and
\[
 \C(n,m)=
 \begin{cases}
  \Emb(D_n,D_m),&n\leq m,\\
  \varnothing,&n>m.
 \end{cases}
\]
Composition is composition of embeddings. For $n\leq m$,
\[
 \C(n,m)\cong
 \mathrm{PU}(kl^m)/(1_{kl^n}\otimes\mathrm{PU}(l^{m-n})),
\]
a compact smooth manifold; in particular $\C(n,n)=\mathrm{PU}(kl^n)$.
Write $G=\Aut(D)$ and $H=\Aut(A)$.

For a topological category $\mathcal D$, write
$\B\mathcal D=\|N\mathcal D\|$ for the fat realization of its nerve.
This realization glues along faces. The ordinary realization $|X_\bullet|$
also makes the degeneracy identifications. Fat realization takes
levelwise weak equivalences to weak equivalences
\cite[Theorem~2.2]{ERW}. For a proper simplicial space the quotient
$\|X_\bullet\|\to |X_\bullet|$ is a homotopy equivalence
\cite[Appendix~A]{Segal74}. For a topological group $G$ we write
$BG$ for the classifying homotopy type of numerable principal bundles
and use the fat nerve of the one-object category $*G$ as its model
\cite[Section~2]{GH}.
Throughout, $\simeq$ denotes an equivalence of weak homotopy types.

\begin{theorem}\label{thm:main}
There is a map in the homotopy category
\[
 \rho:|N\C|\longrightarrow B\Aut(D)
\]
which is a weak equivalence. It is represented by the categorical zigzag
constructed in Section~\ref{sec:K}. Moreover, for
$T=\Ex(N\C)$, formed in the simplicial direction,
\[
 |T|\simeq |N\C|\simeq B\Aut(D),
\]
where the two realizations on the left are ordinary geometric realizations.
\end{theorem}

\begin{theorem}[Classification on finite complexes]\label{thm:intro-classification}
For finite CW complexes $X$, with the definitions and extension from
polyhedra in Sections~\ref{sec:strict}--\ref{sec:coherent}, there are
natural bijections
\[
 \mathsf F(X)\cong\HCraw(X)\xrightarrow[\cong]{\chi_X}\HCgpd(X)
 \cong[X,|N\C|]\cong[X,B\Aut(D)].
\]
Here $\mathsf F$ uses normalized strict diagrams and the equivalence
relation generated by morphisms, refinement, and concordance.
The functor $\HCraw$ uses the specified refinements and concordances of
Definition~\ref{def:raw-coherent}, and $\HCgpd$ uses coherent descent
after groupoid completion (Definition~\ref{def:coherent}).
The first bijection takes a strict diagram to its raw class, and the
composite from $\mathsf F(X)$ to $[X,B\Aut(D)]$ is categorical rectification.
\end{theorem}

Theorem~\ref{thm:main} is proved independently in Sections~\ref{sec:approx}--\ref{sec:Ex}.
The strict classification is proved in Sections~\ref{sec:strict}--\ref{sec:normalization}.
Theorem~\ref{thm:coherent}, Corollary~\ref{cor:strict-coherent}, and
Proposition~\ref{prop:raw-comparison} establish the coherent comparisons.

We use the following external operator-algebraic theorem:
\begin{equation}\label{eq:H}
 H=\Aut(M_{l^\infty})\text{ is contractible in the point-norm topology.}
\end{equation}
It follows from Dadarlat--Pennig~\cite[Theorem~2.3]{DP}, since $A$ is an
infinite-type UHF algebra. Strong self-absorption of $D$ is not a hypothesis.
The finite-stage approximation needed below is proved in the next section.

For a sequence of spaces we use its ordinary mapping telescope as a
model of its homotopy colimit. For closed cofibrations the map from this
telescope to the direct limit is a weak equivalence. The simplicial
mapping spaces used for coherent descent are introduced in
Section~\ref{sec:coherent}.

\section{Relative finite-stage approximation}\label{sec:approx}

\begin{lemma}[Close embeddings]\label{lem:close}
Let $B$ be a unital $C^*$-algebra and let $\varphi,\psi:M_r\to B$ be
embeddings sufficiently close on the standard matrix units. There is a
unitary $u(\psi,\varphi)$, continuous in the two embeddings and close to $1$,
such that $\psi=\Ad(u)\varphi$. It equals $1$ when $\psi=\varphi$.
\end{lemma}
\begin{proof}
Set
\[
 s=\sum_{i=1}^{r}\psi(e_{i1})\varphi(e_{1i}).
\]
Then $s\varphi(e_{ab})=\psi(e_{ab})s$, and $s$ is close to $1$, hence
invertible. The polar unitary $u=s(s^*s)^{-1/2}$ has the same intertwining
property because $s^*s$ commutes with $\varphi(M_r)$. Functional calculus
gives continuity. If the embeddings agree, $s=1$.
\end{proof}

\begin{lemma}[Continuous correction]\label{lem:correction}
For fixed $r$, a tuple $(b_{ij})$ in a unital $C^*$-subalgebra $B\subseteq C$
which is sufficiently close to a system $(p_{ij})$ of unital matrix units in
$C$ can be corrected continuously to unital matrix units $(f_{ij})$ in $B$.
The corrected units approach $(p_{ij})$ uniformly as the error tends to zero,
with the allowed error depending only on $r$. The procedure fixes exact
unital matrix units.
\end{lemma}
\begin{proof}
Here is a construction, including the relative property. Replace diagonal
entries by their self-adjoint parts. For $r>1$, define successively, for
$1\leq i<r$,
\[
 e^{(i)}=1-\sum_{j<i}q_j,\qquad
 q_i=\chi_{(1/2,\infty)}
       \bigl(e^{(i)}b_{ii}e^{(i)}\bigr),\qquad
 q_r=1-\sum_{i<r}q_i,
\]
where the spectral projection is taken in the corner $e^{(i)}Be^{(i)}$.
For sufficiently small error the spectra avoid $1/2$, and these orthogonal
projections are close to $p_{ii}$. Put $v_1=q_1$ and, for $i>1$, put
\[
 a_i=q_i b_{i1}q_1,\qquad
 v_i=a_i(a_i^*a_i)^{-1/2},
\]
where the inverse is taken in $q_1Bq_1$. It exists for small error. Then
$v_i^*v_i=q_1$. Its range projection is a subprojection of $q_i$ close to
$q_i$, so it equals $q_i$: two distinct comparable projections have norm
distance $1$. Therefore $f_{ij}=v_iv_j^*$ are unital matrix units in $B$.

All formulas are continuous on a sufficiently small neighborhood of the
exact relations, using fixed spectral gaps. The bounds can be chosen in
terms of $r$ alone by induction through these finitely many operations.
For exact input every operation returns the original matrix units. For
$r=1$, set $f_{11}=1$.
\end{proof}

\begin{proposition}[Relative families]\label{prop:relative}
Let $P$ be compact Hausdorff and $Q\subseteq P$ closed. Suppose
$\varphi:P\to\Emb(D_n,D)$ is continuous and
$\varphi(Q)\subseteq\Emb(D_n,D_{N_0})$. There is $N\geq\max(n,N_0)$ and a
homotopy relative to $Q$ from $\varphi$ to a map with values in
$\Emb(D_n,D_N)$.
\end{proposition}
\begin{proof}
Let $P_N:D\to D_N$ be the conditional expectation given by the normalized
trace on the tail. These contractions converge in norm to the identity on
each element of $D$, uniformly on compact subsets. Consequently
\[
 b_{ij}(x)=P_N\bigl(\varphi_x(e_{ij})\bigr)
\]
are uniformly close to the matrix units $\varphi_x(e_{ij})$ when $N$ is large.
Lemma~\ref{lem:correction} gives continuous exact matrix units $f_{ij}(x)$ in
$D_N$, hence a continuous family $\psi:P\to\Emb(D_n,D_N)$ close to $\varphi$.
For $x\in Q$, $P_N$ fixes the original matrix units, and the correction fixes
them as well; thus $\psi_x=\varphi_x$.

Use Lemma~\ref{lem:close} to write $\psi_x=\Ad(u_x)\varphi_x$ with $u_x$
uniformly close to $1$. The homotopy
\[
 \varphi_{x,t}=\Ad\bigl(\exp(t\log u_x)\bigr)\varphi_x
 \qquad(0\leq t\leq1)
\]
is well defined using the logarithm near $1$. It fixes $Q$ pointwise.
\end{proof}

\begin{corollary}\label{cor:I4}
The canonical map
\[
 \operatorname*{colim}_{N\geq n}\Emb(D_n,D_N)
 \longrightarrow E_n:=\Emb(D_n,D)
\]
is a weak equivalence, naturally in $n$ with respect to all morphisms of
$\C$. The colimit has the direct-limit topology.
\end{corollary}
\begin{proof}
The connecting maps are closed smooth embeddings of compact manifolds and
therefore cofibrations. Their mapping telescope is weakly equivalent to the
colimit. Maps from spheres into a finite subtelescope can be moved to a
single stage; compactness of the telescope height coordinate reduces every
sphere to this case. Proposition~\ref{prop:relative}, applied to spheres
relative to a basepoint and to disks relative to their boundaries, proves
surjectivity and injectivity on homotopy groups, as well as the assertion on
path components. The maps in the statement commute strictly with
precomposition. No natural choice of the perturbations is asserted or needed.
\end{proof}

\section{The limit algebra and the nerve of the matrix category}\label{sec:K}

\subsection{Restriction to a finite stage}

Let
\[
 I_n=\Iso(D_n\otimes A,D),\qquad
 p_n:I_n\longrightarrow E_n,\quad p_n(e)(a)=e(a\otimes1).
\]

\begin{lemma}\label{lem:restriction}
The map $p_n$ is a locally trivial principal $H$-bundle. In particular,
$p_n$ is a weak equivalence.
\end{lemma}
\begin{proof}
Every embedding $\varphi:D_n\to D$ is unitarily conjugate to the standard
one. Indeed, apply Proposition~\ref{prop:relative} to a point, using its
proof to obtain a close finite-stage embedding $\psi$ and an exact unitary
conjugacy between $\varphi$ and $\psi$. Two unital embeddings of $D_n$ into
the same sufficiently large matrix stage are unitarily conjugate there.
The standard embedding has centralizer isomorphic to $A$ and extends to
$D_n\otimes A\cong D$. Conjugating this extension proves that $p_n$ is onto.

The right action is $e\cdot h=e\circ(\id_{D_n}\otimes h)$. Two extensions
of one embedding differ by a unique such automorphism: an automorphism
fixing $D_n\otimes1$ restricts to an automorphism of its commutant
$1\otimes A$, and the two commuting factors generate the tensor product.

For local triviality, fix an extension $e_0$ of $\varphi_0$. For $\varphi$
near $\varphi_0$, Lemma~\ref{lem:close} gives a local section
$\varphi\mapsto\Ad(u(\varphi,\varphi_0))e_0$. Its associated principal-bundle
trivialization is continuous; composition and inversion of isomorphisms
are continuous in point norm. Finally use~\eqref{eq:H}.
\end{proof}

\subsection{Adjoining the limit algebra}

Let $\widehat{\C}$ have the objects of $\C$ and one further object $D$,
with discrete object topology. Keep all finite-stage morphisms and put
\[
 \widehat{\C}(D_n,D)=E_n,\qquad
 \widehat{\C}(D,D)=G,\qquad
 \widehat{\C}(D,D_n)=\varnothing .
\]
Composition is composition of homomorphisms and is continuous in point
norm. Write $j:\C\to\widehat{\C}$ and $u:*G\to\widehat{\C}$ for the
two full inclusions. We will use their classifying spaces to compare the
finite stages with the automorphisms of their limit.

For a continuous functor $J:\mathcal A\to\mathcal D$, the objects of
$(J\downarrow b)$ are arrows $J(a)\to b$, with the topology inherited
from the arrow and object spaces. Morphisms are commuting triangles.
The dual construction $(b\downarrow J)$ uses arrows $b\to J(a)$.
For unital topological categories with discrete object spaces,
topological Quillen Theorem~A says that $\B J$ is a weak equivalence if
all $\B(J\downarrow b)$ are weakly contractible. The dual criterion uses
$\B(b\downarrow J)$. Indeed the source, target and comma-object
projections required in~\cite[Theorems~4.7 and~4.8, condition~(v)]{ERW}
are maps to discrete spaces and hence Serre fibrations. In the proof
of those theorems, weak contractibility of each comma space suffices:
the level maps in the resolution then have weakly contractible homotopy
fibers and hence are weak equivalences.
The morphism spaces retain their topology throughout.

\begin{proposition}\label{prop:two-inclusions}
Both arrows in the following zigzag are weak equivalences:
\begin{equation}\label{eq:rho}
 \B\C\xrightarrow{\,\B j\,}\B\widehat{\C}
       \xleftarrow{\,\B u\,}BG .
\end{equation}
\end{proposition}
\begin{proof}
\emph{The finite-stage inclusion.}
For a finite target $D_N$, $(j\downarrow D_N)=\C/D_N$ has terminal
object $(D_N,\id)$. The unique arrow to it from $(D_n,f)$ is $f$,
continuously in $f$. Its classifying space is contractible.

For target $D$, a $p$-simplex in the comma nerve is a string
\[
 D_{n_0}\longrightarrow\cdots\longrightarrow D_{n_p}
                   \xrightarrow{\varphi} D .
\]
Write $Y^N_\bullet=N(\C/D_N)$ and
$Y_\bullet=N(j\downarrow D)$. Composition with $D_N\hookrightarrow
D_{N+1}$ gives a sequence $Y^N_\bullet$ over $Y_\bullet$.
For a fixed string of indices, put
$P=\prod_{a=1}^p\C(n_{a-1},n_a)$. On this summand the map
$\operatorname{Tel}_N Y^N_p\to Y_p$ is
\[
 P\times\operatorname{Tel}_{N\ge n_p}\Emb(D_{n_p},D_N)
       \longrightarrow P\times E_{n_p}.
\]
It is a weak equivalence by Corollary~\ref{cor:I4}; the earlier
summands are empty. These maps commute with faces, including the last
face, which composes $\varphi$ with the preceding arrow. Fat realization
preserves the levelwise weak equivalence. It commutes with the telescope:
both constructions use coproducts, products with a simplex or interval,
and the indicated face or cylinder identifications in compactly generated
spaces. Consequently
\[
 \B(j\downarrow D)\ \simeq \operatorname{Tel}_N\B(\C/D_N)\ \simeq\ * .
\]
For the second equivalence, a sphere in the telescope has bounded height.
It moves into its last stage in that finite subtelescope and contracts
there. The telescope is path connected for the same reason. Theorem~A
now applies to $j$.

\emph{The automorphism inclusion.}
At $D$, the objects of $(D\downarrow u)$ are elements of $G$ and there
is a unique arrow from $g$ to $g'$, given by $g'g^{-1}$. The identity
element is a terminal object. At $D_n$, the objects of
$(D_n\downarrow u)$ are embeddings $\varphi:D_n\to D$ and an arrow
$\varphi\to\psi$ is a $g\in G$ with $g\varphi=\psi$.
Lemma~\ref{lem:restriction}, after choosing $D_n\otimes A\cong D$,
identifies the orbit map $G\to E_n$ as a locally trivial principal
$H$-bundle. Thus the action is transitive with contractible stabilizer.

To see what this implies for the comma space, forget the group elements
in a string of $p$ arrows, retaining its $p+1$ embeddings. Locally choose
lifts $s_i\in G$ of these embeddings through the orbit map. Every arrow
from vertex $i-1$ to vertex $i$ is then uniquely
$s_i h_i s_{i-1}^{-1}$, with $h_i\in H$. The forgetful map of nerves
in degree $p$ is therefore a locally trivial bundle over $E_n^{p+1}$
with contractible fiber $H^p$. The target is the nerve of the category
on $E_n$ with a unique arrow between each pair of objects. Choosing any
one embedding gives it a terminal object. Fat realization shows that
$\B(D_n\downarrow u)$ is weakly contractible. The dual Theorem~A
applies to $u$.
\end{proof}

Define $\rho=(\B u)^{-1}\B j$ in the homotopy category, using the
ordinary-realization comparison in Section~\ref{sec:Ex} on the source.
For $k=1$ there is also an immediate check: $D_0=\mathbb C$ is initial
in $\C$, and $G=\Aut(A)$ is contractible. The proof of
Proposition~\ref{prop:two-inclusions} works for every $k\ge1$, $l\ge2$.

\section{Ordinary realization and the space \texorpdfstring{$\Ex(N\C)$}{Ex(NC)}}\label{sec:Ex}

Define $T=\Ex(N\C)$ by
\[
 T_p=\operatorname{Map}_{s\mathrm{Top}}(\Sd\Delta[p],N\C),
\]
where the source is simplicially discrete and the mapping object has its
finite-limit topology. Equivalently, $T_p$ is the space of functors from
the poset of nonempty subsets of $[p]$ to $\C$.

\begin{lemma}\label{lem:proper}
The simplicial spaces $N\C$ and $T$ are proper: their latching maps are
closed cofibrations.
\end{lemma}
\begin{proof}
Let $X_\bullet$ denote either $N\C$ or $T$. In degree $p$, fix an
assignment $\lambda$ of object indices and write $X_{p,\lambda}$ for the
corresponding open and closed summand. For $N\C$ it is a finite
product of compact spaces of matrix-unit systems. For $T$ the indices
are assigned to the nonempty subsets of $[p]$, and the equations for
composition cut out a closed subset of such a product. In matrix-unit
coordinates all these equations are polynomial in real and imaginary
parts. Thus $X_{p,\lambda}$ is compact, Hausdorff, and semialgebraic.

Let $D_{p,\lambda}$ be its degenerate locus. For each $i$, the
identity $d_i s_i=\id$ identifies $\operatorname{im}(s_i)$ with the
fixed locus of $s_i d_i$. On a summand whose assignment is preserved
by $s_i d_i$, this is a closed semialgebraic set. On every other
summand the fixed locus is empty. The finite union $D_{p,\lambda}$ is
therefore closed and semialgebraic.

The simplicial identities identify the latching space set-theoretically
with the union of the degeneracy images. Over $\lambda$, it is a
quotient of finitely many compact pieces. A surjection $[p]\to[q]$
either admits no compatible lower-degree assignment or determines it
uniquely. Hence the continuous bijection from this part of the latching
space to $D_{p,\lambda}$ is a homeomorphism. Compatible semialgebraic
triangulation of $(X_{p,\lambda},D_{p,\lambda})$ makes the latter a
subcomplex~\cite{BCR}. Its inclusion is a closed cofibration, and
taking the coproduct over $\lambda$ gives the latching cofibration.
The same triangulation argument for each individual degeneracy image
also gives the degeneracy cofibrations used in Segal's good-space
formulation of the realization comparison.
\end{proof}

\begin{proposition}\label{prop:Ex}
The natural map $N\C\to\Ex(N\C)$ adjoint to the last-vertex map
induces a weak equivalence on ordinary realizations.
\end{proposition}
\begin{proof}
In degree $p$, the comparison sends a simplex
$\sigma:\Delta[p]\to N\C$ to its composite with the last-vertex map
$\Sd\Delta[p]\to\Delta[p]$. For each singular degree $q$, preservation
of finite limits gives a natural
identification of simplicial sets
\[
 \Sing_q\bigl(\Ex(N\C)\bigr)_\bullet
 \cong \Ex\bigl(\Sing_q(N\C)_\bullet\bigr).
\]
The simplicial-set subdivision theorem says that $Z\to\Ex Z$ is a weak
equivalence for every simplicial set $Z$. Hence the indicated map of
bisimplicial sets is a weak equivalence in each row, and thus on diagonals.
For any simplicial space $X_\bullet$, the fat realization and
$|\diag\Sing X_\bullet|$ have naturally equivalent weak homotopy
types (\cite[pp.~2138--2139, before Lemma~6.9]{ERW}, or by iterated realization
and the degreewise singular evaluation maps). Thus the map on fat
realizations is a weak equivalence. Lemma~\ref{lem:proper} and the
fat-to-ordinary comparison give the result for ordinary realizations.
\end{proof}

In particular $\B\C\to|N\C|$ is a homotopy equivalence. Combining
Propositions~\ref{prop:two-inclusions} and~\ref{prop:Ex} proves
Theorem~\ref{thm:main}.

\section{Strict LABs and categorical rectification}\label{sec:strict}

For the geometric arguments, $X$ is a finite polyhedron. Covers are finite;
they can always be refined to good covers by open stars in a sufficiently
fine subdivision. We first allow arbitrary singleton ranks.

\begin{definition}\label{def:strict}
Let $\mathcal U=\{U_i\}_{i\in I}$ be a finite indexed open cover.
For each nonempty $S\subseteq I$ write $U_S=\bigcap_{i\in S}U_i$.
A strict LAB$'$ of type $(k,l)$ consists of a matrix algebra bundle
$A_S\to U_S$ for each $S$ with $U_S\ne\varnothing$, of fiber
$D_{n_S}$, and unital bundle embeddings
$\mu_{S,T}:A_S|_{U_T}\to A_T$ for $S\subseteq T$. We require
$\mu_{S,S}=\id$ and $\mu_{T,R}\mu_{S,T}=\mu_{S,R}$.
A normalized LAB additionally has $n_{\{i\}}=0$.
A morphism $\iota:\mathfrak A\to\mathfrak B$ is a family of unital
embeddings satisfying $\iota_T\mu^A_{S,T}=\mu^B_{S,T}\iota_S$.
\end{definition}

When $S\cap T\ne\varnothing$, the square associated with
$S\cap T,S,T,S\cup T$ already commutes by composition; it is not an
additional coherence condition. No algebra $A_\varnothing$ is part
of the definition. A refinement
$\alpha:J\to I$ pulls $A_S$ back to $A_{\alpha(S')}$ on the new
intersection indexed by $S'\subseteq J$. Repeated values of $\alpha$
are removed when forming the set $\alpha(S')$.

Let $\mathsf F'(X)$ be the set of LAB$'$s modulo the equivalence relation
\emph{generated} by morphisms, refinements, and the two restrictions of
LAB$'$s on $X\times I$. Define $\mathsf F(X)$ using normalized LABs and
normalized intermediates. A span into a common diagram and a zigzag of
morphisms generate the same relation, since a single morphism is a span
with an identity leg. Endpoint isomorphisms and refinements are allowed.
Using the generated relation avoids assuming a gluing theorem for bare
cylinders before one has proved it. Pullback along a homotopy gives
homotopy invariance of these functors directly.

For a good cover let $\mathcal U_\cap$ have objects $(x,S)$, $x\in U_S$,
and the unique arrow $(x,S)\to(x,T)$ for $S\subseteq T$. Segal's
numerable covering theorem, in its barycentric subdivision model
\cite{Segal68}, gives a weak equivalence $\epsilon:\B\mathcal U_\cap\to X$. Trivializations of all $A_S$ define a
continuous functor $a:\mathcal U_\cap\to\C$ and hence a class
\[
 \Phi(\mathfrak A)=[\B a\circ\epsilon^{-1}]\in[X,|N\C|],
\]
where $\epsilon^{-1}$ is taken in the homotopy category and $\B\C$ is identified with $|N\C|$. Changing trivializations gives a natural isomorphism. A LAB morphism gives
a natural transformation and therefore a homotopy on classifying spaces.
Refinement commutes with the projections to $X$. Two refinement maps
$\alpha,\alpha'$ are compared through the diagram with terms
$A_{\alpha(S')\cup\alpha'(S')}$, which receives morphisms from both
pullbacks. Finally a cylinder gives a class on $X\times I$. Consequently
\begin{equation}\label{eq:Phi-descends}
 \bar\Phi:\mathsf F'(X)\longrightarrow[X,|N\C|]
\end{equation}
is well defined and natural. Define $\bar R=\rho_*\bar\Phi$.

\subsection{Compatibility with stabilized arrows}

Define $\K=\K_{\C}$ with objects $n\in\mathbb N_0$. For $n\leq m$ a
morphism $m\to n$ is a pair $(\iota,\Theta)$ with
\[
 \iota\in\C(n,m),\qquad
 \Theta:D_m\otimes A\xrightarrow{\cong}D_n\otimes A,
 \qquad \Theta(\iota(a)\otimes1)=a\otimes1.
\]
Use the subspace topology from the product of the embedding and
isomorphism spaces. Composable morphisms are composed by composing the
embeddings in the opposite direction and the isomorphisms in their displayed
direction. This preserves the defining equation and is continuous.

\begin{lemma}\label{lem:forget}
The forgetful functor $q:\K\to\C^{\op}$ induces a degreewise weak
equivalence of nerves, hence $\B q:\B\K\simeq\B\C^{\op}$.
\end{lemma}
\begin{proof}
For $\iota:D_n\to D_m$, its centralizer $Z(\iota)$ is a matrix algebra
$M_{l^{m-n}}$, and multiplication identifies
$D_m\cong D_n\otimes Z(\iota)$. An admissible $\Theta$ is exactly an
isomorphism $Z(\iota)\otimes A\cong A$ on the second factor. Its space is
an $H$-torsor. These torsors form a locally trivial bundle over $\C(n,m)$:
local continuous conjugating unitaries for nearby finite-dimensional
embeddings transport a chosen $\Theta$ to a local section. Equivalently,
one can use local sections of the compact Lie-group orbit map.

In nerve degree $p$, over each fixed object string the forgetful map is a
product of $p$ such bundles, with contractible fiber $H^p$. The map is the
identity in degree zero. Fat realization preserves this degreewise weak equivalence.
\end{proof}

Choose isomorphisms $\varepsilon_n:D_n\otimes A\cong D$. There is a strict
continuous functor
\[
 F:\K\longrightarrow *G,\qquad
 F(\iota,\Theta)=\varepsilon_n\Theta\varepsilon_m^{-1}.
\]
The arrows in $\K$ point in the direction of $\Theta$, so this formula is
covariant.

We compare this functor with~\eqref{eq:rho}. Let $q^+:\K^{\op}\to\C$
be the opposite of $q$, and let $F^+:\K^{\op}\to *G$ send the arrow
opposite to $f$ to $F(f)^{-1}$. The maps $s_n:D_n\to D$,
$s_n(a)=\varepsilon_n(a\otimes1)$, define a continuous natural
transformation $j q^+\Rightarrow uF^+$. Its naturality is the identity
\[
 F(\iota,\Theta)^{-1}s_n=s_m\iota,
\]
obtained by applying $\Theta^{-1}$ to $a\otimes1$.
Thus $\rho\,\B q^+=\B F^+$ in the homotopy category.

For clarity about the direction of the transition functions, let
$r_{\mathcal D}:\B\mathcal D\to\B(\mathcal D^{\op})$ reverse each
simplex and its barycentric coordinates. For a group, reversal followed
by inversion on arrows induces a map homotopic to the identity on $BG$.
In homogeneous bar coordinates a point is a weighted tuple
$(g_0,\ldots,g_p;t)$ modulo simultaneous right translation.
The reversed tuple has vertices $(g_p,\ldots,g_0)$ and reversed weights
$t^{\rm rev}$. Concatenate the original and reversed tuples with weights
$((1-v)t,v t^{\rm rev})$, $0\le v\le1$. Deleting zero weights
commutes with concatenation, and the formula is equivariant for right
translation. It therefore descends to the claimed homotopy on the fat
nerve. Applying reversal to the preceding naturality identity now gives
\[
 \rho\,r_{\C^{\op}}\,\B q\simeq\B F .
\]
This identifies the map of Theorem~\ref{thm:main} with the original
stabilized-arrow rectification, including its variance. A change of
$\varepsilon_n$ gives an objectwise natural isomorphism.

\begin{remark}[The homotopy fiber via Quillen Theorem B]
The topological comma category $(F\downarrow *)$ has objects $(n,g_n)$.
An arrow over $f:m\to n$ in $\K$ satisfies $g_m=g_nF(f)$.
Every $h\in G$ acts by $g_n\mapsto hg_n$, an isomorphism of topological
categories. This verifies the transport hypothesis of topological
Theorem~B, and the object-space fibration hypotheses hold as above.
Thus $\B(F\downarrow *)$ models the homotopy fiber of $\B F$
\cite[Theorem~4.9 and its proof, diagram~(4.11)]{ERW}.
Restriction $(n,g_n)\mapsto(n,g_n s_n)$ gives a functor to
$(j\downarrow D)^{\op}$. Over a string of $p$ arrows, choose an
extension $e_i=g_i\varepsilon_{n_i}\in I_{n_i}$ independently at each
vertex. These choices force
all stabilized arrows, since their isomorphisms are the quotients
$e_{\mathrm{target}}^{-1}e_{\mathrm{source}}$ of the chosen extensions.
Thus in each degree the functor is a pullback of a product of the
principal bundles $I_{n_i}\to E_{n_i}$, with fiber $H^{p+1}$.
Its realization is a weak equivalence. The finite-stage telescope in
Proposition~\ref{prop:two-inclusions} therefore also proves weak
contractibility of this homotopy fiber directly.
\end{remark}

For a LAB the construction takes place over its cover. Let
$\K(\mathfrak A)$ have objects $(x,S)$ and arrows $(x,T)\to(x,S)$
the isomorphisms
\[
 \Theta:(A_T)_x\otimes A\xrightarrow{\cong}(A_S)_x\otimes A,
 \qquad \Theta(\mu_{S,T}(b)\otimes1)=b\otimes1.
\]
The arrow bundle over $U_T$ is a locally trivial principal $H$-bundle.
Indeed the embedding of finite-dimensional algebras is locally unitarily
conjugate to a fixed one, and its centralizer is a locally trivial
$M_{l^{n_T-n_S}}$-bundle. An admissible $\Theta$ is a trivialization of
this centralizer after tensoring with $A$. The action on the target tail
is the canonical action of $H$; it does not require a global trivialization
of $A_S$. The same argument as for $\K_\C$ proves continuity of composition.
Forgetting $\Theta$ gives a degreewise weak equivalence to
$N\mathcal U_\cap^{\op}$, so $\pi:\B\K(\mathfrak A)\simeq X$.

Choose $t_S:A_S\cong U_S\times D_{n_S}$ and put
$\varepsilon_S=\varepsilon_{n_S}(t_S\otimes\id_A)$.
The functor $F_{\mathfrak A}(\Theta)=\varepsilon_S\Theta\varepsilon_T^{-1}$
to $*G$ factors strictly through $F:\K_\C\to *G$ via
\[
 \Theta\longmapsto
 \bigl(t_T\mu_{S,T}t_S^{-1},
 (t_S\otimes\id_A)\Theta(t_T^{-1}\otimes\id_A)\bigr).
\]
The forgetful square also commutes. Thus
\begin{equation}\label{eq:rect-compat}
 [\B F_{\mathfrak A}\,\pi^{-1}]
 =\rho_*\bar\Phi(\mathfrak A)=\bar R(\mathfrak A).
\end{equation}
This proves compatibility with the rectification used below, including
its descent to the chosen equivalence relation.

\begin{definition}
A subalgebra-LAB$'$ in a locally trivial $D$-bundle $E$ is a LAB$'$ whose
$A_S$ are subbundles of $E|_{U_S}$ and whose structure maps are inclusions.
The pairs $A_S\subset E$ are required to be locally trivial. This last
condition also follows locally from a continuous framing of $A_S$ and
Lemma~\ref{lem:close}.
\end{definition}

\begin{proposition}\label{prop:ambient}
A subalgebra-LAB$'$ in $E$ has $\bar R(\mathfrak A)=[E]$.
\end{proposition}
\begin{proof}
On a good refinement choose trivializations $t_S:A_S\cong
U_S\times D_{n_S}$ and charts $\tau_S:E|_{U_S}\cong U_S\times D$.
Write $a:\mathcal U_\cap\to\C$ for the resulting finite diagram.
On $S\subseteq T$ the transition $\tau_T\tau_S^{-1}$ defines a functor
$b:\mathcal U_\cap\to *G$ classifying $E$ by the cover construction
\cite[Section~4]{Segal68}. If $i_S:A_S\to E|_{U_S}$ is the given inclusion,
the embeddings $\eta_S=\tau_S i_S t_S^{-1}:D_{n_S}\to D$ form
a natural transformation $ja\Rightarrow ub$. Indeed, on $U_T$,
\[
 (\tau_T\tau_S^{-1})\eta_S
   =\tau_T i_S t_S^{-1}
   =\eta_T\,(t_T\mu_{S,T}t_S^{-1}).
\]
Taking classifying spaces gives $\B j\,\B a\simeq\B u\,\B b$.
The definition $\rho=(\B u)^{-1}\B j$ therefore gives
$\bar R(\mathfrak A)=[E]$.
\end{proof}

\section{Every strict diagram has an ambient bundle}\label{sec:sections}

We use two elementary geometric facts about finite polyhedra. A finite
open cover has a shrinking $V_i\subset C_i\subset U_i$ in which the $V_i$
are open and cover $X$, and all $C_i$ are subcomplexes of one finite
subdivision. To obtain this, first choose a shrinking with compact closures
inside $U_i$, and then subdivide until the closed simplices meeting those
closures lie in $U_i$. Also, a compact subset of an open neighborhood has
a smaller neighborhood covered by finitely many open stars subordinate
to any given cover. These stars form a good cover. Both facts apply
relatively at the two ends of a triangulated cylinder.

\begin{lemma}[Relative section extension]\label{lem:section-extension}
If $(C,L)$ is a finite CW pair and $P\to C$ is a locally trivial principal
$H$-bundle, every section on $L$ extends to $C$.
\end{lemma}
\begin{proof}
Over each characteristic disk the bundle is trivial. The boundary section
is a map from a sphere to the contractible group $H$, hence extends.
Induction over the finitely many cells proves the assertion.
\end{proof}

\begin{lemma}[Compatible sections]\label{lem:functorial-sections}
For the shrinking above, put $C_T=\bigcap_{i\in T}C_i$.
There are admissible sections $\Theta_0(S,T)$ on $C_T$ with
$\Theta_0(T,T)=\id$ and
\begin{equation}\label{eq:theta-composition}
 \Theta_0(R,S)\Theta_0(S,T)=\Theta_0(R,T)
 \quad\text{on }C_T\quad(R\subseteq S\subseteq T).
\end{equation}
\end{lemma}
\begin{proof}
Process the nonempty intersections $T$ in nonincreasing order of $|T|$.
For each $S\subsetneq T$, prescribe on $C_I$, $I\supsetneq T$,
\begin{equation}\label{eq:theta-boundary}
 \Theta_0(S,T)=\Theta_0(S,I)\Theta_0(T,I)^{-1}.
\end{equation}
The right side is admissible: the inverse sends $b\otimes1$ to
$\mu_{T,I}(b)\otimes1$ before the first factor is applied.
On $C_I\cap C_{I'}=C_{I\cup I'}$, earlier prescriptions express both
right sides as
$\Theta_0(S,I\cup I')\Theta_0(T,I\cup I')^{-1}$.
This includes the case that one of $I,I'$ contains the other.
Thus the prescriptions agree on the subcomplex
$L_T=\bigcup_{I\supsetneq T}C_I$.
Lemma~\ref{lem:section-extension} extends each section to $C_T$.
When $R\subsetneq S\subsetneq T$, the later prescription for $(R,S)$
on $C_T\subseteq L_S$ is exactly
$\Theta_0(R,S)=\Theta_0(R,T)\Theta_0(S,T)^{-1}$.
This proves \eqref{eq:theta-composition}; equality cases are immediate.
\end{proof}

\begin{proposition}\label{prop:strict-ambient}
After shrinking the cover, every strict LAB$'$ is isomorphic to a
subalgebra-LAB$'$ in a locally trivial $D$-bundle representing its
rectification.
\end{proposition}
\begin{proof}
For any $S,T$ with $V_{S\cup T}\ne\varnothing$, define
\[
 g_{TS}=\Theta_0(T,S\cup T)\Theta_0(S,S\cup T)^{-1}:
 (A_S\otimes A)|_{V_{S\cup T}}\longrightarrow
 (A_T\otimes A)|_{V_{S\cup T}}.
\]
On a triple intersection, put $Q=R\cup S\cup T$.
Equation~\eqref{eq:theta-composition} gives
$g_{TS}=\Theta_0(T,Q)\Theta_0(S,Q)^{-1}$ there. Writing
$a_J=\Theta_0(J,Q)$ for $J=R,S,T$, we obtain the cocycle identity
\[
 g_{TR}g_{RS}=(a_Ta_R^{-1})(a_Ra_S^{-1})=a_Ta_S^{-1}=g_{TS}.
\]
Ordinary bundle gluing on the open cover
$\{V_S\}$ produces a locally trivial bundle $E$ and identifications
$e_S:(A_S\otimes A)|_{V_S}\to E|_{V_S}$.
For $S\subseteq T$, the gluing equation gives
$e_T(\mu_{S,T}(a)\otimes1)=e_S(a\otimes1)$.
The images of the $A_S$ are therefore strictly nested. Local triviality
of the pairs follows in the $S$-chart. The shrunk cover need not be good;
take a further good refinement before applying
Proposition~\ref{prop:ambient}. Refinement invariance then shows that
the rectification of the original diagram is $E$.
\end{proof}

\section{A relative simultaneous gauge lemma}\label{sec:gauge}

All unitaries in this section have the norm topology. Quantitative constants
are unnecessary: continuity with uniform error control, depending only on
the finite matrix size, suffices.

\begin{lemma}[Relative perturbation]\label{lem:gauge-relative}
Fix $m$. Suppose $C$ is compact metrizable, $L\subset C$ is closed, and
$\varphi:C\to\Emb(M_m,D)$ is continuous. Assume that all its matrix
units have distance at most $\epsilon$ from $D_N$. Suppose
$v:L\to U(D)$ is continuous, $\|v-1\|\leq\eta$, and
$v\varphi(M_m)v^*\subset D_N$ on $L$.
For sufficiently small $\epsilon,\eta$ there is an extension
$u:C\to U(D)$ of $v$ with $u\varphi(M_m)u^*\subset D_N$ everywhere.
Its distance from $1$ tends uniformly to zero as
$\epsilon+\eta\to0$, with bounds depending only on $m$.
\end{lemma}
\begin{proof}
Choose $\eta$ small enough for the logarithm near $1$. By
Dugundji's extension theorem~\cite{Dugundji}, extend $\log v$ to a
continuous function $s:C\to D^{\rm skew}$ with the same uniform norm
bound. Set $w=\exp(s)$, taking $w=1$ when $L$ is empty.
Then $w|_L=v$ and $\|w-1\|_\infty$ tends to zero with $\eta$.

Put $\psi_y=\Ad(w(y))\varphi_y$. Its matrix units are within
$\epsilon+2\|w-1\|_\infty$ of $D_N$, and on $L$ they already belong
to $D_N$. Apply $P_N$ followed by the correction of
Lemma~\ref{lem:correction}. The resulting exact units in $D_N$ agree
with those of $\psi$ on $L$. By Lemma~\ref{lem:close} the conjugating
unitary $z$ from $\psi$ to the corrected system is continuous, uniformly
close to $1$, and equal to $1$ on $L$. Thus $u=zw$ extends $v$ and
puts the original units into $D_N$ everywhere. The error bound tends
to zero with $\epsilon+\eta$, uniformly for the given matrix size.
\end{proof}

\begin{lemma}[Simultaneous gauge]\label{lem:simultaneous}
Let $Y$ be compact metrizable, covered by finitely many closed sets $C_i$.
For each nonempty $C_S$ suppose a framing
$\varphi_S:C_S\to\Emb(M_{kl^{n_S}},D)$ is given. Write
$A_S(y)=\varphi_S(y)(M_{kl^{n_S}})$, and assume
$A_S(y)\subseteq A_T(y)$ for $y\in C_T$, $S\subseteq T$.
For every $\delta>0$ there exist $N$ and a continuous
$u:Y\to U(D)$ with $\|u-1\|<\delta$ such that
$uA_Su^*\subset D_N$ on every $C_S$.
\end{lemma}
\begin{proof}
Compactness and finiteness make all matrix units uniformly close to one
stage $D_N$, with an arbitrarily small error $\epsilon$.
List the nonempty intersections $I_1,\ldots,I_r$ in nonincreasing
cardinality. Inductively extend a gauge from
$B_{a-1}=\bigcup_{b<a}C_{I_b}$ to $B_a=B_{a-1}\cup C_{I_a}$.
Its already defined part on $C_{I_a}$ is the closed set
\[
 L=B_{a-1}\cap C_{I_a}=\bigcup_{b<a}C_{I_b\cup I_a}.
\]
Every nonempty union index $I_b\cup I_a$ occurs earlier, so the old gauge
already puts $A_{I_b\cup I_a}\supseteq A_{I_a}$ into $D_N$ there.
Lemma~\ref{lem:gauge-relative} extends that gauge to $C_{I_a}$,
and closed-set pasting defines it on $B_a$. Nothing previously imposed is
changed. The number of steps and matrix sizes are finite. Choose
$\epsilon$ small enough that the successive error bounds, each tending
to zero at zero, remain below $\delta$ and within every logarithm
neighborhood used. At the end $B_r=Y$.
\end{proof}

\section{Relative existence and uniqueness in an ambient bundle}\label{sec:extension}

The following extension argument uses compact chart domains. It does not
require identifying independently chosen collars of different cells.

\begin{lemma}[One extension step]\label{lem:extension-step}
Let $E\to X$ be a $D$-bundle, $L\subset O\subset X$ with $L$ compact
and $O$ open, and let a subalgebra-LAB$'$ $\mathfrak A$ be given on a
finite good cover $\{U_i\}$ of $O$. Let $K\subset N$ be compact, where
$N$ is an open trivializing chart of $E$.
There is a subalgebra-LAB$'$ on a neighborhood of $L\cup K$ whose
restriction near $L$ is a refinement of $\mathfrak A$ by the original
inclusions in $E$.
\end{lemma}
\begin{proof}
Put $Z=L\cap K$. Choose an open $Q$ containing $Z$ with compact
closure in $O\cap N$. If $Z=\varnothing$, take $Q=\varnothing$.
Choose compact $C_i\subset U_i\cap N$ whose interiors $V_i$ cover
$\overline Q$, allowing empty $C_i$. The old algebras are framed over
each $C_S$ because $U_S$ is contractible. On
$Y^*=\bigcup_i C_i$ apply Lemma~\ref{lem:simultaneous} in the chart
of $E$. It gives $u:Y^*\to U(D)$ within $1/4$ of $1$ and a stage
$D_r$ containing every $uA_Su^*$ on $C_S$. We may increase $r$
to dominate all old exponents. Extend $\log u$ from $Y^*$ to a
skew-adjoint function on $N$ by Dugundji's theorem. Its exponential
$w$ agrees with $u$ on $Y^*$. For $Y^*=\varnothing$, take $w=1$.
The new matrix algebra in this chart will be $w^*D_rw$.

The compact set $K\setminus Q$ is disjoint from $L$. Choose an open
$\nu\supset K\setminus Q$ with $\overline\nu\subset N$ and
$\overline\nu\cap L=\varnothing$, and put $O^-=O\setminus\overline\nu$.
Use the patches
\[
                  U_i\cap O^-,\qquad V_i,\qquad \nu .
\]
They cover a neighborhood of $L\cup K$: the first family covers $L$,
the $V_i$ cover $K\cap Q$, and $\nu$ covers $K\setminus Q$.

On intersections of old patches use the old algebra, removing duplicate
old indices. On every intersection involving $\nu$, including $\nu$
itself, use $w^*D_rw$ in the chart. If an old intersection meets
$\nu$, all of its patches are of the form $V_i$. On that intersection
$w=u$ and the old algebra lies in $w^*D_rw$ by construction.
All structure maps are therefore inclusions. The new algebra is locally
trivial by its displayed framing. Near $L$ the patch $\nu$ is absent,
so the diagram there is a refinement of the original diagram in $E$.
Take a finite good refinement on a smaller neighborhood of $L\cup K$
as in Section~\ref{sec:sections}.
\end{proof}

\begin{theorem}[Relative existence after refinement]\label{thm:relative-existence}
Let $Y\subset X$ be compact. Suppose a subalgebra-LAB$'$ $\mathfrak A_0$
in $E$ is defined over a neighborhood $O_Y$ of $Y$, with a finite good
cover. Then $E$ carries a subalgebra-LAB$'$ on all of $X$ whose restriction
to a smaller neighborhood of $Y$ is a refinement of $\mathfrak A_0$
inside the given bundle $E$.
\end{theorem}
\begin{proof}
Choose finitely many compact sets $K_1,\ldots,K_s$ covering $X$, each
inside a trivializing chart. Starting with $L_0=Y$, apply
Lemma~\ref{lem:extension-step} with $L_j=Y\cup K_1\cup\cdots\cup K_j$.
After each step take the indicated good refinement near $L_j$.
At every step the old diagram restricts by refinement near $Y$.
Intersect the finitely many neighborhoods and compose their index maps.
The last neighborhood contains $X$, proving the assertion.
\end{proof}

\begin{corollary}\label{cor:ambient-unique}
Every $D$-bundle on $X$ has a subalgebra-LAB$'$. Two subalgebra-LAB$'$s
in the same $E$ have the same class in $\mathsf F'(X)$.
\end{corollary}
\begin{proof}
For existence take $Y=\varnothing$ in
Theorem~\ref{thm:relative-existence}.
For uniqueness, work in $\operatorname{pr}_X^*E$ on $X\times I$.
Choose a finite good refinement of each endpoint cover and extend
the refined diagrams constantly on disjoint endpoint collars. The
collar patches form a finite good cover of a neighborhood of
$Y=X\times\{0,1\}$. The relative theorem gives a diagram over the
cylinder whose ends refine the original diagrams, giving equality of
their concordance classes.
\end{proof}

\begin{theorem}[Strict classification before normalization]\label{thm:strict-prime}
For finite polyhedra $X$, categorical rectification is a natural bijection
\[
 \bar R:\mathsf F'(X)\xrightarrow{\cong}[X,B\Aut(D)].
\]
\end{theorem}
\begin{proof}
Surjectivity follows from the existence statement and
Proposition~\ref{prop:ambient}. For injectivity,
Proposition~\ref{prop:strict-ambient} realizes each of two diagrams as a
subalgebra diagram in a bundle representing its rectification.
If their rectifications agree, the numerable bundles are isomorphic.
Transport both diagrams into one bundle and apply
Corollary~\ref{cor:ambient-unique}. Naturality was proved in
Section~\ref{sec:strict}.
\end{proof}

\section{Normalization and the strict comparison theorem}\label{sec:normalization}

\begin{proposition}\label{prop:normalization}
The inclusion of normalized diagrams induces a natural bijection
$\mathsf F(X)\xrightarrow{\cong}\mathsf F'(X)$.
\end{proposition}
\begin{proof}
On a good cover, trivialize each $A_i$ and choose a unital rank-$k$
subbundle $B_i\subset A_i$. Replace singleton algebras by $B_i$, leaving
all terms with at least two indices unchanged. Call the resulting
normalized diagram $\mathfrak A_B$. Its structure maps out of a
singleton are the original maps restricted to $B_i$.

Different choices of cores give the same concordance class. In a fixed
trivialization, cores are maps from the contractible $U_i$ into the
connected homogeneous space of unital $M_k$-subalgebras of
$M_{kl^{n_i}}$. Any two such maps are homotopic: contract the domain,
join the two resulting constant values by a path, and reverse the
contraction. This produces a continuous family of singleton subbundles.
Keep every other algebra fixed on the cylinder. The resulting inclusions
define a normalized concordance.

Let $\iota:\mathfrak A\to\mathfrak B$ be a LAB$'$ morphism. Given
source cores $B_i\subset A_i$, choose $\iota_i(B_i)$ as target cores.
Restrict $\iota_i$ to $B_i$ on singletons and retain $\iota_T$ on
intersections with $|T|\geq2$. For $i\in T$, compatibility follows by
restricting the original morphism equation to $B_i$, which gives
\[
 \iota_T\,\mu^{\mathfrak A}_{\{i\},T}|_{B_i}
 =\mu^{\mathfrak B}_{\{i\},T}\,\iota_i|_{B_i}.
\]
These components give a normalized morphism. Any independently chosen
target cores give a concordant normalized diagram by the preceding
paragraph. Normalization therefore respects morphisms on classes.

Now consider a refinement with index map $\alpha:J\to I$.
Use $B_{\alpha(j)}$ as the new singleton cores. There is a canonical
normalized morphism
\begin{equation}\label{eq:normalization-refinement}
 \alpha^*(\mathfrak A_B)\longrightarrow(\alpha^*\mathfrak A)_B.
\end{equation}
If $|S'|\geq2$ and $\alpha(S')=\{i\}$, its component is the inclusion
$B_i\hookrightarrow A_i$. At a singleton it is the identity on the
core, and at every remaining index it is the identity on the algebra.
These components commute with all structure maps.

Now let $\mathfrak A$ be given on an arbitrary finite cover
$\mathcal U$. Choose two good refinements $\mathcal V,\mathcal W$
of $\mathcal U$ and a common finite good refinement
$\mathcal T=\{T_j\}$ of $\mathcal V$ and $\mathcal W$. Such a
refinement is given by the vertex stars of a sufficiently fine
subdivision. Write $p,q$ for the two composite index maps from
$\mathcal T$ to $\mathcal U$. For $T_S=\bigcap_{j\in S}T_j\ne\varnothing$,
put
\[
 C_S=A_{p(S)\cup q(S)}\big|_{T_S}.
\]
This restriction is defined because $T_S\subseteq U_{p(S)\cup q(S)}$.
The structure maps of $\mathfrak A$ give a diagram $\mathfrak C$ with
terms $C_S$ and morphisms
\[
 p^*\mathfrak A\longrightarrow\mathfrak C
                   \longleftarrow q^*\mathfrak A
\]
on $\mathcal T$. Pull the chosen cores on $\mathcal V$ and $\mathcal W$
back to $\mathcal T$. For each arrow, take the images of its source
cores as the target cores in $\mathfrak C$. This gives two normalized
morphisms whose targets are concordant by the core-choice argument on
$\mathcal T$. Equation~\eqref{eq:normalization-refinement} compares
each normalized source with the pullback of the normalized diagram
on its original good refinement. Thus the two normalizations determine
the same class in $\mathsf F(X)$, for every choice of the index maps.
A morphism on an arbitrary cover can now be normalized on a single
good refinement of its source and target. For a refinement of
arbitrary covers one may first refine the finer cover by a good cover
and use the composite index map to the coarser cover.

For a LAB$'$ cylinder, choose a good refinement whose restrictions to
the two ends are good covers, using a triangulation with the ends as
subcomplexes. Choosing singleton cores on this cylinder produces a
normalized cylinder. The preceding comparison on each end identifies
its restriction with the normalization obtained from any good
refinement of that endpoint cover. This proves compatibility with
concordance.
We have obtained an inverse $[\mathfrak A]\mapsto[\mathfrak A_B]$ to
the inclusion on classes. On a normalized diagram take $B_i=A_i$;
in the other direction $\mathfrak A_B\to\mathfrak A$ is a LAB$'$ morphism.
\end{proof}

\begin{corollary}[Classification of normalized strict LABs]\label{cor:strict-main}
The natural maps
\[
 \mathsf F(X)\xrightarrow{\ \bar\Phi\ }[X,|N\C|]
 \xrightarrow{\ \rho_*\ }[X,B\Aut(D)]
\]
are bijections, and their composite is categorical rectification.
\end{corollary}
\begin{proof}
Combine Theorems~\ref{thm:main} and~\ref{thm:strict-prime},
Proposition~\ref{prop:normalization}, and~\eqref{eq:rect-compat}.
\end{proof}

The use of concordance in Proposition~\ref{prop:normalization} and
Corollary~\ref{cor:ambient-unique} is essential to this argument.
It does not prove the corresponding assertion for the relation generated
by elementary equivalences and refinements alone. It also does not prove
Brown's wedge or Mayer--Vietoris axioms on arbitrary infinite CW complexes.

One can define $\mathsf F$ on finite CW complexes using all finite open
covers and the same generated relation. Pullback makes it homotopy
invariant. Every finite CW complex has the homotopy type of a finite
polyhedron. Transporting the displayed natural classification along such
an equivalence extends it to finite CW complexes; independence of the
chosen model follows from naturality and homotopy invariance. This is a
finite-complex extension, not an infinite-CW representability argument.

\section{Completed coherent descent and the raw comparison}\label{sec:coherent}

We first compute descent after groupoid completion of $\C$ on an arbitrary
indexed open cover. We then take cover diagrams valued in a Reedy fibrant
replacement of $\Sing N\C$. On finite CW complexes the diagonal
identifies their concordance classes, with specified refinements, with
the completed classes of Definition~\ref{def:coherent}.

\subsection{Descent after groupoid completion}

Put $Z=\Sing N\C$, with the singular functor applied degreewise, and
choose a Kan replacement
\[
 \diag Z\longrightarrow K_{k,l}.
\]
Here $K_{k,l}$ is a simplicial set. Its realization has the homotopy
type of $|N\C|$, and hence of $BG$ by Theorem~\ref{thm:main}.
It models the \emph{groupoid completion} of the matrix category:
all its arrows are inverted in the resulting homotopy type. This is
different from Reedy fibrant replacement of a simplicial space, and
from Rezk completion, which does not invert arbitrary noninvertible
embeddings.

For an indexed open cover $\mathcal U=\{U_i\}_{i\in I}$,
let $\check C_\bullet(\mathcal U)$ be its ordinary \v{C}ech simplicial
space, with
\[
 \check C_p(\mathcal U)
 =\coprod_{(i_0,\ldots,i_p)\in I^{p+1}}U_{i_0\cdots i_p},
 \qquad U_{i_0\cdots i_p}=\bigcap_{a=0}^p U_{i_a}.
\]
Repeated indices are allowed. Faces forget an index and degeneracies
repeat one. The descent formula uses all ordered tuples of indices.
The strict construction in Section~\ref{sec:strict} uses chains of
nonempty index sets in the intersection category $\mathcal U_\cap$.
On compact bases finite covers suffice for the construction below.

\begin{definition}[Completed coherent LAB classes]\label{def:coherent}
The space of completed coherent descent data on $\mathcal U$ is
\begin{equation}\label{eq:completed-descent}
 \Desc_{k,l}(\mathcal U)
 :=\holim_{[p]\in\Delta}\;
 \prod_{i_0,\ldots,i_p}
 \mathbf R\Map_{\mathrm{sSet}}
       (\Sing U_{i_0\cdots i_p},K_{k,l}).
\end{equation}
Since every simplicial set is cofibrant and $K_{k,l}$ is Kan, each
$\mathbf R\Map$ here is computed by the ordinary simplicial mapping
space $\Map$. Empty intersections contribute a point. We use the
totalization model for the homotopy limit, as described below.

A completed coherent class on $X$ is represented by a pair
$(\mathcal U,z)$ with $z\in\Desc_{k,l}(\mathcal U)$, modulo the
relation generated by paths, pullback along specified index-map
refinements, and the two restrictions of descent data on $X\times I$.
Denote the set of classes by $\HCgpd(X)$. The superscript is part of
the notation and records the groupoid completion.
\end{definition}

The homotopy limit in~\eqref{eq:completed-descent} records local maps,
homotopies between their restrictions, and all higher compatibilities.
The augmentation to $X$ supplies a natural comparison
\begin{equation}\label{eq:descent-comparison}
 \epsilon_{\mathcal U}^*:
 \mathbf R\Map(\Sing X,K_{k,l})
 \longrightarrow\Desc_{k,l}(\mathcal U).
\end{equation}

\begin{theorem}[Completed coherent classification]\label{thm:coherent}
For every indexed open cover $\mathcal U$ of $X$, the comparison
map~\eqref{eq:descent-comparison} is a weak equivalence. In particular, for finite CW complexes there
are natural bijections
\[
 \HCgpd(X)\cong[X,|N\C|]\cong[X,B\Aut(D)].
\]
The local descent spaces all model
$\mathbf R\Map(\Sing X,K_{k,l})$; their comparison is stronger than
the stated bijection of classes.
\end{theorem}
\begin{proof}
Put $Q_p=\Sing\check C_p(\mathcal U)$ and $K=K_{k,l}$.
The simplex $\Delta^q_{\rm top}$ is connected, so $Q_p$ is the
coproduct of the singular complexes of the intersections
$U_{i_0\cdots i_p}$. The cosimplicial mapping space in
\eqref{eq:completed-descent} is thus $p\mapsto\Map(Q_p,K)$.
The latching inclusions of $Q_\bullet$ are monomorphisms.
Consequently the matching maps of this cosimplicial mapping space are
fibrations: they are restrictions along these inclusions into the Kan
complex $K$. Its totalization therefore computes the homotopy limit
\cite[Part~II]{BK}.

The defining adjunction for totalization gives a natural isomorphism
\[
 \operatorname{Tot}\Map(Q_\bullet,K)
   \cong\Map\!\left(\int^{[p]} Q_p\times\Delta[p],K\right)
   \cong\Map(\diag Q,K).
\]
In simplicial degree $q$, the coend identifies with $(Q_q)_q$ by the
Yoneda lemma. By Dugger and Isaksen's \v{C}ech augmentation theorem
\cite[Theorem~2.1]{DI}, the map
$|\check C_\bullet(\mathcal U)|\to X$ is a weak equivalence for every
indexed open cover. The degeneracies of $\check C_\bullet(\mathcal U)$
include whole coproduct components, so the realization
comparison in Section~\ref{sec:Ex} identifies the augmentation
$\diag Q\to\Sing X$ as a weak equivalence. Mapping into $K$ gives
the required weak equivalence~\eqref{eq:descent-comparison}.
It commutes with every specified refinement and pullback.

Every cover refines the one-element cover $\{X\}$. Thus paths and refinement already identify the set of classes with
$\pi_0\mathbf R\Map(\Sing X,K_{k,l})$.
The projection $X\times I\to X$ is a homotopy equivalence, and its
two sections induce the same map on these components. Adding the
concordance relation therefore imposes no further identifications.
For finite CW $X$, realization and singular adjunction identify the
components with $[X,|K_{k,l}|]$. Finally use
$|K_{k,l}|\simeq|N\C|$ and Theorem~\ref{thm:main}.
\end{proof}

Changing the Kan replacement gives a naturally equivalent descent
theory, since the mapping spaces and homotopy limits are derived.
There is no need to assume that two refinement choices act by a path
on a raw Reedy mapping space. The comparison here takes place after
groupoid completion, where the \v{C}ech augmentation is an equivalence.

\begin{corollary}[Strict representatives of completed classes]\label{cor:strict-coherent}
On finite CW complexes, a normalized strict LAB determines a completed
coherent class, and this construction induces a natural bijection
\[
 \mathsf F(X)\xrightarrow{\cong}\HCgpd(X).
\]
Its composite with classification by $[X,BG]$ is categorical
rectification.
\end{corollary}
\begin{proof}
First let $X$ be a finite polyhedron. On a good refinement, choose
trivializations and write the strict diagram as a continuous functor
$a:\mathcal U_\cap\to\C$. Applying the singular nerve and the
diagonal gives
\[
 \diag\Sing N\mathcal U_\cap\longrightarrow\diag Z
 \longrightarrow K_{k,l}.
\]
The augmentation of its source to $\Sing X$ is a weak equivalence.
The diagram therefore specifies a unique morphism
$\Sing X\to K_{k,l}$ in the homotopy category, and hence a completed
class by Theorem~\ref{thm:coherent}. Under its classification map this
is exactly $\bar\Phi(\mathfrak A)$ from Section~\ref{sec:strict}.
That section proves independence of choices and of all generating
strict equivalences. Corollary~\ref{cor:strict-main} proves bijectivity
and compatibility with $\bar R$. The finite-CW extension follows from
the natural homotopy invariance established in
Section~\ref{sec:normalization}.
\end{proof}

The corollary concerns classes. It does not assert equivalence of the
unlocalized spaces of strict diagrams and completed descent data.
It also does not turn a raw coherent datum on a fixed cover into an
isomorphic strict datum on that cover.

\subsection{Raw Segal cover diagrams and concordance}\label{subsec:raw-coherent}

Let $Z\to\mathbf C$ be a Reedy fibrant replacement of $Z=\Sing N\C$.
The nerve has discrete object space and satisfies the homotopy Segal
condition; hence $\mathbf C$ is a Segal space. Its diagonal has the
same weak homotopy type as $K_{k,l}$. For a fixed cover define
\[
 \mathcal D_{\mathbf C}(\mathcal U)
 =\mathbf R\Map_{\mathrm{Reedy}}
       (\Sing N\mathcal U_\cap,\mathbf C).
\]
The source is Reedy cofibrant (its latching maps are monomorphisms)
and the target is Reedy fibrant, so the ordinary enriched mapping space
computes the displayed $\mathbf R\Map$. In this simplicial enrichment, its mapping-space paths are homotopies in the inner simplicial
direction. They must not be confused with noninvertible natural
transformations in the categorical direction.

\begin{definition}[An ordinary-refinement raw model]\label{def:raw-coherent}
For compact $X$, let $\operatorname{Cov}_{\mathrm{disc}}(X)$ be the
ordinary category of finite indexed open covers and specified
index-map refinements. A morphism $\mathcal U\to\mathcal V$ is a
function $\alpha:I\to J$ with $U_i\subseteq V_{\alpha(i)}$.
Repeated members are allowed. Pullback retains every index, including
those whose inverse image is empty. Index sets may be chosen to be
finite ordinals, making this a small category. All morphism sets are
discrete; parallel refinement choices are not identified.

Take the disjoint union of the sets
$\pi_0\mathcal D_{\mathbf C}(\mathcal U)$ over finite indexed covers.
Identify each component with its pullback along each specified refinement,
and identify the two restrictions of a datum on $X\times I$. The resulting
quotient is $\HCraw(X)$. Equivalently, if
\[
 \mathcal L^{\mathrm{raw}}_{\C}(X)
 =\hocolim_{\operatorname{Cov}_{\mathrm{disc}}(X)^{\op}}
                   \mathcal D_{\mathbf C}(\mathcal U),
\]
then $\HCraw(X)=\pi_0\mathcal L^{\mathrm{raw}}_{\C}(X)/{\sim_c}$,
where $\sim_c$ denotes endpoint concordance. Indeed $\pi_0$ of this
homotopy colimit is the ordinary colimit of the sets of components:
the edges of the simplicial replacement impose exactly the refinement
identifications. Different index maps can induce the same restriction
map. For instance, two index maps that differ only on an empty patch
induce the same functor on intersection categories.
In Example~\ref{ex:raw-refinement}, two index maps send the same
datum to different components of the fixed-cover mapping space.
\end{definition}

Pullback along continuous maps and along homotopies is defined using
the corresponding pullback covers. Thus this is a well-defined
concordance-class functor. Functorial derived mapping spaces also make
it independent, up to natural equivalence, of the chosen Reedy
replacement. This definition is a specified raw model; no equivalence
with another enriched cover construction is assumed.

\paragraph{Raw diagrams on vertex stars.}
More generally, let $W$ be any Reedy fibrant simplicial space
in simplicial sets. Write $\mathsf{Raw}_W(X)$ for the quotient of
Definition~\ref{def:raw-coherent} formed using $W$ as the target.
For $W=\mathbf C$ this is $\HCraw(X)$. Let $V=W_{\bullet,0}$
be the simplicial set of vertices in the inner direction. For a simplicial set $J$
write $cJ$ for the bisimplicial set constant in the inner direction.
A map $cJ\to W$ is determined by a map $J\to V$. The degeneracies in
the inner direction induce a weak equivalence
\begin{equation}\label{eq:raw-vertex-diagonal}
 V=\diag(cV)\longrightarrow\diag W .
\end{equation}
Indeed, in each inner degree $q$, evaluation at the zeroth vertex
$W_{\bullet,q}\to V$ is a trivial Kan fibration of outer simplicial
sets. A lifting problem against $\partial\Delta[p]\hookrightarrow
\Delta[p]$ asks for a lift against $\{0\}\hookrightarrow\Delta[q]$
for the matching fibration $W_p\to M_pW$. Such a lift exists since the
vertex inclusion is a trivial cofibration. The row map of $cV\to W$ is
a section of this trivial fibration. The diagonal lemma for bisimplicial
sets now gives~\eqref{eq:raw-vertex-diagonal}.

Let $K$ be a finite ordered triangulation, $\mathcal U(K)$ its cover by
vertex open stars, and $P(K)$ its poset of nonempty faces. There is a
degreewise weak equivalence
\begin{equation}\label{eq:raw-star-forget}
 q_K:\Sing N\mathcal U(K)_\cap\longrightarrow c\,N P(K),
 \qquad N P(K)=\Sd K.
\end{equation}
In outer degree $p$ this map forgets the point in the intersection
corresponding to a chain $S_0\subseteq\cdots\subseteq S_p$. Its fiber
over that chain is the singular complex of the open star of $S_p$,
which contracts to the barycenter of $S_p$. Both sources are Reedy cofibrant. Mapping~\eqref{eq:raw-star-forget}
into $W$ therefore gives a weak equivalence of mapping spaces.
In particular, every raw datum on this cover is joined by a path to
one represented by a simplicial map $N P(K)\to V$.

Under the usual homeomorphism $|N P(K)|\cong |K|$, the realization
of the forgetful map agrees up to homotopy with the augmentation to
$|K|$. In the topological cover-poset nerve, a weighted chain maps
to the same weighted combination of the barycenters of $S_0,\ldots,S_p$.
This point and the augmentation $x$ lie in the simplex supporting $x$,
which contains $S_p$. Straight-line interpolation gives the homotopy.
The degeneracies of the cover-poset nerve include whole coproduct
components, so the comparison with derived realization carries this
homotopy to the identification used in~\eqref{eq:raw-star-forget}.

\paragraph{Realizing an outer homotopy on a cylinder.}
Suppose $T$ is a finite ordered complex and two maps $N P(T)\to V$
are the ends of a simplicial map $N P(T)\times\Delta[1]\to V$.
The corresponding raw diagrams on $\mathcal U(T)$ are concordant.
Cover $|T|\times I$ by the open sets
\[
 \operatorname{st}_T(v)\times J_e,\qquad
 J_0=[0,2/3),\quad J_1=(1/3,1],\quad e=0,1.
\]
On an intersection indexed by a nonempty set $S$ of pairs $(v,e)$,
send its cover-poset object to
\[
 \bigl(\{v:(v,e)\in S\ \text{for some }e\},
                 \max\{e:(v,e)\in S\}\bigr)\in P(T)\times[1].
\]
The first coordinate is a face because the indicated stars intersect.
Both coordinates preserve inclusions. This continuous functor induces
a map of singular nerves to $c(N P(T)\times\Delta[1])$. Composition
with the given homotopy gives a raw diagram on the cylinder. Its
endpoint restrictions are the two original diagrams after removal of
empty patches by the inclusion of the active indices, a specified
refinement in Definition~\ref{def:raw-coherent}.

\begin{proposition}[Raw concordance and realization]\label{prop:raw-comparison}
For a Reedy fibrant simplicial space $W$ in simplicial sets and a
finite CW complex $X$, the diagonal induces a natural bijection
\[
             \mathsf{Raw}_W(X)\cong[X,|\diag W|].
\]
With $W=\mathbf C$ this gives
$\chi_X:\HCraw(X)\xrightarrow{\cong}\HCgpd(X)$.
\end{proposition}
\begin{proof}
A raw map produces a diagonal map from $\diag\Sing N\mathcal U_\cap$
to $\diag W$. The covering theorem identifies its source with
$\Sing X$ up to weak equivalence. A mapping-space path gives a
homotopy on diagonals, and specified refinements commute with the
cover augmentations. Pullback over $X\times I$ gives concordance
invariance. These facts define the natural class map in the statement.

For surjectivity take a finite ordered triangulation $K$ of a
polyhedral model of $X$. By~\eqref{eq:raw-vertex-diagonal} a class
is represented by $K\to\Ex^\infty V$. This map factors through
some $\Ex^r V$, since $K$ has finitely many nondegenerate simplices.
Let $h:T=\Sd^rK\to V$ be its adjoint. On the vertex stars of $T$
use the diagram represented by $\Sd T\to T\xrightarrow{h} V$,
with the first arrow the last-vertex map. Finite subdivision and the
last-vertex maps give the canonical identifications in the homotopy
category, so this diagram has the prescribed class.

\smallskip\noindent\emph{Injectivity on a finite polyhedron.}
Choose two raw classes with equal images and representatives on finite
covers. This choice is possible since every component of the homotopy
colimit comes from one of its terms. Subdivide a triangulation of $X$
until its vertex-star cover refines both covers, and choose the two
index maps. Call the resulting ordered triangulation $K$ and put
$L=N P(K)=\Sd K$. By~\eqref{eq:raw-star-forget} and paths on this
common cover, the two restrictions have representatives $f_0,f_1:L\to V$.
The domain identification following~\eqref{eq:raw-star-forget} and
the target equivalence~\eqref{eq:raw-vertex-diagonal} identify
their images under the diagonal class map with the usual homotopy classes of
$f_0,f_1$. Thus $[f_0]=[f_1]$ in the Kan--Quillen homotopy category.

For a simplicial set $J$, denote the last-vertex map by
$d_J:\Sd J\to J$ and its $a$-fold iterate by $d_J^{(a)}:\Sd^aJ\to J$,
with $d_J^{(0)}=\id$. Use Kan's fibrant replacement $V\to\Ex^\infty V$
\cite[Sections~3.3.4--3.3.6]{Kerodon}. The mapping space from $L$ into
$\Ex^\infty V$ is Kan, so its two vertices given by $f_0,f_1$ are joined
by an edge. The resulting map $L\times\Delta[1]\to\Ex^\infty V$ factors
through some $\Ex^r V$ because the prism has finitely many nondegenerate
simplices. All transition maps into $\Ex^\infty V$ are monomorphisms,
so the endpoints at this stage are exactly the images of $f_i$.
The subdivision adjunction gives
\[
 h:T=\Sd^r(L\times\Delta[1])\longrightarrow V,\qquad
 h|_{T_i}=f_i\,d_L^{(r)},\qquad
 T_i=\Sd^r(L\times\{i\})\cong\Sd^r L .
\]
The endpoint formula follows because the natural map to $\Ex$ is
adjoint to the last-vertex map. The ordered prism triangulates $X\times I$
with the ends as subcomplexes. On the vertex stars of its subdivision
$T$, take the raw diagram represented by $h\,d_T:N P(T)\to V$.
At the $i$th end the nonempty restricted patches are the vertex stars
of $T_i$. By naturality of the last-vertex map the endpoint diagram
is represented by $f_i\,d_L^{(r+1)}$.

\smallskip\noindent\emph{Comparison with the specified refinements.}
Set $a=r+1$, so $T_i=\Sd^a K$. The vertex map of
$b=d_K^{(a)}:T_i\to K$ gives an index-map refinement
$\mathcal U(T_i)\to\mathcal U(K)$. Indeed each selected vertex $b(w)$
lies in the smallest old face containing $w$. On the open star of
$w$ its old barycentric coordinate is positive. On the new cover the
specified pullback of $f_i$ is represented by $f_i\,N P(b)$, where
$P(b)(\tau)$ consists of the distinct images of the vertices of $\tau$.

Let $c:P(T_i)\to P(K)$ assign to a face $\tau$ its carrier, the
smallest face of $K$ containing $|\tau|$. In one barycentric subdivision
the carrier of a chain of faces is its largest face. Iterating this
observation shows that $Nc=d_L^{(a)}:N P(T_i)\to N P(K)$.
Thus the prism endpoint is $f_i\,Nc$. For each face $\tau$ of $T_i$,
the selected vertices satisfy
\[
                     P(b)(\tau)\subseteq c(\tau).
\]
These inclusions give a natural transformation $P(b)\Rightarrow c$.
Its nerve, followed by $f_i$, is an outer simplicial homotopy
$N P(T_i)\times\Delta[1]\to V$. The two-interval cover constructed
above makes this a raw concordance from the specified pullback to the
prism endpoint. Applying this at both ends shows that the original two
classes coincide. Both class functors are homotopy invariant, so a
finite polyhedral homotopy model gives bijectivity on every finite
CW complex. When $W=\mathbf C$, the zigzag
$\diag\mathbf C\leftarrow\diag Z\to K_{k,l}$ and
Theorem~\ref{thm:coherent} identify the target with $\HCgpd(X)$.
\end{proof}

For a strict LAB on a finite polyhedron, the map of its trivialized
nerve into $\mathbf C$ gives a raw representative. By construction its image
under $\chi_X$ is the class of Corollary~\ref{cor:strict-coherent}.
Injectivity now gives independence of all strict choices and generating
relations. The resulting natural bijection $\mathsf F(X)\cong\HCraw(X)$
extends to finite CW complexes by homotopy invariance. This proves the
chain in Theorem~\ref{thm:intro-classification}.

\begin{example}[Refinement choices before realization]\label{ex:raw-refinement}
Let $X=\{*\}$, $\mathcal U=\{*\}$, and
$\mathcal V=\{*_0,*_1\}$. There are two refinements
$\alpha_0,\alpha_1:\mathcal U\to\mathcal V$, and
$\mathcal V_\cap$ is the poset
\[
 \{0\}\longrightarrow\{0,1\}\longleftarrow\{1\}.
\]
Define a strict diagram by sending these three objects to $D_0,D_1,D_1$,
respectively, with the standard embedding and the identity as arrows.
It defines a point $a$ of $\mathcal D_{\mathbf C}(\mathcal V)$.
For $\mathcal U$, the source is the outer representable $\Delta[0]$,
so $\mathcal D_{\mathbf C}(\mathcal U)\simeq\mathbf C_0$.
Reedy replacement is degreewise a weak equivalence, and therefore
\[
 \pi_0\mathbf C_0\cong\pi_0Z_0\cong\mathbb N_0.
\]
The two restrictions of $a$ belong to components $0$ and $1$.
Thus the two restriction maps cannot be homotopic as maps into this
fixed-cover Reedy mapping space.

Normalization does not remove the obstruction. Take both covers to
have two copies of the point, send both singleton objects to $D_0$
and the double intersection to $D_1$, and compare the identity index
map with the constant index map of value $0$. Evaluation at the
double-intersection object gives ranks $D_1$ and $D_0$ for the two
normalized pullbacks, again excluding a Reedy mapping-space path.
\end{example}

The two restrictions in Example~\ref{ex:raw-refinement} are connected
through the original datum in the homotopy colimit over covers.
The cylinder in Proposition~\ref{prop:raw-comparison} likewise has its
own cover. These changes of cover explain how the concordance
classification accommodates the distinct fixed-cover components in
the example. Berwick-Evans, Boavida de Brito and Pavlov prove a
mapping-space result on smooth manifolds for their enriched cover
construction \cite[Definition~6.4 and Theorem~6.6]{BEBP}.
Comparing the underlying spaces of that model and
Definition~\ref{def:raw-coherent} is a further question.

\section{Frame spaces, homotopy groups, and global bundles}\label{sec:frames}

Put
\[
 \operatorname{Fr}_{k,l^N}=\Emb(M_k,M_{kl^N}),\qquad
 \operatorname{Fr}_{k,l^\infty}
 =\operatorname*{colim}_N\operatorname{Fr}_{k,l^N}.
\]
Corollary~\ref{cor:I4} and Lemma~\ref{lem:restriction} give, independently
of any LAB classification,
\begin{equation}\label{eq:frame-G}
 \operatorname{Fr}_{k,l^\infty}\simeq\Emb(M_k,D)\simeq G.
\end{equation}
An inverse to the restriction equivalence specifies a map $\beta$ in
the homotopy category. Thomsen computed the homotopy groups of the
automorphism group of a UHF algebra \cite{Thomsen}. For the present
type the frame description gives the calculation below and identifies
the finite-stage maps used in the examples.

\begin{proposition}\label{prop:homotopy}
The group $G$ is path connected, and for $j\geq1$,
\[
 \pi_{2j-1}G\cong\mathbb Z/k,\qquad \pi_{2j}G=0.
\]
Consequently $BG$ is simply connected, with
$\pi_{2j}BG\cong\mathbb Z/k$ and $\pi_{2j+1}BG=0$.
\end{proposition}
\begin{proof}
Use the homogeneous-space fibration
\[
 \mathrm{PU}(l^N)\longrightarrow\mathrm{PU}(kl^N)
 \longrightarrow\operatorname{Fr}_{k,l^N}.
\]
On fundamental groups the first map is
$\mathbb Z/l^N\to\mathbb Z/(kl^N)$, $a\mapsto ka$.
Its kernel is zero and its cokernel is $\mathbb Z/k$, giving
$\pi_2\operatorname{Fr}_{k,l^N}=0$ and
$\pi_1\operatorname{Fr}_{k,l^N}=\mathbb Z/k$.
For a fixed degree $d\geq3$, choose $N$ with $d<2l^N$.
The fibration $U(t)\to U(t+1)\to S^{2t+1}$ puts both unitary groups
in the stable range through degree $d$. Passage to the projective
unitary group preserves homotopy groups in degrees at least three,
and both projective unitary groups have vanishing $\pi_2$.
Bott periodicity~\cite{Bott} gives
$\mathbb Z$ in odd degrees and zero in positive even degrees, and the
tensor multiplicity-$k$ map induces multiplication by $k$ on the odd
groups. The long exact sequence gives $\mathbb Z/k$ in odd degrees and
zero in even degrees for the frame space. The stabilization maps multiply
these cyclic groups by $l$. Since $(k,l)=1$, they are isomorphisms.
Pass to the colimit using the cofibration/telescope comparison and then
use~\eqref{eq:frame-G}. The frame spaces are connected, as is $G$.
\end{proof}

Thus the sphere formula
\begin{equation}\label{eq:spheres}
 \mathsf F(S^n)\cong\HCraw(S^n)\cong\HCgpd(S^n)\cong
 \begin{cases}\mathbb Z/k,&n\geq2\text{ even},\\
 *,&n\geq1\text{ odd}
 \end{cases}
\end{equation}
holds for the three concordance-class functors by
Theorem~\ref{thm:intro-classification}. For the simply connected
classifying space, based and unbased sphere classes agree. Formula~\eqref{eq:spheres}
must not be substituted for a computation of a fixed-rank global image.

\paragraph{Fixed-rank and arbitrary finite-rank global classes.}
Using $D_r\otimes A\cong D$, stabilization induces a map
\[
 j_r:B\mathrm{PU}(kl^r)\longrightarrow BG.
\]
On bundle classes over $X$, it sends an $M_{kl^r}$-bundle $B$ to the
$D$-bundle $B\otimes A$. For a finite complex $X$ define the global
fixed-rank image to be
$\operatorname{im}([X,B\mathrm{PU}(k)]\to[X,BG])$.
The arbitrary finite-rank global image is the union of the images of
$j_r$, or equivalently the image of
$[X,\hocolim_r B\mathrm{PU}(kl^r)]$. These are definitions inside the
classified theory; neither is the rigid notion of globality unless an
additional rigid comparison has been proved.

\begin{example}[A fixed-rank counterexample]\label{ex:fixedrank}
Take $k=2$, $l=3$, and $X=S^6$.
Since $\pi_5\mathrm{PU}(2)=\pi_5 S^3=\mathbb Z/2$,
whereas $\pi_5\mathrm{PU}(6)=\mathbb Z$, the map
$\mathrm{PU}(2)\to\mathrm{PU}(6)$ is zero on $\pi_5$.
The stabilization $j_0$ factors through $B\mathrm{PU}(6)$; hence
\[
 \pi_6 B\mathrm{PU}(2)\longrightarrow\pi_6 B\Aut(M_{2\cdot3^\infty})
   =\mathbb Z/2
\]
is zero. The full theory has two classes, but the fixed-rank global image
has one. The classical value $\pi_5S^3$ may be found in~\cite{Toda}.
Here the source group $\pi_5\mathrm{PU}(2)$ lies outside the stable range.
\end{example}

If the matrix rank $kl^r$ is allowed to increase with the dimension,
all sphere classes do lie in a global image. For $S^{2j}$ with
$j\geq2$, choose $r$ with $2j-1<2l^r$. Restriction of
$\mathrm{PU}(kl^r)$ to the standard $M_k$ gives the quotient
$\mathbb Z\to\mathbb Z/k$ on $\pi_{2j-1}$ in the homogeneous-space
calculation above. Every subsequent frame stabilization is an
isomorphism in this degree. The degree-two assertion follows from the analogous
quotient of fundamental groups. This statement about spheres does not
imply that all classes on every finite complex are global.

\begin{proposition}[A class outside every finite-rank global image]\label{prop:nonglobal}
For $k>1$, let $Y=\operatorname{Fr}_{k,l}$ and $X=\Sigma Y$.
There is a class in $\mathsf F(X)$, and hence in both coherent theories,
which is not the stabilization of any $M_{kl^r}$-bundle, for any $r$.
\end{proposition}
\begin{proof}
The inclusion $Y\to\Emb(M_k,D)$ lifts through
$I_0\to\Emb(M_k,D)$: its pullback is a principal contractible-$H$
bundle on the compact manifold $Y$ and is trivial. Identifying $I_0$
with $G$ gives $\beta_Y:Y\to G$. Clutch a $D$-bundle on $\Sigma Y$
using this map. The same homogeneous-space sequence, already in its
$\pi_3$ stable range since $l\geq2$, gives $\pi_3Y=\mathbb Z/k$.
Stabilization to the infinite frame space is an isomorphism on $\pi_3$,
so $(\beta_Y)_*:\pi_3Y\to\pi_3G$ is an isomorphism.

If this bundle came from an $M_{kl^r}$-bundle, trivializing that bundle
on the two cones would produce a clutching map
$c:Y\to\mathrm{PU}(kl^r)$ whose stabilization is homotopic to
$\beta_Y$. But $\pi_3\mathrm{PU}(kl^r)=\mathbb Z$, and every homomorphism
$\mathbb Z/k\to\mathbb Z$ is zero. This contradicts the preceding
isomorphism. Indeed, changing cone trivializations multiplies the
clutching map on the left and right by maps extending over the
cones. Each such map contracts to a constant, and those constants
can be joined to the identity in the path-connected group $G$.
To compare induced \emph{based} homomorphisms, fix $y_0\in Y$.
For a free homotopy $H_t:Y\to G$ the formula
$y\mapsto H_t(y_0)^{-1}H_t(y)$ gives a based homotopy of the
normalized maps. Equivalently, the usual action of $\pi_1(G)$ on
$\pi_3(G)$ is trivial for a topological group. The suspension
has finite CW type, so the class belongs to the classified theories.
\end{proof}

\paragraph{Clutching over a suspension.}
For a connected pointed finite CW complex $Y$, clutching gives
\[
 \HCgpd(\Sigma Y)\cong
 [Y,G]\cong[Y,\operatorname{Fr}_{k,l^\infty}],
\]
and the same holds for $\HCraw$ and $\mathsf F$. A frame map is an embedding,
not an automorphism of a fixed matrix algebra; it does not by itself
give a clutching map for an ordinary matrix algebra bundle.
The proposition shows this distinction even after allowing arbitrary
finite ranks $kl^r$.

\subsection{Stable embeddability and the Grassmannian fiber}

The forgetful map from the stabilized Grassmannian of
Section~\ref{sec:prelim} classifies the tautological $M_k$-bundle.
Its composite with stabilization is null-homotopic, and the resulting
sequence is a homotopy-fiber sequence.

\begin{proposition}\label{prop:grassmannian-fiber}
Let $j:\mathrm{PU}(k)\to G$ send $\alpha$ to
$\alpha\otimes\id_A$, using $D=M_k\otimes A$. There is a homotopy-fiber
sequence
\[
 \operatorname{Gr}_{k,l^\infty}\xrightarrow{\gamma_k}
 B\mathrm{PU}(k)\xrightarrow{Bj}BG.
\]
The first map forgets the stabilized embedding.
\end{proposition}
\begin{proof}
Put $P=\mathrm{PU}(k)$, $E_0=\Emb(M_k,D)$ and
$T=\operatorname{Tel}_N\operatorname{Fr}_{k,l^N}$. Restriction gives
$P$-equivariant weak equivalences $T\to E_0\leftarrow G$ by
Corollary~\ref{cor:I4} and Lemma~\ref{lem:restriction}. The right
actions are precomposition on frames and right multiplication by $j(P)$
on $G$. For a right $P$-space $V$ form the ordinary associated bundle
$(EP\times V)/P\to BP$ using the diagonal right action. A
$P$-equivariant weak equivalence gives a weak equivalence of these
bundles by their homotopy exact sequences.

For $V=T$, the total space is the telescope of the corresponding bundles
for the finite frame spaces. At each stage forgetting $EP$ is a bundle
over $\operatorname{Gr}_{k,l^N}$ with contractible fiber $EP$, by
local triviality of the principal frame bundle. Its projection to $BP$
classifies that frame bundle: the two pullbacks of the principal bundles
are identified by $EP\times\operatorname{Fr}_{k,l^N}$. These maps are
compatible with stabilization, so the telescope projection to $BP$
represents $\gamma_k$.

For $V=G$, inversion in the second coordinate identifies the associated
bundle with the principal $G$-bundle obtained from $EP$ by extension
of structure group along $j$. It is classified by $Bj$; its total
space thus models the homotopy fiber of $Bj$, by pullback of
$EG\to BG$. This gives the asserted sequence over $BP$.
\end{proof}

The fiber sequence and its lifting criterion use only restriction
to $M_k$ and finite-stage approximation, and hold for all $k\ge1$
and $l\ge2$. We return to the coprime case for the conclusion about
the Brauer class below.

\begin{corollary}\label{cor:stable-embeddability}
For an $M_k$-bundle $B$ on a finite CW complex $X$, the following
conditions are equivalent:
\begin{enumerate}[label=(\roman*),leftmargin=*]
 \item $B$ embeds unitally in $X\times M_{kl^N}$ for some $N$;
 \item the $D$-bundle $B\otimes A$ is trivial;
 \item its class in the strict LAB concordance theory $\mathsf F(X)$
 is zero.
\end{enumerate}
\end{corollary}
\begin{proof}
If $f_B:X\to B\mathrm{PU}(k)$ classifies $B$, condition (ii) says that
$Bj\circ f_B$ is null-homotopic. By
Proposition~\ref{prop:grassmannian-fiber}, this is equivalent to a
homotopy lift of $f_B$ to $\operatorname{Gr}_{k,l^\infty}$.
Compactness of $X$ and the telescope model move the lift to a finite
stage. Proposition~\ref{prop:embedded-classification} then gives (i),
using the isomorphism of bundles supplied by the homotopy of their
classifying maps. The converse follows by composing the finite-stage
lift into the homotopy fiber. Equivalence with (iii) is the compatibility
of categorical rectification with the classified theory.
The corresponding raw and completed coherent classes are zero as well,
by Theorem~\ref{thm:intro-classification}.
These conditions also imply Morita triviality. An embedding as in (i)
has an $M_{l^N}$-centralizer $Z$ with $B\otimes Z$ trivial, so
$\delta(B)=-\delta(Z)$ is annihilated by both $k$ and $l^N$.
Coprimality forces $\delta(B)=0$.
\end{proof}

\begin{example}[Morita trivial but nontrivial as a LAB]\label{ex:morita-versus-lab}
For $k>1$, take the $M_k$-bundle on $S^2$ classified by a generator of
$\pi_2B\mathrm{PU}(k)=\mathbb Z/k$. Its Brauer class is zero because
$H^3(S^2;\mathbb Z)=0$, so it is an endomorphism bundle. Its image
under $Bj$ is nonzero. Indeed the map
$\mathrm{PU}(k)\to\mathrm{PU}(kl^N)$ multiplies fundamental groups
by $l^N$, and projection to
$\operatorname{Fr}_{k,l^N}$ reduces modulo $k$.
The composite is multiplication by $l^N$ on $\mathbb Z/k$,
an isomorphism since $(k,l)=1$. By~\eqref{eq:frame-G}, $j$ therefore induces an isomorphism
on $\pi_1$, and $Bj$ on $\pi_2$.
Thus the bundle has nonzero LAB class and is not embeddable in any
allowed trivial stage, despite being Morita trivial.
\end{example}

\begin{corollary}[Tensor complements]\label{cor:tensor-complement}
For an $M_s$-bundle $B$ on a finite CW complex, with $s\geq2$, there
is an $M_{s^n}$-bundle $C$ for some $n\geq0$ such that
$B\otimes C\cong X\times M_{s^{n+1}}$.
\end{corollary}
\begin{proof}
Take $k=l=s$ in the lifting criterion of
Proposition~\ref{prop:grassmannian-fiber}, which holds without a
coprimality assumption. Let $f_B:X\to B\mathrm{PU}(s)$ classify $B$.
The composite $Bj\circ f_B$ classifies $B\otimes M_{s^\infty}$.
This bundle has fiber isomorphic to $M_{s^\infty}$ and is trivial by
Dadarlat--Pennig's theorem. A null-homotopy of $Bj\circ f_B$ gives a
homotopy lift of $f_B$ to the Grassmannian telescope. Since $X$ is
compact, the lift can be moved into a finite stage. This supplies an
embedding into the trivial $M_{s^{n+1}}$-bundle
for some $n\ge0$. Its centralizer is the required $M_{s^n}$-bundle
by Proposition~\ref{prop:finite-centralizer}.
\end{proof}

\section{Changing the auxiliary integer and the matrix type}\label{sec:limits}

\begin{proposition}[Independence of the auxiliary integer]
Let $k\geq1$ and $l,l'\geq2$ satisfy $\gcd(k,l)=\gcd(k,l')=1$.
Then $B\Aut(M_k\otimes M_{l^\infty})$ and
$B\Aut(M_k\otimes M_{(l')^\infty})$ are connected by natural
stabilization weak equivalences through
$B\Aut(M_k\otimes M_{(ll')^\infty})$.
\end{proposition}
\begin{proof}
Tensor an automorphism with $\id_{M_{(l')^\infty}}$ and use the
identification $M_{l^\infty}\otimes M_{(l')^\infty}\cong
M_{(ll')^\infty}$ obtained by pairing the $j$th tensor factors.
Any other identification of these tails differs by an automorphism
of $M_{(ll')^\infty}$, whose automorphism group is contractible
by~\eqref{eq:H}. The resulting stabilization maps are homotopic.
The restriction bundles of
Lemma~\ref{lem:restriction} compare this group map with the induced
map of frame spaces. At stage $N$ the latter tensors by $1_{(l')^N}$,
and on every fixed stable odd homotopy group it multiplies by
$(l')^N$ modulo $k$. Pairing the tensor factors in their stage order
makes these maps compatible with stabilization. On the odd groups the
two composites from stage $N$ to stage $N+1$ have the same multiplier
$(ll')(l')^N=l(l')^{N+1}$. Since $(k,l')=1$, each $(l')^N$ is a unit
modulo $k$, so the induced map on the direct limit is an isomorphism.
The even groups vanish. The same calculation applies to the fundamental
group, and all spaces are connected. Thus the group map and its
classifying-space map are weak equivalences. Interchange $l,l'$ for
the other leg.
\end{proof}

\paragraph{Tensor products are graded.}
Tensor product gives a continuous functor
$\C_{a,l}\times\C_{b,l}\to\C_{ab,l}$, sending object indices
$(n,m)$ to $n+m$ and tensoring embeddings. It induces natural pairings
of the classified theories. To distinguish the type from the
finite-stage index in $D_n=M_{kl^n}$, write $D(a)=M_a\otimes A$.
Interleaving tensor factors gives $A\otimes A\cong A$, hence
\[
             D(a)\otimes D(b)\cong D(ab).
\]
For $k>1$ take a prime $p\mid k$. Since $p\nmid l$, its exponents
in the supernatural numbers of $D(k)$ and $D(k^2)$ are respectively
$v_p(k)$ and $2v_p(k)$. These fibers are therefore nonisomorphic.
The family of types $k^n$ is closed under tensor product.

\section{Tensor group completion}\label{sec:group-completion}

Throughout Sections~\ref{sec:group-completion}--\ref{sec:brauer-relative}
assume $k,l\geq2$ and $\gcd(k,l)=1$. The exponent $n$ in $k^n$ will
record tensor degree. The exponent $r$ in $l^r$ still indexes finite
matrix stages within one degree. Write
\[
 X_n=|N\C_{k^n,l}|,\qquad
 A=M_{l^\infty},\qquad D^{(n)}=M_{k^n}\otimes A,\qquad
 A_\infty=A\otimes M_{k^\infty}.
\]
In particular $D^{(0)}=A$. This superscript convention will only be
used in these sections.

\begin{definition}
Let $\mathcal T_{k,l}=\coprod_{n\geq0}\C_{k^n,l}$ as a topological
category. Its objects are pairs $(n,r)$ representing $M_{k^n l^r}$.
Morphisms preserve $n$. Tensor product on objects is addition in both
indices and on arrows is the tensor product of embeddings. Write
$P=|N\mathcal T_{k,l}|$ and let $\mathbf L_{k,l}$ be the connective
spectrum associated by Segal group completion to this tensor category.
Thus $\Omega^\infty\mathbf L_{k,l}=P^{\mathrm{gp}}$ up to equivalence.
\end{definition}

One can use standard matrix units in lexicographic order for the
tensor identifications. Rebracketing and permuting tensor factors give
the associativity and symmetry isomorphisms, with unit $M_1=\mathbb C$.
Their coherence follows from the corresponding permutations of tensor
bases. This is Segal's construction for a symmetric monoidal category
\cite{Segal74}. The graded bifunctors and the suggestion of an $H$-space
after stabilization occur in \cite[Section~3.1]{ErshovCategory}.
The link with classical constructions in algebraic $K$-theory lies in
Segal's group completion. Here its input is tensor product of unital
matrix algebras: the virtual integer counts tensor factors $M_k$.
Classical additive constructions instead start from direct sum of
modules and record virtual rank.
The disjoint union in the definition matters. Allowing every unital
embedding across all tensor degrees would make $\mathbb C$ initial in
the entire category and would contract its nerve.

\begin{proposition}[The degree-zero component]\label{prop:tensor-telescope}
The degree map identifies $\pi_0\mathbf L_{k,l}$ with $\mathbb Z$.
There are weak equivalences
\begin{equation}\label{eq:tensor-telescope}
 (\Omega^\infty\mathbf L_{k,l})_0
 \ \simeq\ \hocolim_n X_n
 \ \simeq\ \hocolim_n B\Aut(D^{(n)}).
\end{equation}
The transition functors tensor with $M_k$, and the homomorphisms of
automorphism groups tensor with $\id_{M_k}$.
\end{proposition}
\begin{proof}
Each $X_n$ is connected and simply connected by
Theorem~\ref{thm:main} and Proposition~\ref{prop:homotopy}, also for
$n=0$. Thus $\pi_0P=\mathbb N_0$. Apply the McDuff--Segal group-completion
theorem \cite{McDuffSegal} to a well-pointed monoid model of the tensor
$E_\infty$-space $P$ (one may use a derived bar construction).
The component of $t=M_k$ generates $\pi_0P$ and the Pontryagin product
is graded commutative, so the telescope of right multiplication by $t$
maps to $P^{\mathrm{gp}}$ by an integral homology equivalence.
In its component of virtual degree zero this is the map from
$\hocolim_n X_n$ that sends $x\in X_n$ to $x\,t^{-n}$ in the group
completion. The telescope is simply connected by van Kampen. The target
component is an $H$-space, so its fundamental group is abelian. Its
first homology vanishes by the homology comparison, hence it too is
simply connected. The homological Whitehead theorem proves the first
equivalence. The group of components is the group completion of
$\mathbb N_0$.

For the second equivalence use the zigzag in
Proposition~\ref{prop:two-inclusions} in each degree. Tensoring with
$M_k$ acts also on the enlarged categories: it carries a finite-to-limit
embedding $\varphi$ to $\varphi\otimes\id_{M_k}$ and a limit
automorphism $\alpha$ to $\alpha\otimes\id_{M_k}$. Identify the new
matrix factor with the next $k$-factor by the tensor flip past the
tail. These are functors between consecutive zigzags and commute with
the two inclusions. Taking telescopes gives~\eqref{eq:tensor-telescope}.
\end{proof}

\begin{corollary}\label{cor:telescope-classes}
Put $\mathbf L^0_{k,l}=\tau_{\geq1}\mathbf L_{k,l}$, the connected
cover of the spectrum. For a finite CW complex $X$ there is a natural
bijection
\[
 \mathsf F^{\mathrm{st}}_{k,l}(X)
 :=\operatorname*{colim}_{n}\mathsf F_{k^n,l}(X)
 \ \cong\ [X,\Omega^\infty\mathbf L^0_{k,l}].
\]
Here stabilization tensors each term of a LAB with the trivial
$M_k$-bundle. The same statement holds for the two coherent theories
of Section~\ref{sec:coherent}. Tensor product induces the abelian group
operation on $\mathsf F^{\mathrm{st}}_{k,l}(X)$.
\end{corollary}
\begin{proof}
The comparison $\bar\Phi$ uses the nerve of the finite matrix diagram
on a cover. Tensoring that diagram tensors its nerve map, after using
a common refinement for two different covers. Thus the class maps in
Theorem~\ref{thm:intro-classification} respect these operations.
A map from $X$, or a homotopy from $X\times I$, into a telescope has
bounded height and can be moved to its last finite stage. This proves
the bijection by Proposition~\ref{prop:tensor-telescope}. The formula
$(x\,t^{-m})(y\,t^{-n})=(x\otimes y)t^{-(m+n)}$ identifies the
operation.
\end{proof}

\section{Sections of prescribed rank and the relative Picard spectrum}\label{sec:relative-picard}

Put $R=KU[1/l]$ and $\Lambda=\pi_0 R=\mathbb Z[1/l]$. Write
\[
 U_R=GL_1(R)_1,\qquad Y_m=(\Omega^\infty R)_m
 \quad (m\geq1).
\]
The subscripts select the components of $1$ and $m$, respectively.
Multiplication in $R$ makes $U_R$ a grouplike space acting on $Y_m$.
Write $Y_m//U_R$ for the homotopy quotient, the total space of the
associated $Y_m$-bundle over $B U_R$. An \emph{orientable $R$-line
with a chosen reduction} is a family of $R$-modules equivalent to
$R$ whose structure group has been reduced to $U_R$. The reduction
identifies its local system of fiber ranks with the constant group
$\Lambda$. Here an orientable twist is taken with this choice of
reduction whenever a fiber rank is specified. All section spaces, homotopy orbits and fibers in this
section are derived.

\subsection{The scalar section at fixed matrix size}

\begin{theorem}[Fixed-size classification]\label{thm:rank-section}
For every positive $m$ coprime to $l$ there is a weak equivalence
\begin{equation}\label{eq:rank-model}
          |N\C_{m,l}|\ \xrightarrow{\ \simeq\ }\ Y_m//U_R .
\end{equation}
Thus on finite CW complexes the strict and coherent classes of type
$(m,l)$ classify orientable $R$-lines with a chosen reduction and
a section of fiber rank $m$, up to equivalence of pairs.
Tensoring with the constant matrix algebra $M_a$, where $(a,l)=1$,
carries the section $x$ to $ax$ in the rank-$am$ model over $B U_R$.
\end{theorem}
\begin{proof}
Consider pairs $(L,s)$ with $L$ equivalent to $R$ as an $R$-module
and $s:R\to L$ an $R$-module map. A trivialization identifies $s$
with a point of $\Omega^\infty R$. Changing that trivialization
acts by $GL_1(R)$, so the space of pairs and their equivalences is
$(\Omega^\infty R)//GL_1(R)$. The components of $GL_1(R)$ act on
fiber ranks by $\Lambda^\times$. The stabilizer of $m$ in this
discrete group is trivial, since $m\ne0$ and $\Lambda\subset\mathbb Q$.
The component of the pair $(R,m)$ is therefore $Y_m//U_R$. Indeed
the union of components of $\Omega^\infty R$ in the orbit of $m$
is induced from the $U_R$-space $Y_m$. An oriented trivialization
in a family is a local coordinate in which the section has rank $m$.

Use spectrum-valued operator $K$-theory with its lax symmetric
monoidal refinement to $KU$-modules \cite{LandNikolaus},
\cite[Sections~2.2, 5.1 and~5.7]{BLP}. It acts on families of arrows:
the homomorphism spaces in the homotopy localization of complex
$C^*$-algebras carry the point-norm model. For $d=ml^r$ assign to
$M_d$ the pair
\begin{equation}\label{eq:scalar-pair}
     L_d=K(M_d)\otimes_{KU}R,\qquad s_d:R\longrightarrow L_d
\end{equation}
induced by its scalar inclusion. A rank-one corner gives a Morita
coordinate $L_d\simeq R$ in which $s_d$ has rank $d$. An embedding of
multiplicity $l^b$ acts on the lines by an equivalence, since it
multiplies this coordinate by the unit $l^b$. Unitality makes it
an equivalence of pairs. Changing the coordinate by $l^{-r}$ puts
the section in rank $m$. This constructs the map in~\eqref{eq:rank-model}
intrinsically from the pairs.

The orbit map $U_R\to Y_m$, $v\mapsto vm$, is the boundary after
looping the base of $Y_m\to Y_m//U_R\to B U_R$. On positive
homotopy groups it acts by multiplication by $m$ on $\pi_*R$:
translate the additive basepoints $1$ and $m$ to zero to see this.
Bott periodicity and injectivity of $m$ on $\Lambda$ give, for $j\ge1$,
\[
 \pi_{2j}(Y_m//U_R)=\Lambda/m\Lambda,\qquad
 \pi_{2j-1}(Y_m//U_R)=0.
\]
The even group is the cokernel of the orbit map on $\pi_{2j}$.
To determine the comparison on this cokernel, choose $r$ with
$d=ml^r\geq j$. There is a rank-$d$ bundle $V$ on $S^{2j}$ with
$[V]=d+u$, where $u$ is a reduced Bott generator: the inclusion
$U(d)\to U$ is stable on $\pi_{2j-1}$ since $2j-1<2d$.
The endomorphism bundle defines a map into the nerve through
$B\mathrm{PU}(d)$. Embed $V$ in a trivial Hilbert bundle and let
$p$ be the resulting projection. Fiberwise corner inclusion trivializes
the $K$-theory line of $\operatorname{End}(V)$ by Morita invariance.
On section algebras this inclusion reads
\[
  \Gamma(S^{2j},\operatorname{End}(V))
     =p\,C(S^{2j},\mathbb K)\,p\ \lhook\joinrel\longrightarrow\
       C(S^{2j},\mathbb K).
\]
It sends the unit of the corner to $p$, whose $K^0$-class is $[V]$.
The natural comparison between the $K$-theory of the section algebra
and the sections of its $K$-theory local system applies to the finite
complex $S^{2j}$ \cite[Proposition~5.30]{BLP}. Thus the section assigned
by~\eqref{eq:scalar-pair}, expressed in the rank-$m$ coordinate, is
\[
                        l^{-r}[V]=m+l^{-r}u .
\]
Its reduced coefficient $l^{-r}$ generates $\Lambda/m\Lambda$.
The source is connected and has groups $\mathbb Z/m$ in positive
even degrees and zero in odd degrees by Proposition~\ref{prop:homotopy}.
The map is consequently an isomorphism on every homotopy group,
including for $m=1$, when both sides are weakly contractible.

Finally the exterior product $L_d\otimes_R L_e\to L_{de}$ is an
equivalence: it takes the product of minimal projection classes to
a minimal projection class. It carries $s_d\otimes s_e$ to $s_{de}$.
In the Morita coordinate for $M_a$ its unit has rank $a$, giving the
map $x\mapsto ax$ over $B U_R$. Tensor powers of a single chosen
rank-one corner give compatible coordinates at successive stages.
The classification of LAB classes now follows from
Theorem~\ref{thm:intro-classification}.
\end{proof}

\begin{corollary}[The prescribed-rank criterion]\label{cor:prescribed-rank}
Fix $m\ge1$ coprime to $l$. Let $X$ be a connected finite CW complex,
$x\in X$, and let
$\tau:X\to B U_R$ specify an orientable twist with its rank
coordinate. It admits a LAB lift of type $(m,l)$ if and only if
\begin{equation}\label{eq:rank-image}
 m\ \in\ \operatorname{im}\!\left(
 K^0(X;\tau)\xrightarrow{\mathrm{ev}_x}\Lambda\right).
\end{equation}
Here $K^0(X;\tau)=\pi_0\Gamma(X,L_\tau)$ uses the spectrum of
sections of the $R$-line classified by $\tau$.
\end{corollary}
\begin{proof}
A lift along $Y_m//U_R\to B U_R$ is a section of the associated
bundle of rank-$m$ points of $\Omega^\infty L_\tau$. The functor
$\Omega^\infty$ preserves limits, so twisted $K^0$-classes are
homotopy classes of sections of $\Omega^\infty L_\tau$. Fiber rank
is locally constant and the chosen reduction fixes its coordinate.
Thus rank $m$ at $x$ gives rank $m$ on all of $X$.
Apply Theorem~\ref{thm:rank-section}.
\end{proof}

\subsection{Families of LABs and their symmetries}\label{subsec:symmetries}

For a connected finite CW complex $X$ write $G=\Aut(D)$ and let
$\mathcal L_{k,l}(X)$ denote the completed coherent descent space of
Theorem~\ref{thm:coherent}. Its description as a stack means that the
space on $X$ is recovered from the spaces on a cover and all their
overlaps by homotopy-coherent descent. Theorems~\ref{thm:coherent}
and~\ref{thm:main} give
\[
  \mathcal L_{k,l}(X)\simeq \mathbf R\Map(X,|N\C_{k,l}|)
                     \simeq \mathbf R\Map(X,BG).
\]
Its components give the classification set $\mathsf F(X)$. If a point
corresponds to a $D$-bundle $E$, its based loop space is
$\Omega_E\mathcal L_{k,l}(X)\simeq\operatorname{Gau}(E)$, the topological
group of fiberwise algebra automorphisms covering $\id_X$. To see this,
let $f:X\to BG$ classify the frame bundle $P$ of $E$. A loop at $f$
in the mapping space is a section of the pullback along $f$ of the
free-loop fibration $LBG\to BG$. A model for this fibration is
$EG\times_G G\to BG$ with the conjugation action: pulling the path
fibration back along the diagonal restricts the left and right
actions on its fiber $G$ to conjugation. Hence
\[
  \Omega_f\mathbf R\Map(X,BG)\simeq\Gamma(X,P\times_G G)
       =\operatorname{Gau}(E).
\]
This is the gauge-group description of the component of $f$
\cite[Proposition~2.4]{AtiyahBottGauge}.

\paragraph{The circle.}
For $k>1$ there is a single LAB class on $S^1$, since $BG$ is simply
connected. The trivial bundle has gauge group $\Map(S^1,G)$.
Evaluation at $1$ and the assignment $u\mapsto u(1)^{-1}u(-)$
identify this mapping space with $G\times\Omega G$ as a space.
Consequently
\[
  \pi_1\mathcal L_{k,l}(S^1)\cong\pi_0\Map(S^1,G)
              \cong\pi_1 G\cong\mathbb Z/k .
\]
A generator acts on the fiber at $z\in S^1$ by
$\Ad(\diag(z,1,\ldots,1))\otimes\id_{M_{l^\infty}}$.
Indeed its projective unitary loop clutches
$\operatorname{End}(H\oplus\underline{\mathbb C}^{\,k-1})$ on $S^2$,
up to the choice of clutching orientation. Example~\ref{ex:morita-versus-lab}
shows that stabilization sends this loop to a generator of $\pi_1 G$.

\paragraph{The space of lifts of a fixed twist.}
Let $\tau:X\to B U_R$ be an orientable twist with its rank coordinate,
and let $P_\tau$ be the associated principal $U_R$-bundle. Under
Theorem~\ref{thm:rank-section} the map forgetting the section has
homotopy fiber
\[
 \hofib_\tau\!\left(\mathcal L_{k,l}(X)
                \longrightarrow\mathbf R\Map(X,B U_R)\right)
       \simeq \Gamma\!\left(X,P_\tau\times_{U_R}Y_k\right).
\]
A point of the fiber specifies a lift together with an identification of
its underlying twist with $\tau$. The formula follows by pulling the
associated $Y_k$-bundle over $B U_R$ back along $\tau$ and taking its
sections. Its components are exactly the classes of fiber rank $k$ in
$K^0(X;\tau)$. At a chosen section $s$, its $j$th homotopy group is
$K^{-j}(X;\tau)$ for $j\ge1$. Indeed the rank-$k$ sections form a union of
components of $\Omega^\infty\Gamma(X,L_\tau)$, and translation by $-s$
identifies the component of $s$ with that of zero. If only the class
$[\tau]$ is prescribed, the LAB classes above it are the orbits of the
rank-$k$ classes under
$\pi_1(\mathbf R\Map(X,B U_R),\tau)$, by the homotopy exact sequence.
These are the changes of identification of the underlying twist.

For example, trivialize the twist on $S^2$ and let $u$ generate
$\widetilde K^0(S^2)$, so that $u^2=0$. The identified rank-$k$ lifts
are $k+a u$ for $a\in\Lambda$. A change of identification is represented
by a unit $1+b u$ for $b\in\Lambda$. Multiplication gives
\[
                    (1+b u)(k+a u)=k+(a+kb)u .
\]
Thus the LAB classes over the class of the trivial twist are
$\Lambda/k\Lambda\cong\mathbb Z/k$. The Hopf-line example in the
introduction has $a=1$ for the generator $u=[H]-1$.

\subsection{Inverting the matrix size}

Set $R'=R[1/k]$, $\Lambda'=\mathbb Z[1/(kl)]$ and $U_{R'}=GL_1(R')_1$.
For a commutative ring spectrum $S$, its connective Picard spectrum
$\mathrm{pic}(S)$ has as zeroth space the invertible $S$-modules
and their equivalences. The component of $S$ is $BGL_1(S)$
\cite{ABGHR}. In particular $\pi_1\mathrm{pic}(S)=(\pi_0 S)^\times$
and $\pi_i\mathrm{pic}(S)=\pi_{i-1}S$ for $i\geq2$. Put
\[
 E_l=\tau_{\geq2}\mathrm{pic}(R),\qquad
 E_{kl}=\tau_{\geq2}\mathrm{pic}(R'),\qquad
 \mathbf Q_{k,l}=\operatorname{fib}(E_l\longrightarrow E_{kl}).
\]
We have $\Omega^\infty E_l\simeq B U_R$ and
$\Omega^\infty E_{kl}\simeq B U_{R'}$.

\begin{proposition}[The comparison map]\label{prop:pic-map}
The scalar pairs~\eqref{eq:scalar-pair} induce a map of spectra
\[
 \mathbf L_{k,l}\longrightarrow
 \operatorname{fib}\bigl(\mathrm{pic}(R)\to\mathrm{pic}(R')\bigr)
\]
and a map on connected covers
\begin{equation}\label{eq:relative-spectrum-map}
              \Theta:\mathbf L^0_{k,l}\longrightarrow\mathbf Q_{k,l}.
\end{equation}
\end{proposition}
\begin{proof}
If $d=k^n l^r$, the section $s'_d=R'\otimes_R s_d$ is an
equivalence: in a Morita coordinate it multiplies by the unit $d$
of $\Lambda'$. Hence $(L_d,s'_d)$ is an $R$-line equipped with a
trivialization over $R'$, a point of the zeroth space of the full
fiber. The scalar inclusions are monoidally natural in unital
embeddings. The exterior products in the proof of
Theorem~\ref{thm:rank-section} and the unit at $M_1$ are equivalences.
The lax symmetric monoidal $K$-functor supplies their associativity,
unit and permutation coherences simultaneously. It therefore gives
an $E_\infty$-map from the graded matrix category on realization to
this grouplike space of pairs. Group completion gives the first
map, through the connective cover of the full fiber.

The map from $\mathbf Q_{k,l}$ to the full fiber is an isomorphism on
$\pi_i$ for $i\ge2$ by the fiber exact sequences. Both $\pi_1$ groups
vanish, since $\pi_2\mathrm{pic}(R')=0$ and
$\Lambda^\times\to(\Lambda')^\times$ is injective. The fiber of
the displayed $2$-connective covers is itself at least $1$-connective.
It is therefore the connected cover of the full fiber.
\end{proof}

\begin{theorem}[Spectrum comparison]\label{thm:relative-spectrum}
The map $\Theta$ is an equivalence of spectra. For $j\ge1$,
\[
 \pi_{2j}\mathbf L^0_{k,l}\cong\Lambda'/\Lambda
       \cong\mathbb Z[1/k]/\mathbb Z,\qquad
 \pi_{2j-1}\mathbf L^0_{k,l}=0 .
\]
The degree-$n$ space $X_n$, translated by $t^{-n}$, maps on
$\pi_{2j}$ isomorphically onto the subgroup killed by $k^n$.
\end{theorem}
\begin{proof}
The fixed-size models for $m=k^n$ have the same base $B U_R$.
Their fibers form the telescope
\[
        Y_1\xrightarrow{\;k\;}Y_k\xrightarrow{\;k\;}Y_{k^2}
                                   \xrightarrow{\;k\;}\cdots .
\]
On the $n$th fiber, coefficient extension followed by division by
$k^n$ gives a $U_R$-equivariant map $x\mapsto x/k^n$ to $U_{R'}$,
where $U_R$ acts through $U_R\to U_{R'}$. Choose the successive
powers of one inverse to the constant unit $k$ in $R'$; this gives
the compatible maps to the target of the telescope. In positive
degrees the source homotopy groups are obtained from $\pi_*R$ by
inverting $k$, and the map sends a stage-$n$ class $a$ to $a/k^n$.
It is thus an isomorphism onto $\pi_*R'$ in these degrees. Both
spaces are connected. Homotopy orbits commute with telescopes, so
\begin{equation}\label{eq:unit-orbit-fiber}
 \hocolim_n X_n\ \simeq\ U_{R'}//U_R\
   \simeq\ \hofib(B U_R\longrightarrow B U_{R'}).
\end{equation}
The last identification uses the principal $U_{R'}$-bundle over
$B U_R$ associated by extension of structure group. We use the
fiber convention in which a trivialization is an equivalence from
the constant $R'$-line to the extended line. A change of coordinates
then acts by left multiplication on $U_{R'}$.

Under Proposition~\ref{prop:tensor-telescope} this equivalence is the
zeroth-space map of $\Theta$. Indeed an element of $X_n$ enters the
group completion as its pair tensored with the inverse of the constant
degree-$n$ pair. In the Morita coordinate the reference section is
$k^n$. Its inverse gives exactly the map $x\mapsto x/k^n$ in
\eqref{eq:unit-orbit-fiber}. The target is $\Omega^\infty\mathbf Q_{k,l}$
since $\Omega^\infty$ preserves fibers. The connected spectra in
\eqref{eq:relative-spectrum-map} are connective, so their zeroth-space
equivalence detects a spectrum equivalence.

The exact sequence for $\mathbf Q_{k,l}$ gives the even group
$\Lambda'/\Lambda$ and the vanishing odd groups by Bott periodicity.
At stage $m=k^n$ the homomorphism on even groups is
\[
    \Lambda/m\Lambda\longrightarrow\Lambda'/\Lambda,\qquad
             [a]\longmapsto[a/m].
\]
For the sphere bundle of Theorem~\ref{thm:rank-section}, the
matrix size is $d=ml^r=k^nl^r$ and $[V]=d+u$. Its rank-$m$
section has reduced term $l^{-r}u$. Dividing by the reference rank
$m$ gives the based unit and its reduced coefficient
\[
       \frac{l^{-r}[V]}{m}=1+\frac{u}{d},\qquad
       \left[\frac1d\right]\in\Lambda'/\Lambda .
\]
Thus the coefficient of this Bott test is $1/d$ modulo $\Lambda$.
The class of $1/k^n$ is the image of $[1]\in\Lambda/m\Lambda$.
The two classes differ by the unit factor $l^{-r}$ and generate
the same cyclic subgroup $m^{-1}\Lambda/\Lambda$, of order $m$.
This subgroup is the image of the stage, and the transition between
successive cyclic groups is multiplication by $k$.
Finally $\mathbb Z[1/k]\cap\Lambda=\mathbb Z$, and a power
of $l$ is a unit modulo every power of $k$. It follows that the
inclusion $\mathbb Z[1/k]\to\Lambda'$ induces the asserted isomorphism
on quotients.
\end{proof}

\begin{corollary}[The primes and the degree]\label{cor:prime-degree}
Let $E=\tau_{\geq2}\mathrm{pic}(KU)$.
There are natural equivalences
\[
 \mathbf Q_{k,l}\ \simeq\ \operatorname{fib}(E\to E[1/k])
  \ \simeq\ \Sigma^{-1}\bigl(E\wedge
        \operatorname{cofib}(\mathbb S\longrightarrow\mathbb S[1/k])\bigr).
\]
In particular the connected spectrum depends only on the set of primes
dividing $k$. For the full group completion put
$P_{\mathrm{rel}}=\tau_{\geq0}\operatorname{fib}
(\mathrm{pic}(R)\to\mathrm{pic}(R'))$. If
$\delta:\mathbb Z\to\pi_0P_{\mathrm{rel}}$ sends $1$ to the component
of the pair $(R,k)$ after localization, then there is a homotopy
pullback description
\[
   \mathbf L_{k,l}\simeq
        P_{\mathrm{rel}}\times_{H\pi_0P_{\mathrm{rel}}}H\mathbb Z .
\]
Here $\mathbb S$ is the sphere spectrum, $H$ denotes the
Eilenberg--Mac Lane functor and the map from $H\mathbb Z$ is $H\delta$.
\end{corollary}
\begin{proof}
The coefficient map $E\to E_l$ identifies $E_l$ with $E[1/l]$:
its target is $l$-local and it induces localization on homotopy
groups in every degree. Similarly $E_{kl}\simeq E_l[1/k]$.
The natural map $\operatorname{fib}(E\to E[1/k])\to\mathbf Q_{k,l}$
is the quotient isomorphism of Theorem~\ref{thm:relative-spectrum}
in positive even degrees and an isomorphism of zero groups elsewhere.
The cofiber formula follows by smashing with $E$ in spectra.

The exact sequence for the full Picard fiber embeds
$(\Lambda')^\times/\Lambda^\times$ in $\pi_0P_{\mathrm{rel}}$.
The component of the pair with localized scalar $k$ has infinite
order there. The map of Proposition~\ref{prop:pic-map} therefore
induces $\delta$ on $\pi_0$ when the source is indexed by virtual
degree. Taking degree and this map together gives the claimed map
to the pullback. The pullback is connective, has component group
$\mathbb Z$, and agrees with $P_{\mathrm{rel}}$ on positive
homotopy groups by the exact sequence for the zeroth Postnikov
section. Theorem~\ref{thm:relative-spectrum} and stable Whitehead
give the equivalence.
\end{proof}

For a based finite CW complex $X$, the reduced degree-zero theory
therefore extends to all cohomological degrees by
\[
   \widetilde{\mathsf F}^{\,i}_{k,l}(X)=
       [\Sigma^\infty X,\Sigma^i\mathbf Q_{k,l}]_{\mathrm{Sp}},
       \qquad i\in\mathbb Z.
\]

\section{Orientable twists and higher Brauer classes}\label{sec:brauer-relative}

For an infinite-type UHF algebra $J$, let
\[
       \mathrm{Tw}(J)=B\Aut_0(J\otimes\mathbb K),
\]
where $\mathbb K$ is the algebra of compact operators and the
subscript denotes the identity path component. A map into this space
classifies a stabilized algebra bundle with its structure group reduced
to that component. We use tensor product for the infinite loop
structure on orientable twists. The scalar extension
$i:A\to A_\infty$ acts on twists by tensoring with the trivial
$M_{k^\infty}$-bundle.

\begin{unnumberedlemma}[The operator $K$-theory action and tensor extension]
Write $B=M_{k^\infty}$, $A_\infty=A\otimes B$, $R=KU[1/l]$ and
$R'=KU[1/(kl)]$. The action on the operator $K$-theory line gives the
vertical equivalences in a homotopy-commutative square of infinite loop spaces
\[
\begin{array}{ccc}
 \mathrm{Tw}(A) & \xrightarrow{\ B t_1\ } & \mathrm{Tw}(A_\infty)\\
 {\scriptstyle\simeq}\big\downarrow && {\scriptstyle\simeq}\big\downarrow\\
 \Omega^\infty\tau_{\geq2}\mathrm{pic}(R)
 & \longrightarrow & \Omega^\infty\tau_{\geq2}\mathrm{pic}(R').
\end{array}
\]
Here $t_1(\alpha)=\alpha\otimes\id_B$ after shuffling the compact factor.
The bottom map is extension of scalars along $R\to R'$.
\end{unnumberedlemma}
\begin{proof}
Use the same lax symmetric monoidal functor $K$ as in
Theorem~\ref{thm:rank-section}. For $J=A,A_\infty$ put $R_J=R,R'$,
respectively, $C_J=J\otimes\mathbb K$, and $p_J=1_J\otimes e$ for
a rank-one compact projection $e$. In tensor degree $n$ put
\[
 L_{J,n}=K(C_J^{\otimes n})\otimes_{KU}R_J,\qquad
 p_{J,n}=p_J^{\otimes n},\qquad C_J^{\otimes0}=\mathbb C,\quad p_{J,0}=1 .
\]
The corner map $\lambda\mapsto\lambda p_{J,n}$ induces an equivalence
$R_J\to L_{J,n}$. For $n>0$ its projection represents $1$ in the
ordered $K_0$-group of the stabilized UHF algebra, in the coordinate
given by the unit of $J^{\otimes n}$. This group is $\pi_0 R_J$.
Bott periodicity gives the same isomorphism in all even degrees, and
the odd groups vanish. In degree zero of the tensor construction the
map is the identity of $R_J$. Exterior product carries the corners in
degrees $a,b$ to the corner in degree $a+b$, so the tensor comparison
$L_{J,a}\otimes_{R_J}L_{J,b}\to L_{J,a+b}$ is an equivalence.

We describe the tensor coherence of the action. In the Dadarlat--Pennig
model indexed by finite sets and injections, degree $n$ has group
$G_J(n)=\Aut(C_J^{\otimes n})$. For an injection $\sigma$ from $n$
indices to $m$, write $\sigma_*g$ for the automorphism obtained by
adjoining identity automorphisms and permuting the factors. Let
$c_\sigma:C_J^{\otimes n}\to C_J^{\otimes m}$ insert $p_J$ in each
unused position. Concretely, the category assembling the one-object
groupoids $*G_J(n)$ has arrows $(\sigma,g):n\to m$, $g\in G_J(m)$,
with composition $(\tau,h)(\sigma,g)=(\tau\sigma,h\,\tau_*g)$.
The assignment
\[
                (\sigma,g)\longmapsto g\,c_\sigma
\]
is functorial because $c_\tau g=(\tau_*g)c_\tau$. It preserves exterior
products and permutations. Its composite with localized $K$-theory
takes values in the groupoid of $R_J$-lines: each $c_\sigma$ takes the
reference corner to the reference corner, and every $g$ is invertible.
The exterior-product equivalences above and the symmetric monoidal
coherences of $K$ give a symmetric monoidal functor into this groupoid.
Here the groupoid of lines includes the full spaces of equivalences
of the modules.

The category of pairs is the Grothendieck construction of
$n\mapsto *G_J(n)$. Its derived classifying space models the
homotopy colimit of $BG_J(n)$ over finite sets and injections:
apply Thomason's homotopy colimit theorem in each singular degree
and take diagonals \cite{ThomasonHocolim}. To identify its tensor
structure, use the permutative category $\mathcal B_\otimes(J)$ of
\cite[Theorem~4.9]{DPUnits}. Its arrows between tensor powers carry
the reference projection to a projection of invertible $K_0$-class.
Consequently localized $K$-theory sends every such arrow to an
equivalence of lines. The functor $(\sigma,g)\mapsto g\,c_\sigma$
takes values in $\mathcal B_\otimes(J)$ and is the comparison built
from corner stabilization in \cite[proof of Theorem~3.13]{DPUnits}:
an injection inserts the corner, and a one-factor automorphism acts on
that factor. Theorem~4.9 there identifies the indexed and permutative
models as infinite loop spaces. Thus on $B\Aut(C_J)$ with its tensor
infinite loop structure the resulting map classifies the action
$g\mapsto K(g)_{R_J}$ on $L_{J,1}$.

To check the equivalence on simply connected covers, loop this map
and evaluate an automorphism of the free $R_J$-line at its generator.
For a based family $g:S^a\to\Aut(C_J)$, $a\geq1$, the induced
map in degree $a$ is
\[
   [g]\longmapsto
      [\,s\mapsto g(s)(p_J)\,]-[\,s\mapsto p_J\,]
       \quad\text{in }\widetilde K^0(S^a)\otimes\pi_0R_J .
\]
Indeed the composite with the generating corner is the corner of this
continuous projection family. This is the projection calculation in
the proof of \cite[Theorem~4.6]{DPUnits}, where the graded model
multiplies the projection by the universal Bott class. It is an
isomorphism for every $a>0$. Thus the coherent map already constructed
restricts to an equivalence from $B\Aut_0(C_J)$ to
$B(GL_1(R_J)_1)\simeq\Omega^\infty\tau_{\geq2}\mathrm{pic}(R_J)$.
More generally, on families in degree $n$ the projection is
$g(p_{J,n})$. It satisfies
$(g\otimes h)(p_{J,a+b})=g(p_{J,a})\otimes h(p_{J,b})$.

Finally shuffle $C_{A_\infty}^{\otimes n}$ into the order
$C_A^{\otimes n}\otimes B^{\otimes n}$ and set
$i_n(x)=x\otimes1_B^{\otimes n}$ and
$t_n(g)=g\otimes\id_{B^{\otimes n}}$. These maps satisfy
\[
 t_n(g)i_n=i_ng,\qquad
 i_m c_\sigma^A=c_\sigma^{A_\infty}i_n,\qquad
 i_{a+b}=i_a\otimes i_b
\]
under the same shuffles. In particular $t_n$ commutes with adjoining
identities and with permutations. The maps of lines
$R'\otimes_R L_{A,n}\to L_{A_\infty,n}$ induced by $i_n$ are
equivalences since they preserve the reference corners. The displayed
identities make them a symmetric monoidal natural equivalence on the
category of pairs $(\sigma,g)$, including its continuous arrow
families. Taking classifying spaces and simply connected covers gives
the stated square.
\end{proof}

\begin{corollary}[The relative twist and its trivialization]\label{cor:relative-twist}
There is an equivalence of infinite loop spaces
\[
 (\Omega^\infty\mathbf L_{k,l})_0 \ \simeq\
 \hofib_*\bigl(\mathrm{Tw}(A)\longrightarrow\mathrm{Tw}(A_\infty)\bigr).
\]
Consequently, on a finite CW complex, a tensor-stabilized LAB
classifies a pair consisting of an orientable $A$-twist and a
trivialization up to homotopy of its extension to $A_\infty$.
Equivalence of pairs transports the trivialization along the
equivalence of twists.
\end{corollary}
\begin{proof}
The lemma identifies the homotopy fiber of the tensor-extension map with
$\Omega^\infty\mathbf Q_{k,l}$, since $\Omega^\infty$ preserves homotopy
fibers. Apply Theorem~\ref{thm:relative-spectrum} and
Corollary~\ref{cor:telescope-classes}. Mapping $X$ into this fiber
specifies a twist on $X$ together with a nullhomotopy of its extension.
This nullhomotopy is the chosen trivialization.
\end{proof}

\paragraph{What the forgetful map does to a bundle.}
For a representative of tensor degree $n$ let $E$ be its rectified
$D^{(n)}$-bundle. The map to $\mathrm{Tw}(A)$ sends it to the
stabilized bundle $E\otimes\mathbb K$, using a fixed identification
$D^{(n)}\otimes\mathbb K\cong A\otimes\mathbb K$ coming from
matrix stability. To check this description against the map of
spectra, extend the functor of Proposition~\ref{prop:pic-map} over
the enlarged category of Section~\ref{sec:K} in degree $n$.
Assign to its extra object $D^{(n)}$ the $R$-line
$K(D^{(n)})\otimes_{KU}R$ and the localized section given by its
algebra unit. In the Morita coordinate in which a rank-one matrix
corner tensored with $1_A$ has class $1$, the unit has class $k^n$.
A finite-to-limit embedding from $M_{k^n l^r}$ sends a minimal
projection to class $l^{-r}$: its $k^n l^r$ equivalent diagonal
projections sum to the unit. This makes its map of $R$-lines an
equivalence, and unitality preserves the section. Along the other
inclusion in the nerve zigzag this functor is precisely the action
of $\Aut(D^{(n)})$ on its $K$-theory line. The full-corner inclusion
$x\mapsto x\otimes e$, with $e$ a rank-one compact projection,
is equivariant for $\alpha\mapsto\alpha\otimes\id$ and is a
$K$-theory equivalence. Hence the line obtained by forgetting the
section is associated to the stable frame bundle via the one-factor
action $g\mapsto K(g)_R$ of the preceding lemma. Under matrix
stability, a minimal projection in $M_{k^n}$ tensored with $1_A$
and $e$ gives its generating corner. Moreover $\Aut(D^{(n)})$ is
path connected by Proposition~\ref{prop:homotopy}, so the stabilized
structure group lies in the identity component. This identifies the
resulting orientable twist with that of $E\otimes\mathbb K$.
The same proof with $M_m\otimes A$ gives the geometric meaning of the
twist in Corollary~\ref{cor:prescribed-rank} for every prescribed size $m$.

To identify the localized trivialization, tensor $E$ with
$M_{k^\infty}$. Its unital fiber is isomorphic to $A_\infty$,
whose unital automorphism group is contractible
\cite[Theorem~2.3]{DP}. First compare the scalar sections before
translation to the identity component. Set
$J_n=D^{(n)}\otimes M_{k^\infty}$ and let
$i_n:D^{(n)}\to J_n$ be the unital inclusion. There is a commuting
square of $R'$-lines
\[
 \begin{array}{ccc}
 R' & \xrightarrow{\ s'_n\ } & K(D^{(n)})\otimes_{KU}R'\\
 \big\Vert && \big\downarrow{\scriptstyle K(i_n)_{R'}}\\
 R' & \xrightarrow{\ s'_{J_n}\ } & K(J_n)\otimes_{KU}R'.
 \end{array}
\]
Both rows are the scalar inclusions after localization. The lower row
is an equivalence since $J_n$ is unitally isomorphic to $A_\infty$.
The upper row is an equivalence because $k^n$ is invertible in
$\Lambda'$. The square is equivariant for $\alpha$ on the upper
right and $\alpha\otimes\id$ on the lower right. Thus the
unnormalized scalar sections agree over the entire classifying space.

In the stable Morita coordinate, $1_{D^{(n)}}\otimes e$ has class
$k^n$. A unital identification of $J_n$ with $A_\infty$ reads the
corresponding unit as $1$. The resulting change of scale is
$k^{-n}$, a nonidentity component of the units if $n>0$. For the
orientable trivialization compare with the constant degree-$n$
bundle in the same chart. Tensoring with its inverse pair divides
the localized section by its reference value $k^n$. The resulting
section $x/k^n$ has rank one, exactly as in~\eqref{eq:unit-orbit-fiber}.

\begin{theorem}[Forgetting the localized trivialization]\label{cor:brauer-forget}
Let $X$ be a connected finite CW complex. Write
$C_k=\mathbb Z[1/k]/\mathbb Z$ and
$E_l^i(X)=[\Sigma^\infty X_+,\Sigma^iE_l]_{\mathrm{Sp}}$.
There is a natural short exact sequence
\begin{equation}\label{eq:brauer-short-exact}
 0\longrightarrow E_l^{-1}(X)\otimes_{\mathbb Z} C_k
 \longrightarrow\mathsf F^{\mathrm{st}}_{k,l}(X)
 \longrightarrow\mathrm{Br}_A(X)[k^\infty]\longrightarrow0 .
\end{equation}
Here $[k^\infty]$ denotes the subgroup killed elementwise by powers
of $k$. For each $X$ the sequence splits as a sequence of abelian
groups. The kernel is abstractly isomorphic to $C_k^{\,b}$, where
\[
          b=\sum_{j\ge1}\dim_{\mathbb Q}H^{2j}(X;\mathbb Q).
\]
\end{theorem}
\begin{proof}
Applying maps from $\Sigma^\infty X_+$ to the fiber sequence of
Theorem~\ref{thm:relative-spectrum} gives the exact sequence
\[
 E_l^{-1}(X)\longrightarrow E_{kl}^{-1}(X)\longrightarrow
 \mathsf F^{\mathrm{st}}_{k,l}(X)\longrightarrow
 E_l^0(X)\longrightarrow E_{kl}^0(X).
\]
In degree zero the last map extends an orientable twist to
$A_\infty$ by Corollary~\ref{cor:relative-twist}. Since $E_{kl}\simeq
E_l[1/k]$ and $\Sigma^\infty X_+$ is finite, $E_{kl}^i(X)$
identifies with $E_l^i(X)[1/k]$. The kernel in degree zero is
therefore the $k$-power torsion. Under the Dadarlat--Pennig
identification of torsion in orientable twists with $\mathrm{Br}_A(X)$
this is precisely the right term in~\eqref{eq:brauer-short-exact}
\cite[Theorems~2.11 and~2.15]{DPBrauer}. In degree $-1$, right
exactness of tensor product with the sequence
$\mathbb Z\to\mathbb Z[1/k]\to C_k\to0$ identifies the cokernel
with the left term. This term acts simply transitively on the
choices above each Brauer class.

For any abelian group $M$, $M\otimes C_k$ is divisible: every
integer multiplication on $C_k$ is surjective. Divisible abelian
groups are injective, so the kernel inclusion admits a retraction.
To compute its isomorphism type, $E_l$ is $l$-local and its
homotopy groups are finitely generated over the PID $\Lambda$
in each degree. The finite cellular filtration gives the same
finite generation for $E_l^{-1}(X)$. Rational spectra split by
their homotopy groups. In this case there is one rational
coefficient in degrees $3,5,\ldots$, giving free rank $b$ in
cohomological degree $-1$ on $X$. The finitely generated torsion
part over $\Lambda$ is finite as an abelian group and tensors to
zero with $C_k$. Finally $\Lambda\otimes C_k\cong C_k$ since
multiplication by $l$ is an automorphism of $C_k$.
\end{proof}

\begin{example}[A nontrivial trivialization on the two-sphere]\label{ex:relative-half}
Take $k=2$, $l=3$, and choose the Hopf line $H$ with $c_1(H)=1$.
For $V=H\oplus\underline{\mathbb C}$ the bundle $\operatorname{End}(V)$
has trivial stable twist by Morita equivalence. Nevertheless its
class after tensor stabilization is
\[
       \tfrac12 \pmod{\mathbb Z[1/3]}
       \ \in\ \mathbb Z[1/6]/\mathbb Z[1/3] .
\]
Indeed $u=[H]-1$ is a reduced Bott generator and the normalized
section is $[V]/2=1+u/2$. Changing a trivialization before extending
coefficients changes the reduced coefficient by an element of
$\mathbb Z[1/3]$. The half therefore survives. For a direct
finite-dimensional check, suppose $\operatorname{End}(V)$ embeds
unitally into $\operatorname{End}(E_0)$, where
$E_0=\underline{\mathbb C}^{\,2\cdot3^r}$ is the fixed trivial
vector bundle. Give $E_0$ the resulting left module structure and put
\[
 W=\operatorname{Hom}_{\operatorname{End}(V)}(V,E_0).
\]
Locally $E_0$ is the sum of $3^r$ copies of the defining module, so
$W$ is a rank-$3^r$ bundle and evaluation identifies
$V\otimes W\cong E_0$. If one instead starts with a multiplicity
bundle $W_0$ known only through a trivialization of
$\operatorname{End}(V\otimes W_0)$, that trivialization gives
$V\otimes W_0\cong L\otimes E_0$ for a line bundle $L$:
local unitary lifts of the projective trivialization differ by
scalars. Replacing $W_0$ by $W_0\otimes L^*$ absorbs that ambiguity.
In either description, taking first Chern numbers now gives
$3^r+2\langle c_1(W),[S^2]\rangle=0$, which has no integral
solution. This explains the obstruction in the
introduction even though the Brauer class is zero.
\end{example}

\begin{example}[A higher Brauer class on a Moore space]\label{ex:relative-moore}
Let $q\geq2$ divide some power of $k$ and let
$X_q=S^4\cup_q e^5$ be the cofiber of the degree-$q$ self-map of
$S^4$. Applying the cofiber exact sequence and using
$\pi_4\mathrm{Tw}(A)=0$ and $\pi_5\mathrm{Tw}(A)=\Lambda$ gives
\[
 [X_q,\mathrm{Tw}(A)]\cong\Lambda/q\Lambda\cong\mathbb Z/q .
\]
For $A_\infty$ the same computation gives zero because $q$ is
invertible in $\Lambda'$. Furthermore $[X_q,\Omega\mathrm{Tw}(A_\infty)]=0$:
the adjacent groups give the kernel of $q:\Lambda'\to\Lambda'$
and $\pi_6\mathrm{Tw}(A_\infty)=0$. Thus the forgetful map gives
\[
 \mathsf F^{\mathrm{st}}_{k,l}(X_q)
       \ \xrightarrow{\ \cong\ }\ \mathrm{Br}_A(X_q)[k^\infty]
       \ \cong\ \mathbb Z/q .
\]
Here $H^3(X_q;\mathbb Z)=0$, so no degree-three Brauer invariant
detects the generator. The two examples exhibit the kernel and the
image, respectively, of the same forgetful map. Compare the
Moore-space calculation in \cite[Example~3.5]{DPBrauer}.
\end{example}

\begin{corollary}[Prescribed sizes on the Moore space]\label{cor:moore-sizes}
For the space $X_q$ of Example~\ref{ex:relative-moore}, a class
$\tau\in [X_q,\mathrm{Tw}(A)]\cong\mathbb Z/q$ admits a LAB of
type $(m,l)$, $(m,l)=1$, exactly when its order divides $m$.
In particular the least such size for a generator is $q$.
\end{corollary}
\begin{proof}
Use the projection $Y_m//U_R\to B U_R$ of Theorem~\ref{thm:rank-section}
and the identification of $B U_R$ with $\mathrm{Tw}(A)$.
For any of these simply connected targets $Z$, the bottom-cell
inclusion gives an isomorphism $[X_q,Z]\cong(\pi_4Z)[q]$ when
$\pi_5Z=0$, by the cofiber sequence for multiplication by $q$ on
$S^4$. These sets carry abelian group structures functorial in $Z$,
since $X_q$ is a double suspension. Apply this calculation to
$Y_m$ and $Y_m//U_R$. We obtain $[X_q,Y_m]=0$ and
$[X_q,Y_m//U_R]\cong(\Lambda/m\Lambda)[q]$.
The exact sequence for the projection thus injects the latter group
into $[X_q,B U_R]\cong\mathbb Z/q$. Its image has order
$\gcd(m,q)$ and is therefore precisely the subgroup killed by $m$.
\end{proof}

\begin{corollary}[Obstruction to a natural section]\label{cor:no-natural-splitting}
For each $k>1$ the forgetful map in~\eqref{eq:brauer-short-exact}
admits no section natural on connected finite CW complexes, even
as a map of set-valued functors.
\end{corollary}
\begin{proof}
Use $X_q$ of Example~\ref{ex:relative-moore} with $q>1$ dividing
a power of $k$, and let $i:S^4\hookrightarrow X_q$ include its
bottom cell. Forgetting the trivialization on $X_q$ is a bijection
onto $\mathbb Z/q$. On $S^4$ all Brauer restrictions are zero
since $\pi_4\mathrm{Tw}(A)=0$. On the LAB side, the cofiber
sequence defining $X_q$ gives
\[
  0=\pi_5\mathbf Q_{k,l}\longrightarrow
  \mathsf F^{\mathrm{st}}_{k,l}(X_q)
       \xrightarrow{i^*} C_k \xrightarrow{\ \times q\ } C_k.
\]
Thus the $q$ unique lifts on $X_q$ restrict to $q$ distinct
elements of $C_k$ on $S^4$. Naturality of a section would send
all their restrictions to its single value at the zero Brauer
class, a contradiction. On the full subcategory with $b=0$ in
Theorem~\ref{cor:brauer-forget}, the forgetful map is an
isomorphism and its inverse is natural.
\end{proof}

\paragraph{Remaining questions.}
For a finite CW complex $X$, every raw or completed coherent class has
a normalized strict representative, unique modulo morphisms, refinements
and concordances. A question about families remains. Can one construct a
moduli space directly from strict diagrams on varying covers so that,
for each $D$-bundle $E\to X$, the homotopy fiber of rectification over
$E$ is weakly contractible? A point of this fiber includes a chosen
identification with $E$. Forgetting that identification allows the gauge
symmetries of Section~\ref{subsec:symmetries} to act. For embedded
diagrams over finite polyhedra, Theorem~\ref{thm:relative-existence}
already extends collared boundary data over $X\times D^r$ after
refinement. The further issue is to model families of refinements and
their higher compatibilities and to compare their homotopy fibers with
this geometric extension problem. Example~\ref{ex:raw-refinement}
exhibits the role of the choice of cover even for paths.

The role of concordance in Corollary~\ref{cor:ambient-unique} and
Proposition~\ref{prop:normalization} leaves the rigid equivalence
relation as a separate problem. Extending the geometric classification
to arbitrary CW complexes would also require controlling the finite-stage
and finite-cover choices made over compact sets.

\section{Permutation LABs and complex linearization}\label{sec:perm-linearization}

The rank-section description of Theorem~\ref{thm:rank-section} permits a
comparison with finite sets. A decomposition of a finite set into copies
of a smaller set gives an inclusion of matrix algebras on their complex
spans. Which matrix LAB classes can arise in this way? We compare the
classifying spaces first, then calculate the image on $S^2$ and the effect
on two underlying twists. We then determine the least object size needed
to represent a generator on the first odd-primary sphere.

Fix $k\ge1$ and $l\ge2$ with $(k,l)=1$. Retain $R=KU[1/l]$ and
$\Lambda=\mathbb Z[1/l]$ and put $S=\mathbb S[1/l]$, where $\mathbb S$
is the sphere spectrum. Write $j:S\to R$ for the localized unit map.
For $E=S,R$, put
\[
  U_E=\mathrm{GL}_1(E)_1,\qquad Y_E(k)=(\Omega^\infty E)_k,\qquad
  Z_E(k)=Y_E(k)//U_E,\qquad q_E:Z_E(k)\longrightarrow BU_E.
\]
Double slash denotes the homotopy quotient under multiplication.
A map to $BU_E$ classifies an $E$-line with structure group reduced to the identity
component of its units. This reduction fixes the coordinate of fiber degrees in
$\Lambda$. A lift to $Z_E(k)$ adds a section of degree $k$. All bars, section spaces and
fibers below can be taken in derived models.

\begin{definition}[Permutation category and complex linearization]\label{def:perm-lab}
The permutation category $\mathcal P_{k,l}$ has discrete object and morphism spaces.
Its objects are the standard sets of sizes $kl^n$, one for each $n\ge0$. An arrow $X\to Y$ is a bijection $u:X\times D\cong Y$ with $|D|$ a
power of $l$, modulo isomorphism of the multiplicity set $D$. Composition uses Cartesian
product of multiplicities. The matrix category $\C_{k,l}$ has corresponding objects
$\operatorname{End}(\mathbb C[X])$ and all unital $*$-embeddings, with the point-norm
topology. Let $P_u$ denote the unitary induced by $u$. The functor $L$ on arrows is
\[
 L[u](a)=P_u(a\otimes I_{\mathbb C[D]})P_u^{-1}.
\]
Relabeling $D$ acts on its factor alone and leaves this expression unchanged.
For successive arrows represented by $X\times D\cong Y$ and
$Y\times D'\cong Z$, the composite is represented by
$X\times(D\times D')\cong Z$. The canonical identification
$\mathbb C[D\times D']\cong\mathbb C[D]\otimes\mathbb C[D']$ proves
that $L$ preserves composition. For a space $X$, a permutation LAB class means a
homotopy class of maps $X\to B\mathcal P_{k,l}$.
\end{definition}

\begin{theorem}[Complex linearization over the twist spaces]\label{thm:perm-linearization}
There is a homotopy-commutative diagram
\[
\begin{array}{ccc}
 B\mathcal P_{k,l}&\xrightarrow{\ BL\ }&B\C_{k,l}\\[-2pt]
 \big\downarrow{\scriptstyle\simeq}&&\big\downarrow{\scriptstyle\simeq}\\
 Z_S(k)&\xrightarrow{\ Z_j\ }&Z_R(k)\\
 {\scriptstyle q_S}\big\downarrow&&\big\downarrow{\scriptstyle q_R}\\
 BU_S&\xrightarrow{\ B\mathrm{GL}_1(j)_1\ }&BU_R .
\end{array}
\]
The right equivalence uses the scalar section given by the unit of a matrix algebra. For
a finite CW base, the diagram sends a pair $(\mathcal L,s)$ to $(R\otimes_S\mathcal
L,1\otimes s)$. For a matrix LAB with rectified algebra bundle $\mathfrak E$, the
projection $q_R$ classifies the orientable twist of $\mathfrak E\otimes\mathbb K$.
\end{theorem}

The additive group completion of finite sets gives the sphere spectrum
by the Barratt--Priddy--Quillen theorem \cite{Segal74,GMMO}. Tornehave's
multiplicative localization identifies spaces of units using finite sets
of specified cardinalities \cite[Chapter~VII, after Corollary~5.4]{MayTornehave}.
At fixed $k$ one must also keep the space on which these units act.
The proof below identifies that action from the arrows of the category.
The resulting equivalence applies to the ordinary nerve itself. The
group-completion theorem enters the proof for its monoid of multiplicities.

\subsection{A common model for the two nerves}
Use $G_d=\Sigma_d$ in the set case and $G_d=U(d)$ in the matrix case. Identify
successive products by ordered bases. Each space of arrows of multiplicity $t$ is
$G_{dt}/(1_d\otimes G_t)$, where $\otimes$ denotes Cartesian product of permutations or
tensor product of unitaries. Define a monoid and a right module
\[
 P_G=\coprod_{s\ge0}BG_{l^s},\qquad M_G=\coprod_{r\ge0}BG_{kl^r}.
\]
Their multiplication and action use the product of factors. The maps of these spaces
induced by permutation matrices preserve the products. Associativity can be made strict
with ordered models.

\begin{lemma}\label{lem:perm-bar}
The nerve in each case has the homotopy type of the two-sided bar $B(M_G,P_G,*)$. Under
these identifications $BL$ is induced by $\Sigma_d\to U(d)$ on every factor.
\end{lemma}
\begin{proof}
Resolve a string of $p$ arrows by the groupoid of changes of coordinates at its $p+1$
vertices. An element of $G_d$ acts by the corresponding algebra automorphism in the
matrix case. Keep the actual unitaries in the resolution, including scalars.
Equivalently, take the bisimplicial nerve of grids with $p$ horizontal arrows and $q$
successive vertical coordinate changes, requiring commutativity after the coordinate
changes act. Regard the original nerve as constant in the vertical direction.
Inserting identity coordinate changes gives a map $s$ from this nerve to the grid:
it commutes with the faces and degeneracies in both directions. At fixed vertical
degree $q$, top-row evaluation $e_q$ is a horizontal functor with $e_qs_q=\id$.
There is a continuous natural isomorphism $s_qe_q\Rightarrow\id$. At each vertex
its top arrow is the identity, and its arrow in row $j$ is the composite of the
first $j$ vertical changes, acting by conjugation in the matrix case. Thus $s_q$
induces a homotopy equivalence on horizontal classifying spaces for every $q$.
It is the bisimplicial map $s$ that gives the comparison on realization.

For fixed sizes $a,at_1,\ldots,at_1\cdots t_p$, the coordinate group acts transitively
on the space of horizontal strings: one can make every arrow standard in succession. The
stabilizer of the standard string is
\[
                  G_a\times G_{t_1}\times\cdots\times G_{t_p}.
\]
Indeed the successive coordinate changes have the form $h_i=h_{i-1}\otimes b_i$, with
$b_i\in G_{t_i}$ unique. This is also the centralizer calculation for a standard unital
matrix embedding. For a Lie group orbit the quotient here has local sections;
equivalently the Borel construction of the orbit $H/K$ is $EH/K$. Thus inclusion of this
stabilizer induces a weak equivalence on vertical classifying spaces. Omitting the first
vertex merges $h_0$ with $b_1$, an interior face multiplies adjacent multiplicities, and
the last face drops the final multiplicity. Insertion of an identity uses the
distinguished identity of multiplicity one. The standard strings therefore give
precisely the bar above, including the $U(1)$ factor in the unitary case. These two
comparisons of the grid are natural in $\Sigma_d\to U(d)$. Degreewise singular complexes
and diagonals justify realization of the comparisons. The ordinary nerves have the same
weak type here: the degeneracies of the discrete nerve are cofibrations and those of the
matrix nerve insert an identity in a compact smooth morphism manifold on each index
component.
\end{proof}

\subsection{Localization of the bar}
The additive class map $\coprod_d BG_d\to\Omega^\infty E$ is the BPQ map for sets and
the usual vector-bundle class map for unitaries, followed by coefficient localization.
Both respect products and the maps compare under $j$. Indeed $X\mapsto\mathbb C[X]$
takes disjoint union to direct sum and Cartesian product to tensor product, with the
canonical bijections and unitary identifications satisfying the distributivity and
permutation identities. On the resulting ring spectra its map is the unit from the
sphere to connective complex $K$-theory, followed by Bott inversion and localization
\cite{Segal74,GMMO}. Write
\[
  D_E=\coprod_{z\in\mathbb Z}(\Omega^\infty E)_{l^z},
  \qquad W_E=\coprod_{z\in\mathbb Z}(\Omega^\infty E)_{kl^z}.
\]
The space $D_E$ is a grouplike multiplicative subspace of the units. We use $D_E$ acting
on $W_E$, retaining the action of its identity component on $Y_E(k)$.

\begin{proposition}\label{prop:perm-localization}
The following comparisons are natural under complex linearization:
\[
 B(M_G,P_G,*)\longrightarrow B(W_E,D_E,*)
                         \longleftarrow Z_E(k).
\]
Both arrows are weak equivalences.
\end{proposition}
\begin{proof}
The right arrow includes a component and its stabilizer. The action on component labels
is transitive, with stabilizer $U_E$ at $k$. Equivalently $Y_E(k)\times_{U_E}^{h}D_E\to
W_E$ is an equivalence on each component, since multiplication by $l^z$ is an
equivalence. Taking homotopy orbits gives the claim.

\smallskip\noindent\emph{The telescope.} For the other arrow write $P=P_G$, $M=M_G$ and
let $M_\infty=\operatorname{Tel}(M\to M\to\cdots)$ under multiplication on the left by a
fixed $l$-element set or trivial $l$-space. This operation commutes with the right
$P$-action. On $B(M,P,*)$ it induces a map homotopic to the identity: on the category of
Lemma~\ref{lem:perm-bar} the standard arrows $X\to[l]\times X$ (or $a\mapsto I_l\otimes
a$) give a natural transformation. Composition on either side differs only by a
reordering of the multiplicity factors. The same transformation acts on the grid
resolution, commuting with coordinate changes. Taking the telescope inside the bar
therefore preserves its weak type.

Let $c:M\to W_E$ be the additive class map, where $W_E$ is a right $P$-module
via $P\to D_E$. Write $L_l$ for left multiplication by the chosen object of size
$l$, or by its class in $E$. Compatibility of $c$ with products gives a map
of sequential diagrams of right $P$-modules
\[
\begin{array}{ccccccc}
 M   &\xrightarrow{L_l}& M   &\xrightarrow{L_l}& M   &\longrightarrow&\cdots\\
 {\scriptstyle c}\big\downarrow &&{\scriptstyle c}\big\downarrow
     &&{\scriptstyle c}\big\downarrow\\
 W_E &\xrightarrow{L_l}& W_E &\xrightarrow{L_l}& W_E &\longrightarrow&\cdots .
\end{array}
\]
Every lower transition is an equivalence since $l$ is a unit in $\pi_0 E$.
Hence the stage-zero inclusion $i_0:W_E\to\operatorname{Tel}(W_E,L_l)$
is an equivalence of right $P$-modules. In the homotopy category of these
modules the required comparison is
\[
 M_\infty \xrightarrow{\operatorname{Tel}(c)}
 \operatorname{Tel}(W_E,L_l)\xrightarrow{\ i_0^{-1}\ } W_E .
\]
For the lower stage inclusions $i_h$ the telescope cylinders give
$i_h\simeq i_{h+1}L_l$. Thus at integer stage $h\ge0$ the comparison sends
$x$ to $l^{-h}c(x)$ up to homotopy, with compatibility on successive stages
supplied by the cylinders. On a given component the pre-localization
comparison is the telescope of the usual maps
$BG_d\to(\Omega^\infty E_0)_d$, where $d$ increases by a factor $l$ at every step and
$E_0=\mathbb S$ or $KU$. The connecting maps on the latter spaces are $l$-fold addition.
Their telescope has the claimed localized component, by the homotopy group calculation
for a telescope. For $U(d)$ the maps to stable $BU$ become arbitrarily highly connected.
For $\Sigma_d$ the maps are integral homology isomorphisms in a range tending to
infinity by symmetric-group stability and BPQ. Hence $M_\infty\to W_E$ is a homology
equivalence on each component. The same argument for $M=P$, with $k=1$, and the group
completion theorem \cite{McDuffSegal}, show that $P^{\rm gp}\to D_E$ is a weak
equivalence. In detail, this map is an isomorphism on component groups and on the
integral homology of each component. Its source and target are grouplike
$E_\infty$-spaces, so their components are simple and the homology isomorphisms are weak
equivalences. In particular $BP\to BD_E$ is an equivalence. Since $P$ acts on $W_E$ by
equivalences, $B(W_E,P,*)\to B(W_E,D_E,*)$ is an equivalence: each bar projects with
homotopy fiber $W_E$ to the corresponding classifying space.

\smallskip\noindent\emph{Homology of the sign covers.} To recover the
ordinary homotopy type on the permutation side, we need homology with
the local coefficients coming from the fundamental group. If $l$ is odd, the map
$M_\infty\to W_S$ also induces a homology isomorphism on the double covers pulled back
from the universal covers of the components of $W_S$. On a sufficiently large stage the
covering on $B\Sigma_d$ is $BA_d$ for the alternating group $A_d$. Indeed the comparison
on fundamental groups is sign, and multiplication by odd $l$ preserves it. At each stage
$BA_d$ maps to the universal cover of the corresponding component of $Q
S^0=\Omega^\infty\mathbb S$ by a homology isomorphism in degrees $i\le (d-3)/3$
\cite[Theorem 5.2]{RWWahl}. Those universal covers telescope to the target cover, by
their homotopy groups. For even $l$ the components of $W_S$ are simply connected. It
follows that the bar comparison with $B(W_S,P,*)$ is an integral homology equivalence
and also a homology equivalence with coefficients pulled back from a double cover of its
target. Indeed on a term $W_S\times P^p$ the map to the realized bar is homotopic to
projection to the $W_S$-vertex followed by inclusion. The cover on that term is
therefore pulled back from $W_S$; the claim follows by taking products and the skeletal
filtration of the bar.

\smallskip\noindent\emph{The fundamental group of the nerve.} In the groupoid completion
of the permutation category choose the standard arrows between successive sizes as a
tree. Every arrow is a standard arrow followed by a permutation. Transporting loops
along the tree gives the presentation
\[
 \pi_1 B\mathcal P_{k,l} =
  \big(\operatorname*{colim}_{r}\Sigma_{kl^r}\big)
  \big/\big\langle\!\big\langle\,1_a\otimes\Sigma_t\,\big\rangle\!\big\rangle ,
\]
where the colimit uses diagonal replication, and $a=kl^r,\ t=l^s$. The relations are
sufficient: composition uses replication of permutations and equality of representatives
is precisely relabeling the multiplicity. At a stage of size $d=kl^r\ge5$ with $r\ge1$,
take a transposition of the $l$ multiplicity labels for the last arrow. Its image in
$\Sigma_d$ is nonidentity and its normal closure contains the alternating group. This
follows from simplicity of $A_d$ for $d\ge5$ and the trivial centralizer of $A_d$ in
$\Sigma_d$. Replication acts on signs by the $l$th power, and the sign of a relabeled
multiplicity transposition is $(-1)^a$. Thus the quotient is trivial except when $k$ is
even and $l$ is odd, in which case sign gives exactly $\mathbb Z/2$. In $Z_S(k)$ the
same group is the cokernel of $k:\pi_1S\to\pi_1S$, by the orbit sequence. A stage
transposition maps to its stable sign, so the comparison induces the stated isomorphism.
In the exceptional case the inclusion of the rank-$k$ component of $W_S$ induces
the quotient $\pi_1 S\to(\pi_1 S)/k(\pi_1 S)$ on fundamental groups. This is an
isomorphism because $\pi_1 S=\mathbb Z/2$ and $k$ is even. The target universal
cover therefore restricts to the universal double cover on this component,
and multiplication by $l$ transports the identification to the other components.
The preceding cover calculation now gives a homology isomorphism on universal covers. Whitehead's theorem now proves the weak equivalence. The unitary
comparison already follows degreewise from the stable connectivity calculation.
\end{proof}

\subsection{Compatibility with the scalar section}
\begin{proof}[Proof of Theorem~\ref{thm:perm-linearization}]
Lemma~\ref{lem:perm-bar} and Proposition~\ref{prop:perm-localization} give the upper
square, with $Z_j$ induced on both the action space and the acting units by $j$. The bar
projections are induced by $Y_E(k)\to *$, hence commute with this map on every bar
simplex. This gives the lower square with its specified homotopy. The associated line
construction sends a change of coordinates $v\in U_S$ to $j(v)\in U_R$ and the section
to its extension of scalars.

We verify that the right equivalence really uses the scalar section of matrix LABs. In
Morita coordinates for $A_X=\operatorname{End}(V_X)$, $V_X=\mathbb C[X]$, the regular
right module corresponds to $V_X$ and therefore has class $[V_X]$. For $u:X\times D\cong
Y$, extension of scalars on right modules satisfies the natural identity
\[
 (P\otimes_{A_X} A_Y)\otimes_{A_Y} V_Y
   \ \cong\ P\otimes_{A_X}V_Y
   \ \cong\ (P\otimes_{A_X}V_X)\otimes_{\mathbb C} V_D .
\]
This holds with arbitrary unitary identifications of the tensor factors. It is
equivariant in all changes of coordinates and commutes with iterated tensor products by
associativity. Consequently on the groupoid of any string, including its nonconstant
families, the line transition is multiplication by the class of its multiplicity bundle
and the section is the class of its defining vector bundle. These are precisely the
unitary bar maps used above. Selecting the degree-$k$ component divides the stage-$r$
coordinate by $l^r$, as in~\eqref{eq:scalar-pair}. Explicitly, an arrow of multiplicity
$t=l^b$ from stage $r$ to $r+b$ acts in normalized coordinates by
$[V_D]/t$, a unit of rank one. These normalized units multiply under
composition. By the bundle identification in Section~\ref{sec:brauer-relative},
forgetting the scalar section gives the twist of the stabilized rectified
algebra bundle.
\end{proof}

\subsection{The image on the two-sphere}
For a spectrum $E$ let $E/k$ denote the cofiber of multiplication by $k$. The following
identification concerns positive homotopy groups; it only requires a description of a
single loop space.

\begin{lemma}\label{lem:perm-moore}
For $i\ge1$ there are natural isomorphisms $\pi_i Z_E(k)\cong\pi_i(E/k)$ for $E=S,R$.
Under these identifications the linearization map is induced by $j/k:S/k\to R/k$.
\end{lemma}
\begin{proof}
The loop space of the homotopy quotient at the reference point is the homotopy fiber of
$U_E\to Y_E(k)$, $v\mapsto kv$. Translate the source by $-1$ and the target by $-k$ in
the additive infinite loop space. The resulting map is multiplication by $k$ on
$\Omega^\infty\tau_{\ge1}E$. After taking the fiber, the equivalence of based spaces is
with
$\Omega^\infty\Sigma^{-1}((\tau_{\ge1}E)/k)$. Taking $\pi_{i-1}$ for
$i\ge2$ gives $\pi_i((\tau_{\ge1}E)/k)$. The map to $\pi_i(E/k)$ is an isomorphism in
these degrees by the cofiber exact sequences. For $i=1$ the orbit sequence gives the
cokernel of $k$ on $\pi_1 E$. This agrees with $\pi_1(E/k)$ since $\pi_0 E=\Lambda$ is
torsion-free. All translations and the coefficient maps preserve the chosen basepoints
and commute with $j$.
\end{proof}

\begin{theorem}[Permutation lifts on the two-sphere]\label{thm:perm-sphere-image}
On $S^2$ the permutation classes have image exactly
\[
 \operatorname{im}\big(\pi_2 B\mathcal P_{k,l}\longrightarrow\pi_2 B\C_{k,l}\big)
       =(\mathbb Z/k)[2]\ \subseteq\ \mathbb Z/k .
\]
In particular the image is zero for odd $k$ and consists of $0,k/2$ for even $k$. These
are also precisely the unbased matrix classes admitting a permutation lift.
\end{theorem}
\begin{proof}
The unit from the sphere factors through $KO\to KU$. The map $\mathbb S\to KO$ is an
isomorphism on $\pi_1$ and $\pi_2$: the Hopf elements $\eta,\eta^2$ map to the
respective generators. Applying localization and the cofiber sequence for $k$ therefore
gives an isomorphism
\[
       \pi_2(S/k)\ \xrightarrow{\cong}\ \pi_2(KO[1/l]/k).
\]
Use the real--complex Bott cofiber sequence \cite[\S5.3, (5.3.2)]{Rognes}
\[
 \Sigma KO\xrightarrow{\eta}KO\xrightarrow{c}KU
       \xrightarrow{\partial}\Sigma^2 KO,\qquad
       \partial=(\Sigma^2 r)\,\beta^{-1},
\]
where $c$ is complexification, $r$ is realification and $\beta$ the complex Bott
equivalence. Localize and smash with $\mathbb S/k$. In degree two both the source and
target of $\partial_*$ identify with $\Lambda/k\Lambda$: for the target use $\pi_0
KO=\mathbb Z$ and $\pi_{-1}KO=0$. The map $\partial_*$ is multiplication by $2$ up to
sign, since $\beta^{-1}$ sends a Bott generator to complex rank one and realification
doubles its rank. Exactness gives the asserted image of $c_*$ and
Lemma~\ref{lem:perm-moore} identifies it with the categorical image. The target $Z_R(k)$
is simply connected, so passing from based to unbased classes does not change this
image.
\end{proof}

There is an elementary explanation for the upper bound. Complex
conjugation acts on matrix embeddings by
$f\mapsto(a\mapsto\overline{f(\bar a)})$. It fixes the image of
$L$ because every permutation matrix is real. On the target $\pi_2$ it negates the
generator, since it reverses the unitary clutching loop
$\operatorname{diag}(z,1,\ldots,1)$. A permutation image must therefore satisfy $x=-x$.
The Bott argument proves sufficiency as well. In particular the Hopf-line LAB of
Example~\ref{ex:morita-versus-lab} has a permutation representative at $k=2$, and for
every $k>2$ its generator lies outside the image. A representative here means a class in
the permutation nerve; it need not come from a single covering of $S^2$.

\subsection{Two underlying twists at the same size}
Every orientable complex twist on $S^2$ is trivial, since
$\pi_2 B U_R=\pi_1 R=0$. Theorem~\ref{thm:perm-sphere-image} therefore
measures information in the section. To see the effect on twists
themselves, fix $k=4$ and any odd $l\ge3$, for instance $l=3$. If a twist over a reduced
suspension $\Sigma Y$ is clutched by $u:Y\to U_E$, it admits a degree-$k$ section
exactly when $k(u-1)=0$ in $\widetilde E^0(Y)$ for connected based finite $Y$. Indeed
the two cones are contractible. The two constant degree-$k$ sections glue up to homotopy
exactly when $ku$ is homotopic to $k$ on the equator. Any section on a cone of this
degree is homotopic to the constant one.

\begin{example}[A twist detected on $\Sigma\mathbb{RP}^2$]\label{ex:perm-sign}
Let $Y=\mathbb{RP}^2$ and let $E_2\to Y$ be its canonical double cover. The BPQ class of
this two-point family gives $b\in S^0(Y)$ of degree two. Set $\alpha=b-2$ and
$u=1+\alpha$. The cofiber description $Y=S^1\cup_2 e^2$ gives an exact sequence
\[
 0\longrightarrow\operatorname{coker}(2:\pi_2 S\to\pi_2 S)
 \longrightarrow\widetilde S^0(Y)
 \longrightarrow\ker(2:\pi_1 S\to\pi_1 S)\longrightarrow0.
\]
Both outer terms are $\mathbb Z/2$, so $4\alpha=0$. Clutching by $u$ therefore gives a
spherical twist over $\Sigma Y$ with a degree-four lift.

Let $\lambda$ be the tautological real line over $Y$. The complex permutation bundle of
the cover is $\mathbb C\oplus\lambda_{\mathbb C}$, by the decomposition of the two-point
representation into its trivial and sign lines. Hence $j(u)=[\lambda_{\mathbb C}]$. This
degree-one unit is nontrivial: $c_1(\lambda_{\mathbb C})$ is the integral Bockstein of
$w_1(\lambda)$ and generates $H^2(Y;\mathbb Z)=\mathbb Z/2$. On this two-dimensional
complex the first Chern class identifies reduced complex $K$-theory with $H^2$;
localization at odd $l$ preserves the generator. Suspension adjunction now shows that
the clutched $R$-line is nontrivial. Thus linearization detects the \emph{underlying
twist}, and every choice of lift above it has nonzero image as a matrix LAB. Its twist
is the nonzero class in $H^3(\Sigma\mathbb{RP}^2;\Lambda)\cong\mathbb Z/2$: the first
nonzero homotopy group of $BU_R$ is $\pi_3=\Lambda$, and the base has dimension three.
\end{example}

\begin{example}[A twist lost on $S^3$]\label{ex:perm-eta}
Take $Y=S^2$ and let $u$ be the degree-one unit whose reduced part is the stable class
$\eta^2\in\pi_2 S=\mathbb Z/2$. Here $\eta\in\pi_1\mathbb S$ is the stable Hopf
element. Clutching by $u=1+\eta^2$ gives a nontrivial spherical twist on $S^3$. Since $4\eta^2=0$, a
degree-four section exists. Under $j$, the equatorial class $\eta^2$ maps to zero in the
torsion-free group $\pi_2R=\Lambda$. In fact the \emph{entire pair} becomes trivial for
every choice of degree-four lift, since
\[
 \pi_3 Z_R(4)\cong\pi_3(R/4)=0
\]
by $\pi_3R=0$ and injectivity of $4$ on $\pi_2R$. The same example works for every even $k$ and odd $l$ coprime to it.
\end{example}

\subsection{The least object size on the first odd-primary sphere}
\label{subsec:perm-first-size}

Let $p$ be an odd prime with $p\nmid l$. Write $\mathcal P^{\le N}_{p,l}$ for
the full subcategory on $[p],[pl],\ldots,[pl^N]$. A matrix class on a sphere
has a permutation representative of \emph{height at most $N$} if it lies in
the image of $B\mathcal P^{\le N}_{p,l}\to B\C_{p,l}$ on that sphere.
Thus the largest permitted object has $pl^N$ elements. This measures the
sizes in a representative while keeping the prescribed section degree $p$ fixed.

\begin{theorem}[Sharp object size]\label{thm:perm-first-size}
Put $d=2p-3$. The first positive even dimension in which a nonzero matrix
LAB of size $p$ has a permutation representative is $d+1=2p-2$. In that
dimension,
\[
 \operatorname{im}\!\left(\pi_{d+1}B\mathcal P^{\le N}_{p,l}
              \longrightarrow\pi_{d+1}B\C_{p,l}\right)
 =\begin{cases}
  0,& l^N<p,\\
  \mathbb Z/p,& l^N\ge p.
 \end{cases}
\]
For $N_0=\min\{N:l^N\ge p\}$, the two objects $[p]$ and $[pl^{N_0}]$
alone realize the entire image. In particular $pl^{N_0}$ is the least
possible largest object for each nonzero class. The same assertion holds
for free homotopy classes since $B\C_{p,l}$ is simply connected.
\end{theorem}

\begin{proof}
\emph{Detection in the first stem.}
Use $p$-local homology in this proof. The first positive $p$-primary stable
stem is $\pi_d\mathbb S_{(p)}=\mathbb Z/p$, and
$\pi_{d+1}\mathbb S_{(p)}=0$ \cite[Theorem~2.7]{Yamaguchi}.
The orbit sequence for $q_S$ has orbit map multiplication by $p$ after
translation of components, as in Lemma~\ref{lem:perm-moore}. It follows that
\[
 (q_S)_*: \pi_{d+1} Z_S(p)\xrightarrow{\ \cong\ }
                  (\pi_{d+1} BU_S)_{(p)}=\mathbb Z/p ,
\]
where the displayed map includes localization of the target. Indeed the
source is $\pi_{d+1}(S/p)$, which is already $p$-primary, and this map
identifies with the Moore boundary at $p$. Moreover
$BL$ sends a generator to a generator. To see this, write
$M(p)=\mathbb S/p$ and $i:\mathbb S\to M(p)$ for its bottom cell. The Adams
map $v:\Sigma^{d+1}M(p)\to M(p)$ induces an isomorphism in periodic
complex $K$-homology \cite[\S12]{AdamsJIV}. Hence
$v\circ\Sigma^{d+1}i$ has nonzero $KU$-Hurewicz image, since $i$ represents
$1\bmod p$ in $KU_0(M(p))$. By Lemma~\ref{lem:perm-moore} this Hurewicz
map is precisely linearization on the sphere. Inverting $l$ does not
change either mod-$p$ spectrum. The same lemma and the vanishing of the
earlier stable stems at $p$ give zero permutation image in every positive
even dimension below $d+1$.

\smallskip\noindent
\emph{The cycle carried by a multiplicity set.}
Let $t$ be a power of $l$ with $t\ge p$. The proof of
Lemma~\ref{lem:perm-bar} restricted to the two objects $[p],[pt]$
identifies their classifying space $X_t$ with the double mapping cylinder
of
\begin{equation}\label{eq:perm-two-level}
 B\Sigma_p\ \xleftarrow{\ \mathrm{pr}_1\ } B(\Sigma_p\times\Sigma_t)
 \xrightarrow{\ B\rho\ } B\Sigma_{pt},
 \qquad \rho(\sigma,\tau)=\sigma\times\tau .
\end{equation}
There is just one strict size increase. The vertex spaces record
permutations of each object. In the edge space a change of coordinates
also permutes the $t$ multiplicity labels.

For $C_p\subset\Sigma_t$ generated by one $p$-cycle, let
$z\in H_d(B\Sigma_t;\mathbb Z_{(p)})$ be the image of a generator of
$H_d(BC_p;\mathbb Z)=\mathbb Z/p$. This image is nonzero.
The permutation representation restricts to the regular representation of
$C_p$ plus trivial summands. With $u=c_1(\chi)$ for a generating
character, its $(p-1)$st Chern class restricts to
$(p-1)!u^{p-1}=-u^{p-1}\ne0$. For a finite group the positive integral
homology groups are finite, so the natural universal coefficient
isomorphism in degree $d+1$ identifies this restriction with a character
on $H_d(BC_p;\mathbb Z)$. Nakaoka's injectivity theorem for adjoining a
fixed point \cite{Nakaoka} and BPQ \cite{Segal74,GMMO} show that $z$
survives in the stable symmetric group. They also show that all finite
symmetric groups have zero $p$-local homology in degrees $0<i<d$:
the $p$-localized degree-zero component of $\Omega^\infty\mathbb S$ is
$(d-1)$-connected, with first homology group $\mathbb Z/p$ in degree $d$.

Include $z$ in the second factor of the edge in~\eqref{eq:perm-two-level}.
The left endpoint kills it. At the right endpoint the map factors through
the diagonal $\Sigma_t\to(\Sigma_t)^p$ followed by the block inclusion.
In degree $d$ the K\"unneth formula has only the $p$ single-factor
summands, by the preceding vanishing. The diagonal sends $z$ to
$(z,\ldots,z)$ and the block inclusions are conjugate. Its image is
therefore $p$ times the image of $z$ under a single block inclusion, hence zero. Mayer--Vietoris
gives $c\in H_{d+1}(X_t;\mathbb Z_{(p)})$ with edge boundary $z$.

\smallskip\noindent
\emph{Spherical realization and its image.}
On $X_t$, the projection to $BU_S$ is constant on both vertices of the
bar model. Its edge loop is the composite of the second projection with
$\mu_t:B\Sigma_t\to U_S$, the BPQ class of the multiplicity set divided
by $t$. This follows by forgetting the first factor in the bar and using
the normalized coordinates of Theorem~\ref{thm:perm-linearization}.
Collapsing the two vertex spaces of the double mapping cylinder gives the
suspension description of the Mayer--Vietoris boundary. Consequently the
image of $c$ in $H_{d+1}(BU_S;\mathbb Z_{(p)})$ is, up to sign, the
homology suspension of $(\mu_t)_*z$. Normalization by $t$ is an equivalence
on the localized components, so $(\mu_t)_*z\ne0$. Homology suspension
is an isomorphism here: the Hurewicz isomorphisms in the first nonzero
$p$-local degree for $U_S$ and $BU_S$ identify it with
$\pi_d U_S\cong\pi_{d+1}BU_S$.

The space $X_t$ is simply connected. By van Kampen its fundamental group
is the quotient of $\Sigma_{pt}$ by the normal closure of
$1_p\times\Sigma_t$. A transposition in $\Sigma_t$ gives an odd permutation
of $[pt]$, whose normal closure is all of $\Sigma_{pt}$ since $pt\ge5$.
Mayer--Vietoris and the lower homology vanishing above show that
$(X_t)_{(p)}$ is $(d-1)$-connected. Hurewicz is surjective in degree
$d+1$ for such a space ($d\ge3$): in a CW model with no cells in
dimensions $1,\ldots,d-1$, the $d$-skeleton is a wedge of $d$-spheres.
The cellular degrees detect the attaching maps in $\pi_d$ of this wedge,
so a cellular $(d+1)$-cycle lifts by the relative homotopy sequence. Thus $c$ has a spherical representative
after localization. As $X_t$ is simply connected,
$\pi_{d+1}((X_t)_{(p)})=\pi_{d+1}(X_t)\otimes\mathbb Z_{(p)}$.
Clearing a denominator prime to $p$ gives an integral sphere map with
nonzero image in $\mathbb Z/p$. Multiples realize every matrix class.

\smallskip\noindent
\emph{Sharpness with all intermediate objects allowed.}
On the truncated bar for stages $0,\ldots,N$, the map to $BU_S$ factors
through the bar of the multiplicities alone. Over a chain
$i_0\le\cdots\le i_r$ its space is
$\prod_{a=1}^r B\Sigma_{l^{i_a-i_{a-1}}}$: the projection drops the
initial coordinate group, and composition multiplies successive
multiplicities. When $l^N<p$, each factor is $p$-locally acyclic by
transfer. Projection to the stage chain therefore induces a homology
isomorphism on each bar level. The realization spectral sequence gives
the $p$-local homology of the simplex on $0,\ldots,N$, so every sphere
maps trivially to $H_{d+1}(BU_S;\mathbb Z_{(p)})$. By the first-degree
Hurewicz isomorphism and the calculation of $(q_S)_*$ above, its
permutation class in the full category is zero. Taking $t=l^{N_0}$ in
the upper bound completes the proof.
\end{proof}

\subsection{An integral cycle on twelve letters}
\label{subsec:perm-twelve}
For $p=3$ and $l=2$, Theorem~\ref{thm:perm-first-size} gives $N_0=2$.
The two objects $[3]$ and $[12]$ already suffice on $S^4$, whereas
the full subcategory on $[3],[6]$ gives zero. Here the spherical class
can also be specified by an integral cycle in the double mapping cylinder
$X_4$.

\begin{proposition}\label{prop:perm-twelve-cycle}
Let $g=(1\,2\,3)\in\Sigma_4$ and, for any element $h$ of order three, set
\[
 z(h)=[h\mid h\mid h]+[h\mid h^2\mid h]
\]
in normalized integral bar chains. Take any $4$-chain $C$ in $B\Sigma_{12}$
with $\partial C=2z(G)$, where
$G=(1\,2\,3)(5\,6\,7)(9\,10\,11)$. Such a chain exists. Capping the
edge cylinder on $2z(g)$ by $C$, with matching boundary orientations,
gives an integral $4$-cycle in $X_4$. Its homology class is represented by
a map $S^4\to X_4$, and each such map linearizes to a generator of
$\pi_4 B\C_{3,2}=\mathbb Z/3$.
\end{proposition}
\begin{proof}
The periodic resolution of $C_3$ sends its degree-three generator to
$z(g)$ in the bar resolution after passage to coinvariants. In
\eqref{eq:perm-two-level} this edge cycle goes to zero at the left endpoint
and to $z(G)$ at the right. The replication argument above kills
$[z(G)]$ at $3$. Its integral order divides $3$ because it comes from
$H_3(BC_3;\mathbb Z)$. Hence it vanishes integrally, which gives the
claimed filling. Every choice of filling has Mayer--Vietoris boundary
$\pm2[z(g)]$, a generator at $3$.

To realize this integral class by a sphere it suffices to check that
$X_4$ is $2$-connected. Simple connectivity was proved above. On $H_1$
the map from the edge to its two ends is an isomorphism: the first end
records $\operatorname{sgn}(\sigma)$ and the second records
$\operatorname{sgn}(\tau)$, since the two exponents in the product sign
are $4$ and $3$. On $H_2$ the target is
$H_2(B\Sigma_3;\mathbb Z)\oplus H_2(B\Sigma_{12};\mathbb Z)
=0\oplus\mathbb Z/2$ by the Schur multiplier calculation
\cite[\S3]{LassueurMazza}. The commuting transpositions $(1\,2),(3\,4)$
in the multiplicity group define a torus in $B\Sigma_4$. After replication
its real permutation bundle has three sign lines in each torus direction
and the other summands are trivial. For the two degree-one generators
$x,y$ modulo $2$, its total Stiefel--Whitney class is
$(1+x)^3(1+y)^3$; its $w_2$ evaluates to $1$ on the torus. Thus the
integral fundamental class of this torus hits the generator of the target
$H_2$. Mayer--Vietoris gives $H_2(X_4;\mathbb Z)=0$ and Hurewicz gives
$\pi_2X_4=0$. Now $\pi_4X_4\to H_4(X_4;\mathbb Z)$ is surjective.
The edge boundary computed above detects the linearization by the proof
of Theorem~\ref{thm:perm-first-size}.
\end{proof}

The accompanying file \texttt{permutation\_lab\_s4\_certificate.py}
specifies one such $C$ as a signed sum of $152$ distinct nondegenerate
bar simplices (sum of absolute coefficients $188$). It checks the integer
boundary using right-to-left composition of permutations and prints the
individual simplices with \texttt{--terms}. This specifies the cycle
explicitly; Proposition~\ref{prop:perm-twelve-cycle} supplies its
spherical realization.

\end{document}